\ifdefined\pdfoutput\pdfoutput=1\fi
\documentclass[11pt]{article}
\usepackage[T1]{fontenc}
\usepackage[utf8]{inputenc}
\usepackage{lmodern}
\usepackage{amsmath,amssymb,amsthm}
\usepackage[margin=25mm]{geometry}
\usepackage[numbers,sort&compress]{natbib}
\usepackage{graphicx}
\usepackage{subfig}
\usepackage{booktabs}
\usepackage{multirow}
\usepackage[hidelinks]{hyperref}

\hypersetup{pdftitle={Variably saturated flow on polygonal and cut meshes using a mass-conservative hybrid virtual element}, pdfauthor={Francesco S. Liguori}}

\newcommand{\bq}{\mathbf q}
\newcommand{\bn}{\mathbf n}
\newcommand{\bx}{\mathbf x}
\newcommand{\bX}{\mathbf X}
\newcommand{\br}{\mathbf r}
\newcommand{\bh}{\mathbf h}
\newcommand{\Nq}{\mathbf N_q}
\newcommand{\Nh}{\mathbf N_h}
\newcommand{\Oe}{\Omega_e}
\newcommand{\dE}{\partial\Oe}
\newcommand{\Ae}{A_e}

\title{Variably saturated flow on polygonal and cut meshes using a mass-conservative hybrid virtual element}
\author{Francesco S. Liguori\\[2mm]
{\normalsize DIMES, University of Calabria, 87030 Rende (Cosenza), Italy}\\
{\normalsize\texttt{francesco.liguori@unical.it}}}
\date{}

\begin{document}

\maketitle

\begin{abstract}
A hybrid virtual element formulation (HVEM) is presented for transient flow in variably saturated soils governed by Richards' equation. The nodal heads are the only global degrees of freedom, and a Darcy flux is assumed in each polygonal element and condensed. The flux combines a divergence-free part with a particular field of linear divergence that absorbs storage and sources, so that the water balance of every element holds identically, at every iterate and for any time step, and that of the whole domain at convergence. To guarantee rank-sufficiency with no stabilisation parameter, the divergence-free fields are harmonic polynomials ordered by their behaviour under rotation, giving a closed-form rule for the number of flux parameters, derived on regular polygons and verified numerically on distorted ones and on elements with hanging nodes. Darcy's law is written in each element for the Kirchhoff potential, which removes the conductivity from the diffusion while the head stays continuous across material interfaces. A Kirchhoff mean of the conductivity makes hydrostatic equilibrium and gravity drainage exact on any mesh. The storage is lumped on the nodes with weights supplied by the assumed flux. All the element operators reduce to boundary integrals, so that the nonlinear iterations need only nodal values. The formulation is assessed on solutions in closed form, on published benchmarks and on rainfall over homogeneous and layered slopes, where independent solutions provide a further check. It converges at second order in the head on polygonal and locally refined meshes, also where gravity dominates, and closes the water budget in every test. On a layered slope, cut grids refined in a band along the surface reach the accuracy of a uniform grid twice as fine with about $40\%$ of the nodes. Its structure, shared with the HVEM for elastoplasticity, suits coupled analyses of slopes.

\medskip\noindent\textbf{Keywords:} Hybrid virtual element method; unsaturated flow; Richards' equation; polygonal meshes; mass conservation; stabilisation-free
\end{abstract}

\section{Introduction}\label{sec:intro}
Richards' equation \cite{Richards1931} is the standard model of water flow in variably saturated soils and underlies most analyses of infiltration, drainage and rainfall-induced slope instability. The equation is a mass balance in which Darcy's law is written with a conductivity and a water content that depend on the pressure head \cite{vanGenuchten1980,Mualem1976}. The speed of its solution matters, especially in early warning systems, where the flow model is coupled with rainfall forecasts and with the stability analysis of the slope \cite{CapparelliVersace2011}, and simplified methods are therefore often preferred \cite{ConteTroncone2012}.

Several features of Richards' equation make its numerical solution delicate \cite{Farthing2017Numerical,BergamaschiPutti1999}. The conductivity changes by orders of magnitude between wet and dry soil, so that wetting fronts are sharp and advance into a nearly impervious medium. The storage term depends on the unknown and vanishes at saturation, where the equation loses its parabolic character. Gravity drives the flow as much as the pressure gradient does, and in layered soils it acts across interfaces where the conductivity jumps. On a seepage face the boundary condition is not known in advance, and it switches between a prescribed head and an impervious boundary. Each of these features weighs on the choice of the discretisation.

Conservation of mass has guided much of the work on the discretisation, both in the storage term and in the flux. Celia et al. \cite{Celia1990A} showed that writing the storage as an increment of water content makes the scheme conservative, and that a lumped capacity prevents spurious oscillations ahead of a front. Finite volume schemes reached the same conclusions on heterogeneous soils and unstructured grids \cite{Cumming2011A,ManziniFerraris2004}. Mixed and mixed hybrid finite elements, usually in Raviart--Thomas spaces \cite{RaviartThomas1977}, approximate the flux independently and make it locally conservative \cite{BergamaschiPutti1999,Bause2004Computation}.

The complex geometry of slopes makes the use of polygonal meshes attractive. Their geometric freedom allows a slope to be described by a regular grid cut by the surface and by the interfaces, which produces polygons only where the geometry demands them, as in the interface-fitted meshes obtained by cutting a Cartesian grid \cite{Chen2017Interface,RangarajanSukumar2026}. Alternatively, the dual of a Delaunay triangulation follows an irregular domain with few elements. Local refinement needs no transition elements, since a hanging node is simply one more vertex \cite{vanHuyssteen2025}. The virtual element method (VEM) \cite{Beir2012Basic} provides the natural setting for such meshes. Primal head-based virtual elements have been developed for parabolic problems \cite{Vacca2015Virtual}, and applied to saturated--unsaturated seepage \cite{Sun2024Seepage} and to the soil in a coupled soil--root model \cite{Berrone2026}. Mixed virtual elements have been proposed for Darcy flow \cite{Veiga2015Mixed}, and a mixed virtual element for Richards' equation, based on a Kirchhoff transformation, has been developed and analysed by Adak, Manzini and Natarajan \cite{Adak2021The}. In the primal VEM the consistent part of the element is complemented by a stabilisation term, which is not unique and whose scaling is left to the analyst. Stabilisation-free formulations remove this term by choosing the projection so that the consistent part alone is rank-sufficient \cite{Berrone2023First,ChenSukumar2023,Sun2025}. The gradient is projected onto the gradients of harmonic polynomials \cite{Berrone2023First} or onto divergence-free fields \cite{BerroneDivFree2024}, the latter recently applied to the hydro-mechanical analysis of slopes \cite{BerroneSlope2026}.

Polygonal discretisations of Richards' equation are not confined to virtual elements. Hybrid finite volume and mimetic schemes, with degrees of freedom (dofs) on the elements and on the faces, have been unified in a single family \cite{Droniou2010Unified} and analysed for Richards' equation as gradient schemes \cite{Eymard2014Gradient}. Vertex-based schemes extend the control-volume finite element method to general meshes. The vertex approximate gradient scheme adds element dofs to the nodal ones \cite{Eymard2012VAG}, and vertex-based compatible discrete operators are built on a discrete Hodge operator \cite{BonelleErn2014}. Each of these schemes leaves some freedom in its construction, in the control volumes, in the weights between elements and vertices or in the discrete Hodge operator.

Hybrid formulations offer another route to stability. In the hybrid principle of Pian \cite{Pian1964} an independent field is assumed within each element and condensed, and the global dofs lie on the element boundaries. Hybrid-Trefftz elements \cite{Jirousek1978,JirousekWroblewski1996,Qin2000} choose the assumed field among exact solutions of the governing equation and connect the elements through a frame on their boundaries. The hybrid virtual element method (HVEM) \cite{Liguori2024A,Liguori2024B,Liguori2025Plates,Liguori2026Shells} belongs to this family. Its assumed field is drawn from a divergence-free polynomial basis rich enough to make the element stable without any stabilisation term, and the boundary trace keeps it conforming on polygons of any shape. The method was developed for elasticity, elastoplasticity, plates and shells, and proved accurate even on coarse meshes.

This paper extends the HVEM to transient flow in variably saturated soils. The nodal heads are the sole global dofs, and the head is assumed only on the boundary of each polygon, linear along each edge. Within each polygon a Darcy flux is assumed as an independent field and eliminated at the element level. The flux receives the greater care, since it carries the mass balance and produces the seepage forces that govern the stability of a slope. %

The storage rate depends on the unknown head. The assumed flux is therefore split into a divergence-free part, which carries the free parameters, and a particular part of linear divergence, which absorbs the mean and the first moments of storage and sources. The balance of every element then holds identically. It holds at every iterate of the nonlinear solver and for any time step, as a property of the flux rather than of the convergence of the iteration.

The stabilisation-free character of the HVEM rests on the choice of the flux basis. In mechanics a hierarchical selection was proposed \cite{Liguori2024A}, in which the stress basis is enriched by increasing the polynomial degree until the rank condition is met. Here an alternative, automatic selection is introduced. The divergence-free fields are harmonic polynomials ordered by their behaviour under rotation rather than by degree. The ordering yields a closed-form rule for the number of flux parameters on a polygon with any number of edges, and a simple correction for the collinear edges produced by hanging nodes. Since the rule is derived on regular polygons, its effectiveness on severely distorted ones is verified numerically, as in \cite{Liguori2025Plates,Liguori2026Shells}. For a unit conductivity the diffusion operator coincides with that of the harmonic stabilisation-free element \cite{Berrone2023First} and with that of the hybrid-Trefftz element for the Laplace equation, and the rule explains why the smallest harmonic space suffices.

The element defines no head inside the polygon, whereas the retention law needs one. The storage is therefore lumped on the nodes in the conservative form of \cite{Celia1990A}, and each node receives the storage of its own head with a weight supplied by the assumed flux. This is the lumped capacity that keeps a front advancing into dry soil free of oscillations, obtained here from the structure of the element rather than from a specific quadrature rule.

The conductivity may change by orders of magnitude within a single element, so that any discretisation must decide which conductivity acts between two nodes. Finite difference and finite volume schemes use an average of the two nodal values, arithmetic, geometric or integral, or the upstream value \cite{HaverkampVauclin1979,Szymkiewicz2009}. Here Darcy's law is written for the Kirchhoff potential of each element, which removes the conductivity from the diffusion operator. Between two nodes this naturally leads to the integral mean, that is, the mean of the conductivity over the range of heads between them, so that the average follows from the formulation and is not chosen arbitrarily. The transformation is applied within the element only. The head remains the continuous global unknown, and layered soils need no special treatment at their interfaces, where a global Kirchhoff variable would be discontinuous.

The gravity term, however, still contains the conductivity, and in a slope it drives a large part of the flow. The same integral mean already used in the diffusion term is adopted there as well. Diffusion and gravity then balance exactly at rest, and hydrostatic equilibrium and gravity drainage are reproduced on any mesh. On triangles the element coincides with the control-volume finite element scheme that uses the integral mean, and on polygons it generalises that scheme with no stabilisation.

Since every field of the element is the gradient of a known polynomial, all the element operators reduce to boundary integrals computed once. The nonlinear iterations then need only nodal values, which keeps the assembly inexpensive.

The element is assessed against primal and mixed hybrid finite elements, a stabilised and a stabilisation-free virtual element and, as a reference for a Gardner soil, a hybrid-Trefftz element with exact Trefftz functions. The tests range from solutions in closed form to infiltration into a dry layered soil, a seepage face resolved by local refinement and rainfall on a homogeneous and on a layered slope, on polygonal, cut and locally refined meshes. Where no solution in closed form exists, an independent solution computed with a commercial code is used as a further reference. The tests are chosen to show where the formulations differ, in conservation, in the treatment of gravity and in cost, and where they agree.

The numerical results show that the element is free of spurious modes on distorted polygons and on elements with hanging nodes, and that it converges at second order in the head, also where gravity dominates. The water balance closes in every test, on each element at every iterate and on the whole domain at convergence, and the nonlinear iterations remain inexpensive. On a layered slope, cut grids refined in a band along the surface, where the rain acts, reach a given accuracy with a fraction of the nodes of uniform refinement, and the solutions agree with those computed independently. To the author's knowledge, grids of this kind have not been used for unsaturated flow in layered slopes, and the refinement, here prescribed in advance, could be made adaptive. The main novelties of the work are thus the particular flux that makes the element balance exact, the storage lumped with weights supplied by the assumed flux, the Kirchhoff formulation within each element with a single integral mean for diffusion and gravity, and the ordering of the flux basis by rotation, with its closed-form stability rule.

The work is organised as follows. Section~\ref{sec:continuum} states the continuum problem. Section~\ref{sec:hvem} derives the element, and Section~\ref{sec:interp} addresses the interpolation on a polygon, the stability condition, hanging nodes and boundary conditions. Section~\ref{numerical.s} presents the numerical results, and Section~\ref{sec:concl} draws the conclusions.

\section{Continuum problem}\label{sec:continuum}

This section recalls Richards' equation and the constitutive laws of the soil, and completes the problem with its initial and boundary conditions.

\subsection{Governing equations}

A porous medium occupies the bounded domain $\Omega\subset\mathbb R^2$, referred to a Cartesian frame $(x,z)$ with the $z$ axis pointing upwards, so that $\mathbf e_z$ is the upward unit vector and the elevation head is $z$. The skeleton is rigid, except for its contribution to the specific storage, and the flow is followed over the interval $[0,T]$. The unknown fields are the pressure head $h$, the Darcy flux $\bq$ and the volumetric water content $\theta$. The soil is characterised by the retention curve $\Theta(h)$, the conductivity $K(h)$, the residual and saturated water contents $\theta_r$ and $\theta_s$, the saturated conductivity $K_s$, the porosity $\phi$ and the specific storage $S_s$, which may vary in space, piecewise, between materials. The balance of mass, Darcy's law and the retention law read
\begin{equation}\label{eq:richards}
\frac{\partial\theta}{\partial t}
+S_s\,S_w\,\frac{\partial h}{\partial t}
+\nabla\!\cdot\bq=s,
\qquad
\bq=-K(h)\bigl(\nabla h+\mathbf e_z\bigr),
\qquad
\theta=\Theta(h),
\end{equation}
where $s$ is a volumetric source and $S_w=\theta/\phi$ the degree of saturation. The second term, the elastic storage, accounts for the compressibility of water and skeleton and keeps the problem parabolic once the retention curve has become flat at saturation. The conductivity is isotropic in all the applications. An anisotropic saturated conductivity $\mathbf K_s$, constant within each element, can be reduced to it by an affine change of coordinates in each element, the element operators being computed on the mapped polygon, where gravity acts along a fixed but no longer vertical direction.

The head is the primary unknown in the saturated zone as well, and the storage is kept in the conservative form $\partial_t\Theta(h)$.

\subsection{Constitutive laws}

The retention curve $\Theta(h)$ and the conductivity $K(h)$ of Eq.~\eqref{eq:richards} are both written in terms of the effective saturation
\begin{equation}\label{eq:Se}
S_e=\frac{\theta-\theta_r}{\theta_s-\theta_r}\in[0,1],
\end{equation}
so that the retention curve is
\begin{equation}\label{eq:storage}
\Theta(h)=\theta_r+(\theta_s-\theta_r)\,S_e(h),
\end{equation}
and the conductivity is $K(h)=K_s\,K_r(h)$, with $K_r\in(0,1]$ the relative conductivity.

Three models are used. The van Genuchten--Mualem model \cite{vanGenuchten1980,Mualem1976} reads
\begin{equation}\label{eq:vg}
S_e(h)=\bigl[1+|\alpha h|^{\,n_{\mathrm{vG}}}\bigr]^{-m_{\mathrm{vG}}},
\qquad
K_r=S_e^{1/2}\Bigl[1-\bigl(1-S_e^{1/m_{\mathrm{vG}}}\bigr)^{m_{\mathrm{vG}}}\Bigr]^{2},
\qquad m_{\mathrm{vG}}=1-\frac1{n_{\mathrm{vG}}},
\end{equation}
for $h<0$, and $S_e=K_r=1$ for $h\ge0$; the data give $\alpha$ through $h_\alpha=-1/\alpha$, of the order of the air-entry value, and a larger $n_{\mathrm{vG}}$ gives a steeper curve and a sharper front.

Gardner's model \cite{Gardner1958} is
\begin{equation}\label{eq:gardner}
S_e(h)=e^{\alpha h},\qquad K_r=e^{\alpha h}.
\end{equation}

The model of Bergamaschi and Putti \cite{BergamaschiPutti1999} reads
\begin{equation}\label{eq:bp}
S_e(h)=\bigl[1+\bigl(\alpha(h_s-h)\bigr)^{\beta}\bigr]^{-\gamma},
\qquad
K_r=S_e^{\,\mu}
\quad\text{or}\quad
K_r=10^{\,aS_e^2+(b-2a)S_e+a-b},
\end{equation}
for $h<h_s$, with $h_s\le0$ the air-entry head, $\beta$ and $\gamma$ playing the part of $n_{\mathrm{vG}}$ and $m_{\mathrm{vG}}$, and $\mu$, $a$ and $b$ fitting coefficients; its data are given through $\phi=\theta_s$ and the residual saturation $S_{wr}=\theta_r/\phi$. When $\alpha$ and $\gamma$ are negative the saturation is defined only above a lower bound of the head, which the implementation enforces by clipping the Newton iterate.

Nothing in the formulation depends on the particular model. It is required that $\Theta$ be bounded and non-decreasing and that $K$ be bounded and positive; the inverse of $K$ never appears, so that a conductivity vanishing in very dry soil needs no lower bound. The derivatives of $\Theta$ and $K$, needed only for the linearisation, must exist almost everywhere.

\subsection{Initial and boundary conditions}

The problem is completed by the initial condition $h(\bx,0)=h_0(\bx)$ and by conditions on three disjoint parts of the boundary, $\partial\Omega=\Gamma_D\cup\Gamma_N\cup\Gamma_S$:
\begin{equation}\label{eq:bcs}
h=\overline h \ \ \text{on }\Gamma_D,
\qquad
\bq\cdot\bn=\overline q_n \ \ \text{on }\Gamma_N,
\end{equation}
with $\bn$ the outward normal and $\overline q_n$ positive outwards. A rainfall of intensity $R$ per unit horizontal area falling on an inclined surface gives $\overline q_n=-R\,|\bn\cdot\mathbf e_z|$. The rain is assumed to infiltrate entirely, which holds as long as the intensity does not exceed the infiltration capacity of the soil; otherwise the rainfall boundary can be treated as a potential seepage face, which limits the head to zero and lets the excess run off. On a potential seepage face $\Gamma_S$ water leaves the domain at atmospheric pressure where the soil behind it is saturated, and the face is impervious elsewhere:
\begin{equation}\label{eq:seepage}
h\le 0,\qquad \bq\cdot\bn\ge0,\qquad
h\,(\bq\cdot\bn)=0
\qquad\text{on }\Gamma_S .
\end{equation}
The exit point between the two regimes is part of the solution, and the condition of Eq.~\eqref{eq:seepage} turns the problem into a free boundary problem.

\section{The hybrid virtual element for unsaturated flow}\label{sec:hvem}

The hybrid virtual element is derived in this section. The discrete fields and the time discretisation are introduced first, and the element equations follow, together with the treatment of gravity and storage and the nodal equations solved at each step.

\subsection{Discrete fields}\label{sec:fields}

The domain $\Omega$ is partitioned into a mesh of polygonal elements $\Oe$, so that $\Omega=\bigcup_e\Oe$, and $\sum_e$ denotes the sum over the elements. A generic element $\Oe$ has $n$ vertices and $n$ straight edges, area $\Ae$, boundary $\dE$, centroid $\bx_G=(x_G,z_G)$ and outward normal $\bn$, and $\bX=(X,Z)=\bx-\bx_G$ denotes the position relative to the centroid. The quantities of the element carry the subscript $e$, the global ones none.

Following the VEM, the explicit interpolation of the head $h$ inside the element domain is not given. Only its trace $\widetilde h$ along the element boundary $\dE$ is interpolated as
\begin{equation}\label{eq:trace}
\widetilde h=\Nh\,\bh_e\qquad\text{on }\dE,
\end{equation}
where $\bh_e\in\mathbb R^{n}$ collects the nodal heads of $\Oe$ and $\Nh$ the corresponding interpolation functions along $\dE$. The nodal heads are the only global dofs, collected in $\bh$.

Within each element an independent Darcy flux $\bq$ is assumed. It is polynomial, discontinuous across the element boundaries and eliminated at the element level. It is the sum of a divergence-free part, which carries the free parameters, and of a particular part, which carries the divergence,
\begin{equation}\label{eq:flux}
\bq=\Nq(\bx)\,\br+c_e\,\mathbf p(\bx)+\mathbf P(\bx)\,\mathbf c_e,
\qquad
\nabla\!\cdot\Nq=\mathbf 0,
\end{equation}
where $\Nq$ is a $2\times n_r$ matrix of linearly independent divergence-free polynomial fields, $\br$ the vector of the $n_r$ flux parameters, and $c_e$ and $\mathbf c_e=(c_{e,1},c_{e,2})$ three scalars. Each column of $\Nq$ is the gradient of a harmonic polynomial, and the choice of these fields, which depends on the geometry of the polygon alone, governs the stability of the element. The particular fields are gradients of radial polynomials centred at the centroid,
\begin{equation}\label{eq:pfield}
\mathbf p=\nabla\frac{|\bX|^2}{4}=\frac{\bX}2,
\qquad
\mathbf P=\bigl[\,\mathbf p_1\ \big|\ \mathbf p_2\,\bigr]
=\Bigl[\,\nabla\frac{X|\bX|^2}{8}\ \Big|\ \nabla\frac{Z|\bX|^2}{8}\,\Bigr],
\end{equation}
with $\nabla\!\cdot\mathbf p=\tfrac12\nabla\!\cdot\bX=1$, $\nabla\!\cdot\mathbf p_1=X$ and $\nabla\!\cdot\mathbf p_2=Z$, so that the divergence of the assumed flux is the linear function
\begin{equation}\label{eq:divq}
\nabla\!\cdot\bq=c_e+c_{e,1}\,X+c_{e,2}\,Z .
\end{equation}
All the fields of the element are gradients, and therefore irrotational, and none of them introduces a preferred direction.

\subsection{Time discretisation}\label{sec:time}

Over the time step $[t_n,t_{n+1}]$ the backward Euler rule applied to both storage terms of Eq.~\eqref{eq:richards} gives the storage rate
\begin{equation}\label{eq:euler}
\sigma(h)=\frac{\Theta(h)-\Theta(h^{n})}{\Delta t}
+\frac{S_s}{\phi}\,\Theta(h)\,\frac{h-h^{n}}{\Delta t},
\end{equation}
where $h$ is the head at $t_{n+1}$ and $h^n$ the converged head at $t_n$; flux, conductivity, source and boundary data are evaluated at $t_{n+1}$. Each step is therefore the nonlinear stationary problem
\begin{equation}\label{eq:stepproblem}
\nabla\!\cdot\bq=s-\sigma(h),
\qquad
\bq=-K(h)\bigl(\nabla h+\mathbf e_z\bigr),
\end{equation}
with the boundary conditions of Eqs.~\eqref{eq:bcs} and~\eqref{eq:seepage}. The scheme is fully implicit and first-order accurate in time. The capillary part of Eq.~\eqref{eq:euler} is the increment of a function of the state, so that its sum over the steps telescopes and the water held in the pores is conserved exactly; the elastic part is not, and its volume is accumulated step by step in the water budget.

\subsection{Balance of the element}\label{sec:flux}

The divergence of the assumed flux is fixed by the first of Eqs.~\eqref{eq:stepproblem}, the balance of mass. By the divergence theorem, whatever $\br$,
\begin{equation}\label{eq:divthm}
\int_{\dE}\bq\cdot\bn\,d\Gamma=\int_{\Oe}\nabla\!\cdot\bq\,d\Omega=\Ae\,c_e ,
\end{equation}
since $X$ and $Z$ have zero mean over the element. The balance of the element requires the outflow through its boundary to equal the source minus the storage,
\begin{equation}\label{eq:elembal}
\Ae\,\overline\sigma_e+\int_{\dE}\bq\cdot\bn\,d\Gamma=\Ae\,\overline s_e ,
\end{equation}
with $\overline\sigma_e$ and $\overline s_e$ the means of the storage rate and of the source over the element. By Eq.~\eqref{eq:divthm} it holds if
\begin{equation}\label{eq:cE}
c_e=\overline s_e-\overline\sigma_e,
\qquad
\mathbf c_e=\mathbf s_e-\boldsymbol\sigma_e ,
\end{equation}
where the second condition, with $\mathbf s_e$ and $\boldsymbol\sigma_e$ the linear coefficients of the projections of the source and of the storage rate onto the linear functions, makes the first moments of the balance hold as well. The three scalars are therefore not additional unknowns, and the balance of Eq.~\eqref{eq:elembal} holds identically, for any flux parameters, hence at every iterate of the nonlinear solver and for any time step. The linear part of the divergence costs nothing, the fields $\mathbf p_1$ and $\mathbf p_2$ entering only through moments computed once, and it distributes the source and the storage within the element consistently with their variation, where a uniform divergence would see only their mean.

The balance concerns the flux of the element itself. The fluxes of two adjacent elements are independent, and the continuity of their normal components on the shared edge is imposed separately, in weak form. HVEM thus reverses the roles of the mixed hybrid $RT_0$ scheme, whose flux is $H(\mathrm{div})$-conforming but whose element balance is a discrete equation, satisfied only at convergence. What remains to be determined is the divergence-free part of the flux, through Darcy's law.

\subsection{Darcy's law for the Kirchhoff potential}\label{sec:kirchhoff}

The divergence-free part of the flux is determined by Darcy's law, the second of Eqs.~\eqref{eq:stepproblem}, which has to be imposed inside the element, where no head is interpolated and the conductivity, a function of the head, cannot be evaluated. Darcy's law is therefore rewritten in a variable that removes the conductivity from the diffusive term. In each material the Kirchhoff potential
\begin{equation}\label{eq:kirchhoff}
\Phi(h)=\int_{-\infty}^{h}K(\eta)\,d\eta
\end{equation}
is a strictly increasing function of the head, with $d\Phi/dh=K$, so that $K(h)\nabla h=\nabla\Phi(h)$ and Darcy's law splits into a part linear in $\Phi$ and a gravity part,
\begin{equation}\label{eq:darcykir}
\bq+\nabla\Phi(h)+K(h)\,\mathbf e_z=\mathbf 0 .
\end{equation}
The diffusion no longer contains the conductivity, which survives only in the gravity term, and the diffusion operator of the element will depend on its geometry alone. Across a material interface $\Phi$ is discontinuous even where $h$ is continuous, which is the difficulty of the methods that take $\Phi$ as global unknown \cite{Adak2021The}; here the global unknown remains the head, continuous across the element boundaries, and $\Phi$ is used only inside each element, computed with the retention and conductivity laws of its own material. The trace of the potential on the boundary of $\Oe$ is taken linear on each edge between the nodal values,
\begin{equation}\label{eq:phitrace}
\widetilde\Phi_e=\Nh\,\boldsymbol\Phi_e,
\qquad
\Phi_{e,i}=\Phi_e(h_i),
\end{equation}
which is as accurate as the trace of Eq.~\eqref{eq:trace} of the head and is exact where $\Phi$ varies linearly, as it does in steady one-dimensional flow without gravity. For van Genuchten's law $\Phi$ has no closed form and is tabulated once per material, from the dry end so that the small values of dry soil keep their relative precision, and interpolated by cubic Hermite polynomials whose derivative is $K$ itself.

\subsection{Weak forms}\label{sec:weak}

Darcy's law is a constitutive relation and is imposed in weighted form within each element, and the continuity of the normal flux between elements is imposed in weighted form on the element boundaries.

The residual of Eq.~\eqref{eq:darcykir} is made orthogonal, element by element, to every admissible variation $\delta\bq$ of the assumed flux:
\begin{equation}\label{eq:wr1}
\int_{\Oe}\delta\bq\cdot\bigl(\bq+\nabla\Phi+K\,\mathbf e_z\bigr)\,
d\Omega=0 .
\end{equation}
The potential enters through its interior gradient, which is not available. The divergence theorem gives
\begin{equation}\label{eq:byparts}
\int_{\Oe}\delta\bq\cdot\nabla\Phi\,d\Omega
=\int_{\dE}(\delta\bq\cdot\bn)\,\Phi\,d\Gamma
-\int_{\Oe}(\nabla\!\cdot\delta\bq)\,\Phi\,d\Omega ,
\end{equation}
and the last integral vanishes if the variations are divergence-free,
\begin{equation}\label{eq:soltest}
\nabla\!\cdot\delta\bq=0\qquad\text{in }\Oe .
\end{equation}
This is the case for $\delta\bq=\Nq\,\delta\br$, the variation of the only free part of the flux, the three scalars of the particular part being fixed by the balance of the element. The interior potential then disappears exactly, and, with the trace of Eq.~\eqref{eq:phitrace}, Eq.~\eqref{eq:wr1} becomes
\begin{equation}\label{eq:weak}
\int_{\Oe}\delta\bq\cdot\bq\,d\Omega
+\int_{\dE}\widetilde\Phi_e\,(\delta\bq\cdot\bn)\,d\Gamma
+\int_{\Oe}\delta\bq\cdot K\,\mathbf e_z\,d\Omega=0
\qquad\forall\,\delta\bq=\Nq\,\delta\br .
\end{equation}
The choice of divergence-free fields for the free part of the flux serves precisely this purpose. The first term of Eq.~\eqref{eq:weak} contains no conductivity, and depends on the geometry of the element alone.

The exact flux has a normal component continuous across the interior edges and equal to $\overline q_n$ on $\Gamma_N$. For any function $v$ defined on the element boundaries, taking the same values on the two sides of each interior edge and vanishing on $\Gamma_D$, the contributions of each interior edge to the sum over the elements cancel, and
\begin{equation}\label{eq:transcont}
\sum_e\int_{\dE}(\bq\cdot\bn)\,v\,d\Gamma
=\int_{\Gamma_N}\overline q_n\,v\,d\Gamma
+\int_{\Gamma_S}(\bq\cdot\bn)\,v\,d\Gamma .
\end{equation}
This identity is imposed on the assumed flux with $v$ the variation of the head trace, $\delta\widetilde h=\Nh\,\delta\bh_e$, which has the same values on the two sides of each interior edge, since the nodal heads are shared, and vanishes where the head is prescribed,
\begin{equation}\label{eq:wr2}
\sum_e\int_{\dE}(\bq\cdot\bn)\,\delta\widetilde h\,d\Gamma
=\int_{\Gamma_N}\overline q_n\,\delta\widetilde h\,d\Gamma
\qquad\forall\,\delta\bh ,
\end{equation}
the nodes of a seepage face being either of prescribed head or free, according to their state. Since the nodal values are shared by adjacent elements, Eq.~\eqref{eq:wr2} gives one equation for each node whose head is not prescribed. As the assumed flux is discontinuous, it imposes the continuity of the normal flux only in a weak sense, weighted by the trace, which acts as the Lagrange multiplier of the continuity and allows the flux to be condensed element by element.

For a unit conductivity, for which $\widetilde\Phi_e=\widetilde h$, and for a given divergence of the flux and a given gravity term, Eqs.~\eqref{eq:weak} and~\eqref{eq:wr2} are the stationarity conditions of
\begin{equation}\label{eq:functional}
\Pi(\bq,\widetilde h)=\sum_e\Bigl[
\frac12\int_{\Oe}\bq\cdot\bq\,d\Omega
+\int_{\Oe}\bq\cdot K\,\mathbf e_z\,d\Omega
+\int_{\dE}\widetilde\Phi_e\,(\bq\cdot\bn)\,d\Gamma
\Bigr]
-\int_{\Gamma_N}\overline q_n\,\widetilde h\,d\Gamma ,
\end{equation}
the counterpart, written for the Kirchhoff potential, of the hybrid-stress functional of Pian \cite{Pian1964}. With a nonlinear conductivity the trace of the potential depends on the head through $K$, which on a material interface differs between the elements sharing a node, and the divergence and the gravity term depend on the head as well, so that the discrete equations are not the stationarity conditions of any functional and their tangent operator is not symmetric.

\subsection{Element equations in condensed form}\label{sec:condensation}

Introducing the assumed flux of Eq.~\eqref{eq:flux} into Eq.~\eqref{eq:weak} with $\delta\bq=\Nq\,\delta\br$ gives, for any $\delta\br$,
\begin{equation}\label{eq:condsys}
\mathbf M_q\,\br+\mathbf B\,\boldsymbol\Phi_e+\mathbf k_e
+c_e\,\mathbf m_p+\mathbf M_P\,\mathbf c_e=\mathbf 0 ,
\end{equation}
with the gravity load $\mathbf k_e=\int_{\Oe}\Nq^{\!\top}K\,\mathbf e_z\,d\Omega$ and the geometric matrices
\begin{equation}\label{eq:arrays}
\mathbf M_q=\int_{\Oe}\Nq^{\!\top}\Nq\,d\Omega,\quad
\mathbf m_p=\int_{\Oe}\Nq^{\!\top}\mathbf p\,d\Omega,\quad
\mathbf M_P=\int_{\Oe}\Nq^{\!\top}\mathbf P\,d\Omega,\quad
\mathbf B=\int_{\dE}\Nq^{\!\top}\bn\,\Nh\,d\Gamma ;
\end{equation}
the entry $(k,j)$ of $\mathbf B$ is the work of the normal component of the $k$-th flux field on the interpolation function of the $j$-th vertex. Since the fields of $\Nq$ are linearly independent, $\mathbf M_q$ is positive definite, and the flux parameters follow by static condensation,
\begin{equation}\label{eq:cond}
\br=-\mathbf M_q^{-1}\bigl[\mathbf B\,\boldsymbol\Phi_e+\mathbf k_e
+c_e\,\mathbf m_p+\mathbf M_P\,\mathbf c_e\bigr].
\end{equation}
The flux enters Eq.~\eqref{eq:wr2} through the nodal outflows of the element,
\begin{equation}\label{eq:bnodal}
\begin{gathered}
\mathbf f_e=\int_{\dE}\Nh^{\!\top}(\bq\cdot\bn)\,d\Gamma
=\mathbf B^{\!\top}\br+c_e\,\mathbf b_p+\mathbf B_P\,\mathbf c_e,\\
\mathbf b_p=\int_{\dE}\Nh^{\!\top}(\mathbf p\cdot\bn)\,d\Gamma,\qquad
\mathbf B_P=\int_{\dE}\Nh^{\!\top}(\bn^{\!\top}\mathbf P)\,d\Gamma ,
\end{gathered}
\end{equation}
since $\delta\widetilde h=\Nh\,\delta\bh_e$ turns the left-hand side of Eq.~\eqref{eq:wr2} into $\sum_e\delta\bh_e^{\!\top}\mathbf f_e$. Substituting Eq.~\eqref{eq:cond},
\begin{equation}\label{eq:subst}
\mathbf f_e=-\mathbf A_e^0\,\boldsymbol\Phi_e-\mathbf B^{\!\top}\mathbf M_q^{-1}\mathbf k_e
+c_e\,\mathbf d_e+\mathbf D_P\,\mathbf c_e ,
\end{equation}
with
\begin{equation}\label{eq:AE}
\mathbf A_e^0=\mathbf B^{\!\top}\mathbf M_q^{-1}\mathbf B,\qquad
\mathbf d_e=\mathbf b_p-\mathbf B^{\!\top}\mathbf M_q^{-1}\mathbf m_p,\qquad
\mathbf D_P=\mathbf B_P-\mathbf B^{\!\top}\mathbf M_q^{-1}\mathbf M_P ,
\end{equation}
all three depending on the geometry alone.

The matrix $\mathbf A_e^0$ is symmetric and positive semi-definite. For a uniform head, $\bh_e$ with all entries equal to $h$, the trace is $\widetilde h=h$, and since the columns of $\Nq$ are divergence-free
\begin{equation}\label{eq:B1}
\mathbf B\,\bh_e=h\int_{\dE}\Nq^{\!\top}\bn\,d\Gamma
=h\int_{\Oe}\nabla\!\cdot\Nq^{\!\top}\,d\Omega=\mathbf 0 .
\end{equation}
Hence the rows of $\mathbf A_e^0$ sum to zero and a uniform head is in its kernel. Any other kernel depends on $\Nq$ alone, not on the conductivity. Moreover, since $\nabla\!\cdot\mathbf p=1$ and $\nabla\!\cdot\mathbf p_1$, $\nabla\!\cdot\mathbf p_2$ have zero mean,
\begin{equation}\label{eq:dsum}
\sum_{i=1}^{n}d_i=\Ae ,
\qquad
\sum_{i=1}^{n}D_{P,ij}=0 ,\quad j=1,2 ,
\end{equation}
with $d_i$ and $D_{P,ij}$ the entries of $\mathbf d_e$ and $\mathbf D_P$. The entries of $\mathbf f_e$ then sum to $\Ae\,c_e$ for any nodal heads, which is the element balance of Eq.~\eqref{eq:elembal}. A uniform divergence $c_e=1$ produces the nodal outflows $\mathbf d_e$, which distribute the area of the element among its vertices.

Since the rows of $\mathbf A_e^0$ sum to zero, the diffusive part of the outflows is a sum of exchanges between pairs of nodes,
\begin{equation}\label{eq:pairdiff}
\bigl[\mathbf A_e^0\,\boldsymbol\Phi_e\bigr]_i
=\sum_{j\ne i}\bigl(-A^0_{e,ij}\bigr)\bigl(\Phi_i-\Phi_j\bigr)
=\sum_{j\ne i}\bigl(-A^0_{e,ij}\bigr)\,K_{ij}\,\bigl(h_i-h_j\bigr),
\end{equation}
where
\begin{equation}\label{eq:Kpair}
K_{ij}=\frac{\Phi_e(h_i)-\Phi_e(h_j)}{h_i-h_j}
=\frac{1}{h_i-h_j}\int_{h_j}^{h_i}K(\eta)\,d\eta
\end{equation}
is the mean of the conductivity over the range of heads between the two nodes. Each pair of nodes exchanges a diffusive flux driven by its difference of head, with a conductance given by the geometry, $-A^0_{e,ij}$, and by the Kirchhoff mean $K_{ij}$. When the two heads differ by less than a relative tolerance of $10^{-8}$ the quotient is replaced by its expansion, $K_{ij}=(K_i+K_j)/2$, with $\partial K_{ij}/\partial h_i=\tfrac12\,(dK/dh)(h_i)$ and likewise for $h_j$, which avoids the cancellation and keeps the tangent continuous. Among the internodal conductivities of finite difference and finite volume schemes \cite{HaverkampVauclin1979,Szymkiewicz2009,Belfort2013} this is the integral mean, which gives the exact steady flux between two nodes in the absence of gravity.

\subsection{Gravity}\label{sec:gravity}

The gravity load $\mathbf k_e$ of Eq.~\eqref{eq:condsys} contains the conductivity inside the element, where no head is defined. Its discretisation is guided by two states of the continuum that a scheme should reproduce exactly. The first is the hydrostatic equilibrium, with $h+z$ uniform and no flow, which holds in every slope at rest. The second is the gravity drainage at unit gradient, with $h$ uniform and $\bq=-K(h)\,\mathbf e_z$, which holds below a sustained infiltration. A single value of $K$ over the element reproduces the second state but not the first, because in the diffusion term of Eq.~\eqref{eq:pairdiff} each pair of nodes has its own mean $K_{ij}$, and no single value balances all of them. The gravity term is therefore written with the same pairwise means.

Only divergence-free variations enter Eq.~\eqref{eq:weak}, so the gravity field $K\,\mathbf e_z$ acts on them through the trace of a potential $Y$ with $\nabla Y\approx K\,\mathbf e_z$,
\begin{equation}\label{eq:gravload}
\delta\br^{\!\top}\mathbf k_e=\int_{\Oe}\delta\bq\cdot K\,\mathbf e_z\,d\Omega
\approx\int_{\dE}\widetilde Y\,(\delta\bq\cdot\bn)\,d\Gamma
=\delta\br^{\!\top}\mathbf B\,\mathbf y_e ,
\end{equation}
where $\widetilde Y=\Nh\,\mathbf y_e$ is linear on each edge. The nodal values $\mathbf y_e$ are chosen so that the difference of $Y$ between two vertices equals the Kirchhoff mean of the pair times their difference of elevation. The condition cannot hold for every pair at once and is imposed in the least-squares sense,
\begin{equation}\label{eq:ygrav}
\mathbf y_e=\arg\min_{\mathbf y}\ \sum_{i<j}w_{ij}\,
\bigl[y_i-y_j-K_{ij}\,(z_i-z_j)\bigr]^2,
\qquad
w_{ij}=\bigl|A^0_{e,ij}\bigr| .
\end{equation}
The matrix of the normal equations is the weighted graph Laplacian of the element, which depends on the geometry alone. The constant left free in $\mathbf y_e$ is annihilated by $\mathbf B$. The fit is exact in the two states above. At a hydrostatic state $\mathbf y_e=-\boldsymbol\Phi_e$, and the gravity load cancels the diffusive one. At a uniform head $\mathbf y_e=K\,\mathbf z_e$, and $\mathbf B\,\mathbf y_e=K\int_{\Oe}\Nq^{\!\top}\mathbf e_z\,d\Omega$ exactly.

The weights $w_{ij}$ are the absolute values of the entries of $\mathbf A_e^0$. On obtuse triangles and on polygons these entries have both signs. With the signed entries as weights, two large entries of opposite sign would multiply slightly different means, and on a degenerate element the flux would grow without bound. With positive weights $\mathbf y_e$ remains a smooth potential on any polygon.

With gravity no symmetric mean is exact. Upstream weighting of the conductivity \cite{Forsyth1995} favours monotonicity on coarse grids, but has no place here. The diffusion term contains no conductivity to weight, and applied to the gravity term alone the upstream value would weight the diffusion and the gravity of each pair by different conductivities and break the hydrostatic equilibrium. The same mean is therefore kept in both terms.

With $\mathbf k_e=\mathbf B\,\mathbf y_e$, the flux parameters of Eq.~\eqref{eq:cond} and the nodal outflows of Eq.~\eqref{eq:subst} become
\begin{equation}\label{eq:gravsubst}
\br=-\mathbf M_q^{-1}\bigl[\mathbf B\,(\boldsymbol\Phi_e+\mathbf y_e)
+c_e\,\mathbf m_p+\mathbf M_P\,\mathbf c_e\bigr],
\qquad
\mathbf f_e=-\mathbf A_e^0\,\bigl(\boldsymbol\Phi_e+\mathbf y_e\bigr)
+c_e\,\mathbf d_e+\mathbf D_P\,\mathbf c_e ,
\end{equation}
and diffusion and gravity combine in the pairwise exchanges
\begin{equation}\label{eq:pairflux}
\bigl[\mathbf A_e^0\,(\boldsymbol\Phi_e+\mathbf y_e)\bigr]_i
=\sum_{j\ne i}\bigl(-A^0_{e,ij}\bigr)\,
\bigl[(\Phi_i+y_i)-(\Phi_j+y_j)\bigr]
=\sum_{j\ne i}\bigl(-A^0_{e,ij}\bigr)\,K_{ij}\,
\bigl(H_i-H_j\bigr),
\end{equation}
with $H_i=h_i+z_i$ the total head at the node. The second equality is exact where the off-diagonal entries of $\mathbf A_e^0$ are non-positive, since there $w_{ij}=-A^0_{e,ij}$ and the normal equations of Eq.~\eqref{eq:ygrav} are precisely these nodal sums; elsewhere it holds approximately. Each pair of nodes exchanges a flux driven by its difference of total head, with a conductance given by the geometry and by the Kirchhoff mean of the conductivity. At a hydrostatic state, $H_i$ uniform, the flux vanishes on any mesh and for any soil, and a domain at rest stays at rest. At a uniform head the flux is the exact gravity drainage. Both properties follow from the exactness of the fit of Eq.~\eqref{eq:ygrav} in these states.

\subsection{Storage}\label{sec:storage}

The mean storage rate $\overline\sigma_e$ and its linear coefficients $\boldsymbol\sigma_e$ in Eq.~\eqref{eq:cE} require a head inside the element, which the formulation does not define. The storage rate is therefore evaluated at the nodes, $\sigma_i=\sigma(h_i)$ from Eq.~\eqref{eq:euler}, and distributed over the element with the weights $\mathbf d_e$ of Eq.~\eqref{eq:AE}, which sum to $\Ae$ by Eq.~\eqref{eq:dsum},
\begin{equation}\label{eq:lump}
\Ae\,\overline\sigma_e=\sum_i d_i\,\sigma_i ,
\qquad
\Ae\,\overline\theta_e=\sum_i d_i\,\Theta(h_i).
\end{equation}
The coefficients $\boldsymbol\sigma_e$ are those of the linear function with the same first moments as the nodal distribution, $\sum_i d_i\,\sigma_i\,\bX_i$. The water balanced by Eq.~\eqref{eq:elembal} is then $\Ae\,\overline\theta_e$. It differs from the integral of the water content of the nodal interpolant by a term of second order in the mesh size.

Through the flux, Eq.~\eqref{eq:gravsubst} would distribute the storage of the element among its nodes as $\overline\sigma_e\,\mathbf d_e+\mathbf D_P\,\boldsymbol\sigma_e$. The storage is instead lumped, and each node receives the storage $d_i\,\sigma_i$ of its own head. The two distributions have the same total, $\sum_i d_i\sigma_i=\Ae\,\overline\sigma_e$ by Eq.~\eqref{eq:dsum}, and differ by a vector of zero sum. With the lumped storage the storage part of the tangent is the diagonal matrix $\mathrm{diag}(d_i\,\partial\sigma_i/\partial h_i)$, positive as long as the weights are. This is the lumped capacity that Celia et al.~\cite{Celia1990A} found necessary to avoid oscillations ahead of a front advancing into dry soil. Here the weights are supplied by the assumed flux. The source, if present, is the only quantity integrated over the element. A sink depending on the head, such as root uptake, is treated like the storage.

The weights have a physical meaning, which also explains their sign. With $c_e=1$, $\mathbf c_e=\mathbf 0$ and no potential and gravity, Eq.~\eqref{eq:cond} gives $\br=-\mathbf M_q^{-1}\mathbf m_p$, and the flux
\begin{equation}\label{eq:qd}
\bq_d=\mathbf p-\Nq\,\mathbf M_q^{-1}\mathbf m_p
\end{equation}
is the field of unit divergence orthogonal in $L^2(\Oe)$ to the divergence-free fields of $\Nq$. By Eq.~\eqref{eq:bnodal} its nodal outflows are $\mathbf d_e$. A unit source spread over the element leaves it through the boundary, and $d_i$ is the part collected by vertex $i$. On any polygon the weights sum to $\Ae$, Eq.~\eqref{eq:dsum}, and their first moments about the centroid vanish, $\sum_i d_i X_i=\sum_i d_i Z_i=0$, because the trace reproduces $X$ and $Z$ exactly and $\Nq$ contains the uniform fields. The nodal rule $\sum_i d_i f(\bx_i)$ therefore integrates linear functions exactly, and the lumping is of second order.

The weights are positive on most polygons, but not on all. If $\Nq$ contained the gradients of all the harmonic functions, $\bq_d$ would be the steady drainage of a uniform source towards a boundary at uniform potential, whose outflow is nowhere negative, and every weight would be positive. With the few harmonic polynomials of $\Nq$ this field is only approximated. The error is largest where the outflow is smallest, at an acute vertex or at a vertex between short edges, and there the weight can become slightly negative. Such weights are rare and small, and on a mesh they are compensated by those of the neighbouring elements, but nothing excludes a negative nodal capacity.

On an element with a negative weight, the weights of the lumped storage and of the source are therefore replaced by
\begin{equation}\label{eq:dsafe}
\mathbf d_e(t)=(1-t)\,\mathbf d_e+t\,\mathbf d_e^{\Delta},
\end{equation}
where $\mathbf d_e^{\Delta}$ are the lumped weights of the linear triangles of a triangulation of the element on its vertices, each vertex receiving one third of the area of every triangle that contains it, and $t\in(0,1]$ is the smallest value for which no weight is negative. Both vectors sum to $\Ae$ and integrate linear functions exactly, and so does their combination. The replacement keeps the balance of the element and the order of the lumping. It changes only the distribution of the storage among the nodes, by a vector of zero sum as the lumping itself does, and leaves unchanged every element whose weights are positive. With it the weights, and hence the lumped capacity, are non-negative on any polygon.

\subsection{Nodal equations}\label{sec:cond}

With the coefficients of Eq.~\eqref{eq:cE}, the nodal outflows of Eq.~\eqref{eq:gravsubst} with the lumped storage and the assembly $\mathbf f_N$ of $\int_{\Gamma_N}\Nh^{\!\top}\overline q_n\,d\Gamma$ over the edges of $\Gamma_N$, Eq.~\eqref{eq:wr2} gives the nodal equations
\begin{equation}\label{eq:res}
\mathcal R(\bh)=\sum_e
\Bigl[\mathbf A_e^0\,\bigl(\boldsymbol\Phi_e+\mathbf y_e\bigr)
+\mathrm{diag}(\mathbf d_e)\,\boldsymbol\sigma_e^{\mathrm{nod}}
-\overline s_e\,\mathbf d_e-\mathbf D_P\,\mathbf s_e\Bigr]
+\mathbf f_N=\mathbf 0 ,
\end{equation}
restricted to the free nodes, where each element adds its vector to the equations of its own nodes and $\boldsymbol\sigma_e^{\mathrm{nod}}$ are the nodal storage rates. At a node of prescribed head the residual is minus the outflow through the boundary.

The flux $\bq$ is recovered from Eq.~\eqref{eq:gravsubst} with the coefficients of Eq.~\eqref{eq:cE}, so that it satisfies the element balance of Eq.~\eqref{eq:elembal} identically. Eq.~\eqref{eq:flux} then gives it at any point of the element, and at the centroid $\bq(\bx_G)=\Nq(\bx_G)\,\br$, the value used to assess the flux in the numerical tests. Its nodal outflows differ from those entering Eq.~\eqref{eq:res} by the vector of zero sum due to the lumping. The nodal equations therefore enforce the weak continuity of the normal flux between elements up to this redistribution of the storage within each element, and, summed over the nodes, the balance of the whole domain, with the outflow at the prescribed heads given by the nodal reactions. Both hold once the nodal equations are satisfied, whereas the balance of each element holds at every iterate.

The pairwise fluxes in the first form of Eq.~\eqref{eq:pairflux} are antisymmetric in $i$ and $j$ on any element, so that, once the nodal equations are satisfied, the water is balanced not only on every element, identically, but also on the control volumes surrounding the nodes, with explicit fluxes between them. The off-diagonal entries of $\mathbf A_e^0$ are non-positive on triangles without obtuse angles and on rectangles with aspect ratio up to $\sqrt2$; there, with the conductances $K_{ij}$ frozen, the operator of Eq.~\eqref{eq:pairflux} plus the lumped storage is an M-matrix acting on the total head, as for the control-volume schemes that preserve the maximum principle. On general polygons some of these entries are positive, as for every scheme of this kind, and monotonicity is not guaranteed.

\subsection{Element operators as boundary integrals}\label{sec:boundary}

Every field of the element is the gradient of a known polynomial, the flux fields of $\Nq$ of harmonic polynomials $\psi_a$, and $\mathbf p$, $\mathbf p_1$, $\mathbf p_2$ and $\mathbf e_z$ of $w=|\bX|^2/4$, $v_1=X|\bX|^2/8$, $v_2=Z|\bX|^2/8$ and $z$. Green's identities then reduce every matrix of Eq.~\eqref{eq:arrays} to an integral over the boundary, for instance
\begin{equation}\label{eq:green}
\begin{gathered}
\int_{\Oe}\nabla\psi_a\cdot\nabla\psi_b\,d\Omega
=\int_{\dE}\psi_a\,\partial_n\psi_b\,d\Gamma,
\qquad
\int_{\Oe}\nabla\psi_a\cdot\nabla w\,d\Omega
=\int_{\dE}w\,\partial_n\psi_a\,d\Gamma,\\
\int_{\Oe}\nabla\psi_a\cdot\mathbf e_z\,d\Omega
=\int_{\dE}\psi_a\,n_z\,d\Gamma ,
\end{gathered}
\end{equation}
the integrals over $\Oe$ of the Laplacian of a harmonic function vanishing, and the remaining integrals of homogeneous polynomials of degree $k$ in $\bX$ follow from Euler's theorem,
\begin{equation}\label{eq:euler2}
\int_{\Oe} f\,d\Omega=\frac1{k+2}\int_{\dE}f\,(\bX\cdot\bn)\,d\Gamma .
\end{equation}
With Gauss rules on the edges these integrals are exact; they coincide with those of an interior quadrature to round-off, and need neither a sub-triangulation nor the convexity of the element. Within the nonlinear loop the element requires only the nodal values of $\Phi$, $K$ and $\Theta$: there is no quadrature at all, inside the element or on its boundary.

\subsection{Linearisation}\label{sec:tangent}

Eq.~\eqref{eq:res} is solved by Newton's method with a backtracking line search on the norm of the residual. The head enters the element residual through the nodal potentials, the pairwise conductances and the nodal storage rates, and the element tangent is
\begin{equation}\label{eq:tangent}
\frac{\partial\mathcal R_e}{\partial\bh_e}
=\mathbf A_e^0\,\mathrm{diag}\bigl(K(h_i)\bigr)
\;+\;\mathbf A_e^0\,\frac{\partial\mathbf y_e}{\partial\bh_e}
\;+\;\mathrm{diag}\Bigl(d_i\,\frac{\partial\sigma_i}{\partial h_i}\Bigr),
\end{equation}
where the three terms are related to the diffusion, gravity and storage, respectively, and the derivatives of the pairwise means follow from Eq.~\eqref{eq:Kpair},
\begin{equation}\label{eq:dKpair}
\frac{\partial K_{ij}}{\partial h_i}
=\frac{K(h_i)-K_{ij}}{h_i-h_j},
\qquad
\frac{\partial\sigma_i}{\partial h_i}=\frac1{\Delta t}
\Bigl[C_i+\frac{S_s}{\phi}\bigl(C_i\,(h_i-h_i^{n})
+\Theta_i\bigr)\Bigr],
\end{equation}
with $C_i$ and $\Theta_i$ the values at $h_i$ of the moisture capacity $C=d\Theta/dh$ and of the water content. The tangent is not symmetric, the first term being a geometric matrix scaled column by column, and is assembled from nodal quantities only. The complete tangent is used throughout.

\section{Interpolation on polygonal elements}\label{sec:interp}

This section specifies the interpolation of the head and of the flux on a polygon, and the number of flux fields that makes the element stable. Triangles, elements with hanging nodes and the boundary conditions are then discussed.

\subsection{Head trace}

\begin{figure}[h]
\centering
\subfloat[element, nodes and edges]{\makebox[0.44\linewidth]{%
\includegraphics{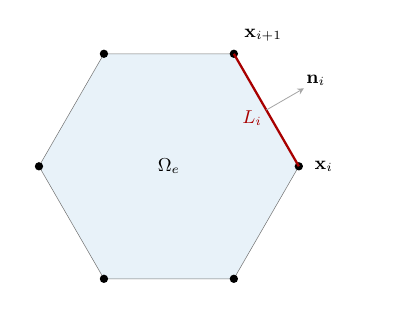}}}
\subfloat[the trace $\widetilde h$]{\makebox[0.44\linewidth]{%
\includegraphics{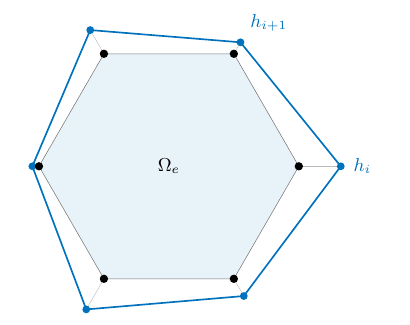}}}
\caption{Interpolation of the head on a polygon with vertices $\bx_i$ and edges of length $L_i$ and outward normal $\bn_i$. The trace varies linearly on every edge between the nodal values, the only global dofs; in (b) the nodal value is drawn as an offset normal to the boundary.}\label{trace.f}
\end{figure}

On each edge the head is linear and continuous across elements, $\widetilde h=\Nh\bh_e$ (Eq.~\eqref{eq:trace}): on a triangle this is the trace of the linear element, on a quadrilateral that of the bilinear one, and on a general polygon a continuous, piecewise linear function of the arc length (Fig.~\ref{trace.f}). The interior of the element is never interpolated.

\def\vlist{{0.95,0.30,0.62,0.15,0.80,0.45}}

\subsection{Stability condition}\label{sec:stab}

Since $\mathbf M_q$ is positive definite, $\bh_e^{\!\top}\mathbf A_e^0\bh_e=0$ if and only if $\mathbf B\bh_e=\mathbf 0$: the kernel of $\mathbf A_e^0$ is that of $\mathbf B$. It always contains the uniform head, by Eq.~\eqref{eq:B1}, and the element is stable when it contains nothing else,
\begin{equation}\label{eq:rank}
\operatorname{rank}\mathbf A_e^0=\operatorname{rank}\mathbf B=n-1 .
\end{equation}
A further kernel vector would be a nodal pattern to which no flux field responds, that is a zero-energy mode. In a primal virtual element the consistency term uses a polynomial projection of the head, which is blind to some patterns, and the stabilisation term controls them; for instance, on the unit square the hourglass pattern $(1,-1,1,-1)$ has zero mean on every edge and its projection onto the linear functions vanishes identically. Here no projection of the head is involved, and such modes are prevented by choosing $\Nq$ rich enough for $\mathbf B$ to have full rank on the complement of the constants, a condition on geometry and basis alone, since the conductivity does not enter $\mathbf A_e^0$.

\begin{figure}[h]
\centering
\subfloat[$\ell=0$: the kernel]{\makebox[0.235\linewidth]{%
\includegraphics{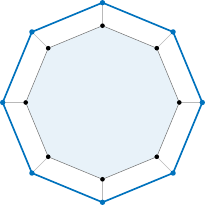}}}
\subfloat[$\ell=1$]{\makebox[0.235\linewidth]{%
\includegraphics{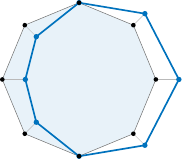}}}
\subfloat[$\ell=2$]{\makebox[0.235\linewidth]{%
\includegraphics{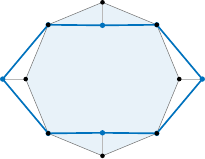}}}
\subfloat[$\ell=3$]{\makebox[0.235\linewidth]{%
\includegraphics{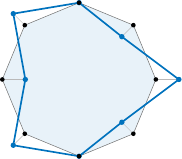}}}
\caption{Nodal head patterns on a regular octagon, drawn as in Fig.~\ref{trace.f}; only the cosine component of each mode of Eq.~\eqref{eq:fourier} is shown. The uniform mode $\ell=0$ is the constant head; each of the others must be reached by some column of $\Nq$, or it would be a spurious zero-energy mode.}\label{headmodes.f}
\end{figure}

On a regular $n$-gon with vertices at the polar angles $\vartheta_j=2\pi j/n$, any nodal vector decomposes into the discrete Fourier modes
\begin{equation}\label{eq:fourier}
h_j=\cos(\ell\vartheta_j),\quad h_j=\sin(\ell\vartheta_j),
\qquad \ell=0,1,\dots,\lfloor n/2\rfloor ,
\end{equation}
shown in Fig.~\ref{headmodes.f} for an octagon; only the constant, $\ell=0$, is an admissible kernel. Counting does not tell which flux fields control these modes. The divergence-free polynomial fields of degree at most $p$ span $2,5,9,14,\dots$ dimensions, so that a truncation with at least $n-1$ fields is always available, but a truncation by degree retains fields that couple with modes already controlled and omits others, and leaves spurious modes. What decides which field reaches which mode is the behaviour under rotation.

\subsection{Harmonic fields ordered by rotation}\label{sec:rot}

The flux fields are built from the harmonic polynomials. In the polar coordinates $(r,\vartheta)$ of the local coordinates $\boldsymbol\xi=(\bx-\bx_G)/\sqrt{\Ae}$, centred at the centroid, these are
\begin{equation}\label{eq:psi}
\psi_k^{c}=\frac{r^k\cos k\vartheta}{k},\qquad
\psi_k^{s}=\frac{r^k\sin k\vartheta}{k},\qquad k=1,2,\dots,
\end{equation}
that is $\xi_1$ and $\xi_2$ for $k=1$, $(\xi_1^2-\xi_2^2)/2$ and $\xi_1\xi_2$ for $k=2$, $(\xi_1^3-3\xi_1\xi_2^2)/3$ and $(3\xi_1^2\xi_2-\xi_2^3)/3$ for $k=3$. Their gradients
\begin{equation}\label{eq:rotmode}
\nabla_{\boldsymbol\xi}\psi_k^{c},\qquad \nabla_{\boldsymbol\xi}\psi_k^{s},\qquad k=1,2,\dots,
\end{equation}
are polynomial fields of degree $k-1$, divergence-free since $\Delta\psi_k=0$ and irrotational since they are gradients, that is harmonic. Up to degree $p$ they are $2(p+1)$, against the $2,5,9,14,\dots$ divergence-free fields of \S\ref{sec:stab}; the fields left out, such as the rigid rotation $(-\xi_2,\xi_1)$, have non-zero curl. The restriction is not arbitrary: in Eq.~\eqref{eq:darcykir} the diffusive flux $-\nabla\Phi$ is a gradient, the gravity term enters the weak form only through the trace of a potential (\S\ref{sec:gravity}), and the divergence is carried by the gradient fields $\mathbf p$, $\mathbf p_1$ and $\mathbf p_2$, so that the homogeneous part of the flux is harmonic. That every field is the gradient of a known polynomial is also what reduces the element operators to boundary integrals (\S\ref{sec:boundary}).

On a circle centred at the centroid, the normal component of these fields is
\begin{equation}\label{eq:rotindex}
\frac{\partial\psi_k^{c}}{\partial r}=r^{k-1}\cos k\vartheta,\qquad
\frac{\partial\psi_k^{s}}{\partial r}=r^{k-1}\sin k\vartheta ,
\end{equation}
and oscillates $k$ times around the element: $k$ is the \emph{rotational index} of the field. A rotation of the frame shifts $\vartheta$ and turns each of the two fields of index $k$ into a combination of the same two, never into fields of another index, so that the index does not depend on the orientation of the element.

\begin{figure}[h]
\centering
\subfloat[$\mathbf N_0$, index $k=1$]{\makebox[0.30\linewidth]{%
\includegraphics{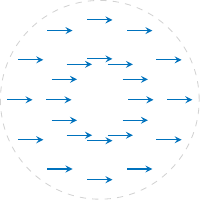}}}
\subfloat[$\mathbf N_1$, index $k=2$]{\makebox[0.30\linewidth]{%
\includegraphics{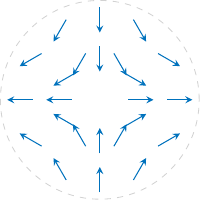}}}
\subfloat[$\mathbf N_2$, index $k=3$]{\makebox[0.30\linewidth]{%
\includegraphics{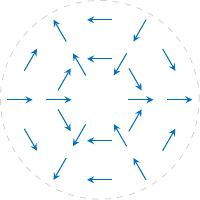}}}
\caption{The first harmonic field of the blocks $\mathbf N_0$, $\mathbf N_1$ and $\mathbf N_2$ of Eq.~\eqref{eq:Nm}, drawn as directions on two circles centred at the centroid. The field of $\mathbf N_m$ completes $m$ turns while the point completes one, and its normal component oscillates $k=m+1$ times, Eq.~\eqref{eq:rotindex}; the second field of each block is the same pattern rotated by a right angle.}\label{fluxmodes.f}
\end{figure}

On a regular $n$-gon, a head mode of Eq.~\eqref{eq:fourier} of index $\ell$ likewise oscillates $\ell$ times around the element. The entry of $\mathbf B$ that couples it with a flux field of index $k$, $\int_{\dE}(\bq\cdot\bn)\,\widetilde h\,d\Gamma$, is the sum of the contributions of the $n$ edges, and since the edges differ only by a rotation of $2\pi/n$, the contribution of the $j$-th edge is that of the first with the phases of the flux and of the head advanced by $2\pi jk/n$ and $2\pi j\ell/n$. Products of cosines and sines of these phases reduce to cosines and sines of $2\pi j(k\mp\ell)/n$, which, summed over the $n$ edges, vanish unless $k\mp\ell$ is a multiple of $n$, as in the discrete Fourier transform. The entry therefore vanishes unless
\begin{equation}\label{eq:sel}
k\equiv\pm\ell \pmod n ,
\end{equation}
that is unless the flux, seen edge by edge, oscillates as the head mode does. To reach every mode $\ell=1,\dots,\lfloor n/2\rfloor$ the basis must contain a field satisfying Eq.~\eqref{eq:sel} for each of them; the lowest indices are $k=\ell$, and retaining them all is also sufficient provided the couplings permitted by the rule do not vanish, which is the generic case. The adopted basis is
\begin{equation}\label{eq:iso}
\Nq=\bigl[\ \nabla_{\boldsymbol\xi}\psi_k^{c}\ \big|\ \nabla_{\boldsymbol\xi}\psi_k^{s}\
\bigr]_{k=1}^{\lfloor n/2\rfloor},
\qquad n_r=\dim\br=2\lfloor n/2\rfloor .
\end{equation}
It has $n-1$ parameters for odd $n$ and $n$ for even $n$, in which case the highest mode is one-dimensional and one parameter is redundant, retained to keep the basis independent of the orientation of the element. The two fields of index $k$, of degree $m=k-1$, form the block
\begin{equation}\label{eq:Nm}
\mathbf N_m=\bigl[\ \nabla_{\boldsymbol\xi}\psi_{m+1}^{c}\ \big|\ \nabla_{\boldsymbol\xi}\psi_{m+1}^{s}\ \bigr],
\qquad
\Nq=\bigl[\ \mathbf N_0\ \big|\ \mathbf N_1\ \big|\ \cdots\ \big|\
\mathbf N_{\lfloor n/2\rfloor-1}\ \bigr],
\end{equation}
the first being
\begin{equation}
\mathbf N_0=\begin{bmatrix}1&0\\0&1\end{bmatrix},\qquad
\mathbf N_1=\begin{bmatrix}\xi_1&\xi_2\\-\xi_2&\xi_1\end{bmatrix},\qquad
\mathbf N_2=\begin{bmatrix}\xi_1^2-\xi_2^2&2\xi_1\xi_2\\
-2\xi_1\xi_2&\xi_1^2-\xi_2^2\end{bmatrix}.
\end{equation}
While the point goes once around the centroid, the first field of $\mathbf N_m$ turns $m$ times in the opposite sense, and its normal component oscillates $m+1$ times (Fig.~\ref{fluxmodes.f}), as the head mode of the same index (Fig.~\ref{headmodes.f}): the blocks pair one by one with the head modes. The basis is given by a single formula for every $n$, needs no local reference frame, and for $n=3$ reduces to the constant fields which, with $\mathbf p$, span the lowest-order Raviart--Thomas space. The fields used on each polygon are listed explicitly in \ref{app:fields}.

\subsection{The element on a triangle}\label{sec:p1}

On a triangle $\Nq=\mathbf N_0=\mathbf I$, and the matrices of Eq.~\eqref{eq:arrays} are available in closed form: $\mathbf M_q=\Ae\,\mathbf I$, and by the divergence theorem $\mathbf B=\int_{\dE}\bn\,\Nh\,d\Gamma=\Ae\,\nabla\Nh$, with $\nabla\Nh$ the constant gradients of the linear shape functions. Hence
\begin{equation}\label{eq:p1stiff}
\mathbf A_e^0=\Ae\,(\nabla\Nh)^{\!\top}\nabla\Nh ,
\end{equation}
the stiffness matrix of the primal $P_1$ element for a unit conductivity, whose off-diagonal entries are the cotangent weights. Since $\mathbf p$ has zero mean, $\mathbf m_p=\mathbf 0$ and $\mathbf d_e=\mathbf b_p$; on each edge $\mathbf p\cdot\bn$ is half the distance of the edge from the centroid, and each vertex receives $\Ae/6$ from each of its two edges, so that
\begin{equation}\label{eq:p1lump}
\mathbf d_e=\frac{\Ae}{3}\,(1,1,1)=\int_{\Oe}\Nh^{\!\top}\,d\Omega ,
\end{equation}
the lumped distribution. On a triangle the nodal equations of HVEM are therefore those of the control-volume finite element scheme on the barycentric dual, with lumped storage and the Kirchhoff mean of the conductivity between each pair of nodes, for diffusion and gravity alike. The assumed flux adds nothing to the operator where the linear element is already optimal; it adds the element balance and a flux that satisfies it, and on more general elements it supplies the operator itself.

\subsection{Collinear edges and hanging nodes}\label{sec:hang}

\begin{figure}[h]
\centering
\subfloat[configuration]{\makebox[0.31\linewidth]{%
\includegraphics{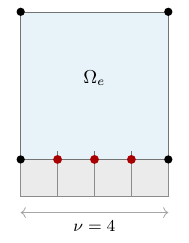}}}
\subfloat[pattern, $\nu=4$]{\makebox[0.31\linewidth]{%
\includegraphics{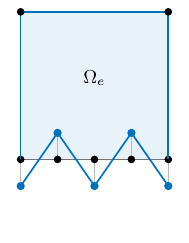}}}
\subfloat[pattern, $\nu=6$]{\makebox[0.31\linewidth]{%
\includegraphics{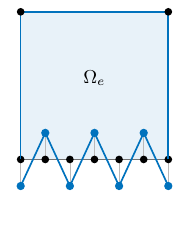}}}
\caption{An element carrying hanging nodes. In (a) the element meets neighbours $\log_2\nu$ levels finer, whose corners are $\nu-1$ additional vertices, in red, on its lower side; in (b) and (c) the trace of a nodal pattern confined to that side. The shape of the element is unchanged, while the variation that the normal flux must reproduce along the side grows with $\nu$, which the selection of Eq.~\eqref{eq:iso} does not control.}\label{hanggeom.f}
\end{figure}

The selection of Eq.~\eqref{eq:iso} is indexed by the number of edges, which is appropriate only as long as no two consecutive edges are aligned. Where a coarse element meets finer neighbours, its side carries their corners: the hanging node is an ordinary vertex of an ordinary polygon, and requires neither a transition element nor a constraint equation, which is what makes virtual elements well suited to adaptive refinement \cite{vanHuyssteen2025}. On a side made of $\nu$ collinear sub-edges, however, the trace has $\nu-1$ interior nodes on one straight line, and the normal component of the flux must vary enough along that line to distinguish them (Fig.~\ref{hanggeom.f}).

A direct computation of the rank of $\mathbf A_e^0$ on a square whose lower side is split into $\nu$ equal sub-edges shows that the selection of Eq.~\eqref{eq:iso} leaves a spurious mode from $\nu=6$ onwards, and so does the complete divergence-free space of the same degree. The basis is therefore sized on the collinear runs as well. On an element with collinear edges, let $\nu\ge2$ be the largest number of consecutive sub-edges lying on one straight side. The highest degree retained in Eq.~\eqref{eq:Nm} is then
\begin{equation}\label{eq:hangrule}
m_{\max}=\max\Bigl(\Bigl\lfloor\frac n2\Bigr\rfloor-1,\;
\Bigl\lceil\frac{\nu}2\Bigr\rceil,\;\nu-2\Bigr),
\qquad
\Nq=\bigl[\ \mathbf N_0\ \big|\ \cdots\ \big|\ \mathbf N_{m_{\max}}\bigr],
\end{equation}
and on an element without collinear edges it is that of Eq.~\eqref{eq:iso}. The last two terms were established by the same computation, for an increasing number of sub-edges, and restore Eq.~\eqref{eq:rank} for every $\nu$ up to eight. Unlike the first term, which follows from the rotational analysis of \S\ref{sec:rot} on regular polygons, they are an empirical rule, verified only up to $\nu=8$. All three are rank conditions, and do not bound the smallest non-zero eigenvalue of $\mathbf A_e^0$ independently of the shape. Eq.~\eqref{eq:hangrule} is the only flux interpolation used in this work. A quadtree balanced $2\!:\!1$ or $4\!:\!1$ gives $\nu\le4$, within the range of Eq.~\eqref{eq:iso}, and only an unbalanced jump of three levels calls for the enrichment.

\subsection{Boundary conditions}\label{sec:bc}

A prescribed head is imposed on the nodal values and a prescribed flux enters through $\mathbf f_N$. The seepage condition of Eq.~\eqref{eq:seepage} is imposed by an active-set iteration nested in the Newton loop: a node of the face whose head becomes positive is set to zero and added to the active set, and an active node is released when its nodal equation indicates inflow, that is when the residual, which at a node of prescribed head is minus the outflow (\S\ref{sec:cond}), is positive. Every point at which a boundary condition changes and every material interface must coincide with nodes of the mesh at all levels of refinement: a strip of prescribed flux ending inside an edge changes the applied flux by an amount of the order of the edge length, an error that does not decrease under refinement and against which the element balance offers no protection. Such points are therefore passed to the mesh generator as constraints, and domains with several materials are meshed region by region.

\section{Numerical results}\label{numerical.s}

The proposed approach is herein validated in some numerical tests. First, the stability properties of HVEM are assessed. Then the solution of HVEM is compared with analytical solutions under steady and transient conditions, measuring accuracy, local refinement and cost. Finally, four problems are considered: infiltration into a dry layered soil, a seepage face resolved by local refinement, and rainfall on a homogeneous and on a layered slope.

The element is compared with the discretisations of Table~\ref{elements.tab}, all implemented in the same code with the same boundary conditions, tolerances and time integration. Except for the mixed hybrid scheme, which follows with its own solver the implementation of \cite{BergamaschiPutti1999} and has no seepage condition, all the schemes are solved by the Newton iteration of the present formulation. The primal elements T3 and Q4 write the storage in the conservative form with a lumped capacity. Two virtual elements of lowest order are implemented as published. The first (VEM) \cite{Berrone2026} writes the storage in the capacitive form $C(h)\,\partial_t h$ and adds the \textsc{dofi-dofi} stabilisation, scaled by the element mean of $K$ in the diffusion and of $C$ in the storage. The second (SF), the stabilisation-free element of a recent hydro-mechanical model of slopes \cite{BerroneSlope2026}, projects the gradient onto divergence-free polynomial fields of the lowest degree that makes the projection injective \cite{BerroneDivFree2024} and lumps the capacitive storage. Both are published with a Picard iteration, which solves the same discrete equations, so that the common Newton iteration confines the comparison to the discretisations. In every scheme the tangent matrix is assembled with the derivatives of the constitutive laws in closed form, including, for the capacitive storage of SF and VEM, the second derivative of the retention curve. The hybrid-Trefftz element (HT) \cite{Jirousek1978,JirousekWroblewski1996,Qin2000} assumes within each element exact solutions of the governing equation; these exist for a Gardner soil, where HT serves as a reference. Additionally, independent solutions are computed with the commercial code Abaqus \cite{Abaqus2020}. The definition of each problem, the geometries, the meshes, the reference solutions, the data behind the tables and figures and the input files of the Abaqus analyses are openly available, as stated in the Data availability section.

The code is written in Julia~\cite{Bezanson2017}, version~1.11. Computing times are measured on one thread of a workstation with an Intel Core i9-12900K processor and $64$~GB of memory, under Windows~11.

\begin{table}[t]
\footnotesize
\centering
\caption{Discretisations compared.}\label{elements.tab}
\begin{tabular}{llll}
\toprule
name & formulation & global dofs & meshes\\
\midrule
\textbf{T3}   & primal $P_1$ \cite{Celia1990A}                 & one head per vertex & triangles\\
\textbf{Q4}   & primal $Q_1$ \cite{Paniconi1991}               & one head per vertex & quadrilaterals\\
\textbf{MHFE} & mixed hybrid $RT_0$ \cite{BergamaschiPutti1999} & one multiplier per edge & triangles\\
\textbf{VEM}  & primal virtual element, $k=1$, stabilised \cite{Beir2012Basic,Berrone2026} & one head per vertex & polygons\\
\textbf{SF}   & primal virtual element, stabilisation-free \cite{BerroneDivFree2024,BerroneSlope2026} & one head per vertex & polygons\\
\textbf{HT}   & hybrid-Trefftz, reference, Gardner soil \cite{Jirousek1978} & one head per vertex & polygons\\
\textbf{HVEM} & present & one head per vertex & polygons\\
\bottomrule
\end{tabular}
\end{table}

Accuracy is measured by the relative error
\begin{equation}\label{eq:relerr}
\mathcal E(f)=\frac{\|f_h-f\|}{\|f\|},
\end{equation}
where $\|\cdot\|$ is the Euclidean norm of the values at a set of sample points, $f$ the reference field and $f_h$ the computed one. The water budget is defined as the imbalance between the change of the stored water, as each scheme defines it, and the net inflow over the whole analysis.

\subsection{Stability assessment}\label{stability.s}

\begin{figure}[t]
\centering
\subfloat[$n=7$, $d=0$]{\makebox[0.18\linewidth]{\includegraphics{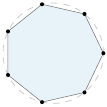}}}%
\subfloat[$n=7$, $d=0.5$]{\makebox[0.18\linewidth]{\includegraphics{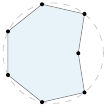}}}%
\subfloat[$n=7$, $d=1$]{\makebox[0.18\linewidth]{\includegraphics{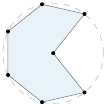}}}%
\subfloat[$\nu=2$, $n=5$]{\makebox[0.18\linewidth]{\includegraphics{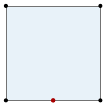}}}%
\subfloat[$\nu=6$, $n=9$]{\makebox[0.18\linewidth]{\includegraphics{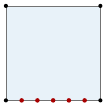}}}
\caption{Stability assessment: (a)--(c) a heptagon at three levels of distortion; (d), (e) a square whose lower side carries $\nu$ sub-edges, as left by a quadtree, with the $\nu-1$ hanging nodes in red.}\label{stabgeo.f}
\end{figure}

The stability properties of HVEM are here tested. The element is stable if the matrix $\mathbf A_e^0$ has exactly one null eigenvalue, for any number of edges and any shape of the polygon. To this end, the eigenvalues of $\mathbf A_e^0$ are evaluated on regular and distorted $n$-gons, as well as on elements with hanging nodes, as shown in Fig.~\ref{stabgeo.f}.
Table~\ref{stability.tab}a shows that HVEM is stable for every $n$ from $3$ to $10$ and every $d$, with the smallest non-zero eigenvalue of order one.
Since $\mathbf A_e^0$ depends on the geometry alone, the soil properties have no influence on these results. It is interesting to notice that the selection rule of Eq.~\eqref{eq:sel} remains effective when the symmetry on which it rests is lost.

\begin{table}[t]
\footnotesize
\centering
\caption{Stability assessment: number of flux parameters and of zero eigenvalues of $\mathbf A_e^0$ (a) on the polygons of Fig.~\ref{stabgeo.f}a--c, with the extreme non-zero eigenvalues, and (b) on the elements with hanging nodes of Fig.~\ref{stabgeo.f}d--e, for the basis sized on the number of edges, Eq.~\eqref{eq:iso}, the complete divergence-free space of the same degree and the adopted rule, Eq.~\eqref{eq:hangrule}.}\label{stability.tab}
\subfloat[distorted polygons]{%
\begin{tabular}{ll*{8}{c}}
\toprule
 & & \multicolumn{8}{c}{number of edges $n$}\\
\cmidrule(l{3pt}r{3pt}){3-10}
  &  & 3 & 4 & 5 & 6 & 7 & 8 & 9 & 10 \\
\midrule
\multicolumn{2}{l}{flux parameters $\dim\br$}
   & 2 & 4 & 4 & 6 & 6 & 8 & 8 & 10 \\
\multicolumn{2}{l}{zero eigenvalues}
   & 1 & 1 & 1 & 1 & 1 & 1 & 1 & 1 \\
\midrule
\multirow{2}{*}{$d=0$} & $\lambda_{\max}$ & 0.8660 & 1.0000 & 0.9511 & 0.9238 & 1.0031 & 1.0149 & 0.9933 & 1.0016 \\
  & $\lambda_{\min}$ & 0.8660 & 0.6667 & 0.7189 & 0.8248 & 0.6661 & 0.7071 & 0.6356 & 0.5878 \\
\multirow{2}{*}{$d=0.5$} & $\lambda_{\max}$ & 1.2990 & 1.3974 & 1.2701 & 1.2329 & 1.2364 & 1.2503 & 1.2408 & 1.2387 \\
  & $\lambda_{\min}$ & 0.5774 & 0.6026 & 0.6014 & 0.6309 & 0.4559 & 0.5366 & 0.3736 & 0.4561 \\
\multirow{2}{*}{$d=1$} & $\lambda_{\max}$ & 2.5981 & 2.3187 & 1.8639 & 1.7204 & 1.5996 & 1.5588 & 1.5153 & 1.4936 \\
  & $\lambda_{\min}$ & 0.2887 & 0.4313 & 0.3991 & 0.3921 & 0.2386 & 0.2854 & 0.1667 & 0.2022 \\
\bottomrule
\end{tabular}}\\[2mm]
\subfloat[hanging nodes\label{hangstab.tab}]{%
\begin{tabular}{ll*{8}{c}}
\toprule
 & & \multicolumn{8}{c}{sub-edges on one side, $\nu$}\\
\cmidrule(l{3pt}r{3pt}){3-10}
 & & 1 & 2 & 3 & 4 & 5 & 6 & 7 & 8\\
\multicolumn{2}{l}{vertices $n$}
   & 4 & 5 & 6 & 7 & 8 & 9 & 10 & 11\\
\midrule
\multirow{2}{*}{by edge count, Eq.~\eqref{eq:iso}}
   & $\dim\br$ & 4 & 4 & 6 & 6 & 8 & 8 & 10 & 10\\
   & zero eig. & 1 & 1 & 1 & 1 & 1 & \textbf{2} & \textbf{2} & \textbf{3}\\
\midrule
\multirow{2}{*}{complete, same degree}
   & $\dim\br$ & 5 & 5 & 9 & 9 & 14 & 14 & 20 & 20\\
   & zero eig. & 1 & 1 & 1 & 1 & 1 & \textbf{2} & \textbf{2} & \textbf{3}\\
\midrule
\multirow{2}{*}{adopted, Eq.~\eqref{eq:hangrule}}
   & $\dim\br$ & 4 & 4 & 6 & 6 & 8 & 10 & 12 & 14\\
   & zero eig. & 1 & 1 & 1 & 1 & 1 & 1 & 1 & 1\\
\bottomrule
\end{tabular}}
\end{table}

With hanging nodes (Fig.~\ref{stabgeo.f}d--e), when some element edges are collinear, the basis in Eq.~\eqref{eq:iso} loses rank, and so does the complete divergence-free space of the same polynomial degree, with twice as many parameters. Table~\ref{stability.tab}b shows that the rule proposed in Eq.~\eqref{eq:hangrule} restores the rank in every case.

\begin{figure}[t]
\centering

\subfloat[regular]{\includegraphics{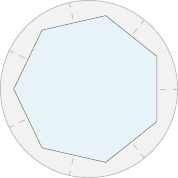}}\hspace{8mm}
\subfloat[random]{\includegraphics{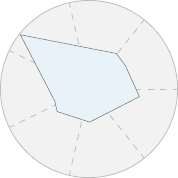}}\hspace{8mm}
\subfloat[random]{\includegraphics{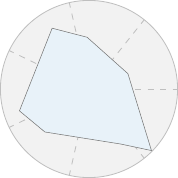}}
\caption{Stability assessment: a regular heptagon and two random ones of the stochastic test, each vertex lying in its own sector of the annulus $0.1\le R<1$.}\label{mcgeo.f}
\end{figure}

To exclude spurious modes on arbitrarily distorted polygons, the stochastic test of \cite{Liguori2025Plates,Liguori2026Shells} is repeated. Each of the $n$ vertices is randomly drawn in sectors of amplitude $2\pi/n$ of the annulus $0.1\le R<1$ (Fig.~\ref{mcgeo.f}), which gives star-shaped, often strongly non-convex polygons. The dimensionless spectrum of $\mathbf A_e^0$ is then computed on $10^3$ polygons for each $n$ from $3$ to $10$, in a stochastic search for the smallest non-zero eigenvalue \cite{Sun2025}.

\begin{figure}[t]
\centering
{\mbox{%
  \raisebox{2pt}{\includegraphics{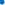}}\;{\footnotesize $\log_{10}\lambda_{\max}$}}\quad\mbox{%
  \raisebox{2pt}{\includegraphics{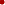}}\;{\footnotesize $\log_{10}\lambda_{\min}$, smallest non-zero}}\quad\mbox{%
  \raisebox{2pt}{\includegraphics{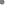}}\;{\footnotesize $\log_{10}|\lambda_0|$, constant head}}\quad}\\[1mm]
\includegraphics{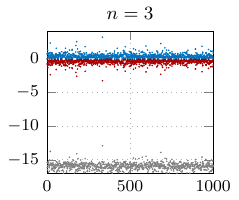}%
\includegraphics{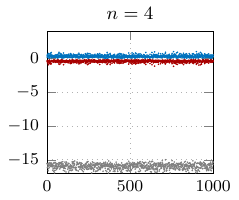}%
\includegraphics{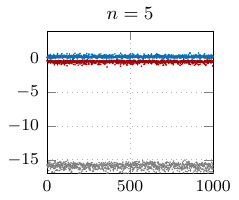}%
\includegraphics{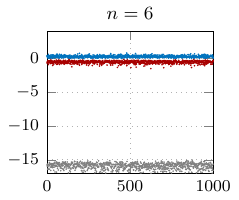}\\[1mm]
\includegraphics{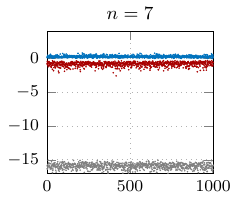}%
\includegraphics{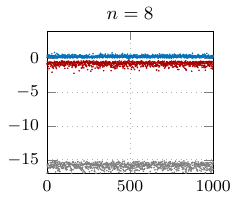}%
\includegraphics{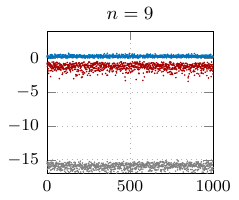}%
\includegraphics{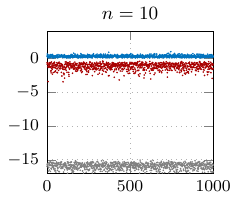}
\caption{Stability assessment: largest and smallest non-zero eigenvalue of $\mathbf A_e^0$, and eigenvalue of the constant head, which measures the numerical zero, for each of the $10^3$ random polygons of the stochastic test, in abscissa.}\label{mc.f}
\end{figure}

Figure~\ref{mc.f} shows the results of this procedure. In every polygon exactly one eigenvalue vanishes, to round-off, and the others stay many orders of magnitude above it. The smallest of them decreases as the number of edges grows, since more edges admit head patterns that oscillate more along the boundary \cite{Liguori2026Shells}. No zero-energy mode appears, although a test of this kind gives no bound uniform in the number of edges. The weights $\mathbf d_e$ of Eq.~\eqref{eq:AE} are positive on all the polygons of Fig.~\ref{stabgeo.f}, at the hanging nodes as well.

\begin{figure}[t]
\centering
\subfloat[$n=3$]{\makebox[0.235\linewidth]{\includegraphics{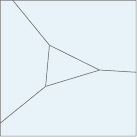}}}\hfill\subfloat[$n=4$]{\makebox[0.235\linewidth]{\includegraphics{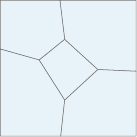}}}\hfill\subfloat[$n=5$]{\makebox[0.235\linewidth]{\includegraphics{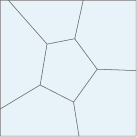}}}\hfill\subfloat[$n=6$]{\makebox[0.235\linewidth]{\includegraphics{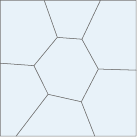}}}\\[1mm]
\subfloat[$n=7$]{\makebox[0.235\linewidth]{\includegraphics{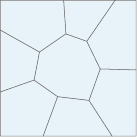}}}\hfill\subfloat[$n=8$]{\makebox[0.235\linewidth]{\includegraphics{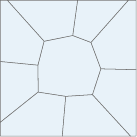}}}\hfill\subfloat[$n=9$]{\makebox[0.235\linewidth]{\includegraphics{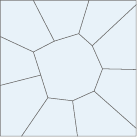}}}\hfill\subfloat[$n=10$]{\makebox[0.235\linewidth]{\includegraphics{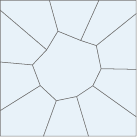}}}
\caption{Stability assessment: the eight squares of the patch test, each with a distorted central polygon of $n$ edges.}\label{patch.f}
\end{figure}

Finally, consistency is verified by the patch test on the eight squares of Fig.~\ref{patch.f}, each containing a distorted polygon with three to ten edges among elements of four to six edges. If linear heads are prescribed on the boundary and gravity is included, nodal head and flux are reproduced to round-off.

\subsection{Steady flow in a Gardner soil}\label{square.s}

\begin{figure}[t]
\centering
\subfloat[steady flow in a Gardner soil]{\makebox[0.235\linewidth]{%
\includegraphics{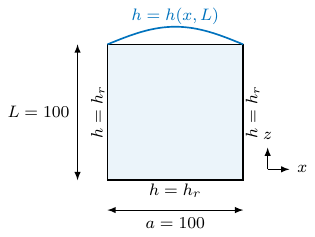}}}\hspace{0.06\linewidth}
\subfloat[quadtree, two levels]{\makebox[0.235\linewidth]{%
\includegraphics{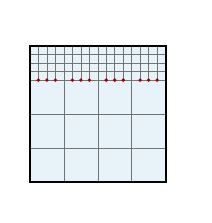}}}\\[1mm]
\subfloat[quadrilaterals]{\makebox[0.235\linewidth]{\includegraphics{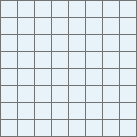}}}\hfill
\subfloat[distorted quadrilaterals]{\makebox[0.235\linewidth]{\includegraphics{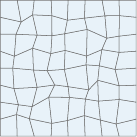}}}\hfill
\subfloat[Voronoi]{\makebox[0.235\linewidth]{\includegraphics{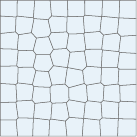}}}\hfill
\subfloat[Delaunay dual]{\makebox[0.235\linewidth]{\includegraphics{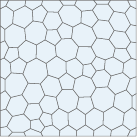}}}
\caption{Steady flow in a Gardner soil: (a) geometry in centimetres and boundary conditions; (b) quadtree mesh of the hanging-node test, refined by two levels towards the top side, with the hanging nodes in red; (c)--(f) the four mesh families, with $n=8$ elements per side.}\label{sqgeo.f}
\end{figure}

The first comparison with an analytical solution concerns the steady problem of Tracy \cite{Tracy2006}, defined on the square $(0,a)\times(0,L)$, $a=L=100$~cm, in a Gardner soil, Eq.~\eqref{eq:gardner}, with $K_s=20.464$~cm/d (Fig.~\ref{sqgeo.f}a). The head is equal to $h_r$ on the bottom and on the vertical sides, while on the top side it rises from $h_r$ at the corners to $0$ at the midpoint according to $h(x,L)=\alpha^{-1}\log[e^{\alpha h_r}+(1-e^{\alpha h_r})\sin(\pi x/a)]$. Letting $u_r=e^{\alpha h_r}$, the solution reads
\begin{equation}\label{eq:tracy}
h(x,z)=\frac1\alpha\log\Bigl[u_r+(1-u_r)\,
  \sin\frac{\pi x}{a}\;e^{\alpha(L-z)/2}\,
  \frac{\sinh\beta z}{\sinh\beta L}\Bigr],\qquad
\beta=\sqrt{\frac{\alpha^2}4+\Bigl(\frac\pi a\Bigr)^2},
\end{equation}
and the exact Darcy flux follows by differentiation. Two sets of parameters are considered. In the first, $\alpha=0.001$~cm$^{-1}$ and $h_r=-1000$~cm, so that $\alpha L=0.1$, gravity is weak and the conductivity varies by a factor of $e$ over the domain. In the second, $\alpha=0.05$~cm$^{-1}$ and $h_r=-100$~cm, so that $\alpha L=5$, gravity dominates and the conductivity varies by a factor of $150$. For this soil the Kirchhoff variable $u=e^{\alpha h}$ satisfies a linear equation, and the function $w=e^{\alpha z/2}u$ obeys $\Delta w=(\alpha/2)^2w$, whose solutions $I_k(\alpha r/2)\{\cos,\sin\}(k\vartheta)$ are exact Trefftz functions. The HT element is built with as many of these functions as the flux fields of HVEM. It is therefore the same hybrid element with exact local solutions in place of harmonic polynomials, and it indicates how much could still be gained from richer fields when, as in this case, a transformation makes the equation linear.

\begin{figure}[t]
\centering
{\mbox{%
  \raisebox{2pt}{\includegraphics{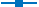}}\;{\footnotesize quadrilaterals}}\quad\mbox{%
  \raisebox{2pt}{\includegraphics{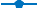}}\;{\footnotesize Voronoi}}\quad%
 \mbox{%
  \raisebox{2pt}{\includegraphics{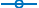}}\;{\footnotesize Delaunay dual}}\quad}\\[1mm]
\subfloat[head, $\alpha L=0.1$]{\includegraphics{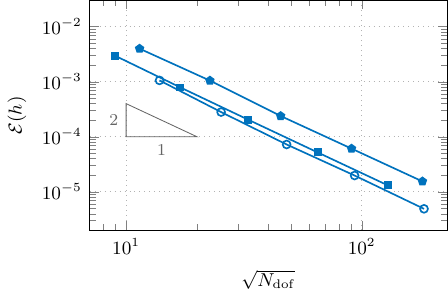}}\hfill
\subfloat[flux, $\alpha L=0.1$]{\includegraphics{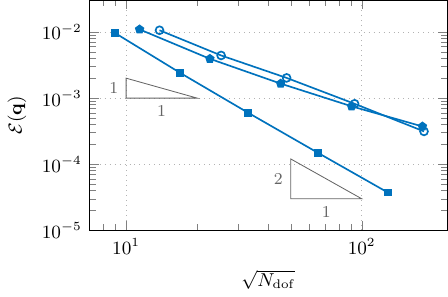}}\\
\subfloat[head, $\alpha L=5$]{\includegraphics{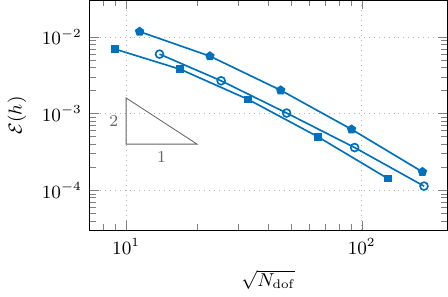}}\hfill
\subfloat[flux, $\alpha L=5$]{\includegraphics{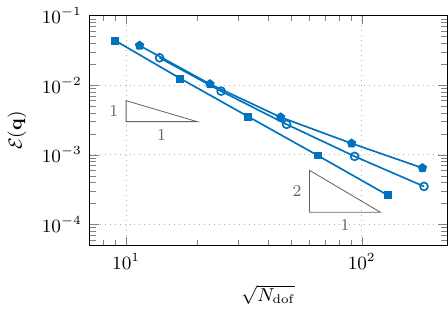}}
\caption{Steady flow in a Gardner soil: relative error of the nodal head and of the flux at the centroids of HVEM against $\sqrt{N_{\mathrm{dof}}}$ on three mesh families, with weak ($\alpha L=0.1$) and strong ($\alpha L=5$) gravity; the triangles mark the slopes of first and second order.}\label{tracy.f}
\end{figure}

The convergence of the proposed HVEM on three of the mesh families of Fig.~\ref{sqgeo.f}c--f is shown in Fig.~\ref{tracy.f} for both values of $\alpha L$. The error of the nodal head decreases at second order on every family. The error of the flux decreases at second order on the structured quadrilaterals, where the flux superconverges and its curve has the same slope as that of the head, and at a rate between first and second order on the polygons, about $1.5$ on the dual meshes and decreasing towards first order on the Voronoi meshes.

\begin{figure}[t]
\centering
{\mbox{%
  \raisebox{2pt}{\includegraphics{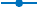}}\;{\footnotesize HVEM}}\quad\mbox{%
  \raisebox{2pt}{\includegraphics{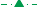}}\;{\footnotesize SF}}\quad\mbox{%
  \raisebox{2pt}{\includegraphics{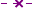}}\;{\footnotesize VEM}}\quad\mbox{%
  \raisebox{2pt}{\includegraphics{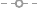}}\;{\footnotesize HT, exact Trefftz functions (reference)}}\quad}\\[1mm]
\subfloat[head, Voronoi]{\includegraphics{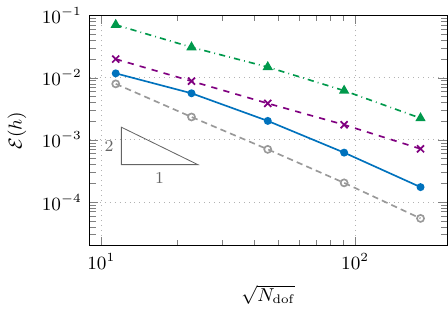}}\hfill
\subfloat[flux, Voronoi]{\includegraphics{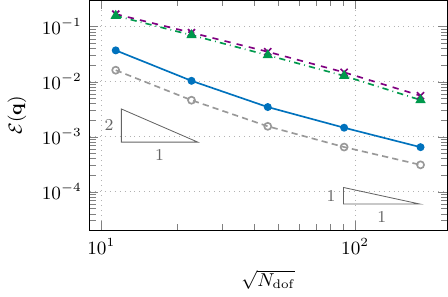}}\\
\subfloat[head, Delaunay dual]{\includegraphics{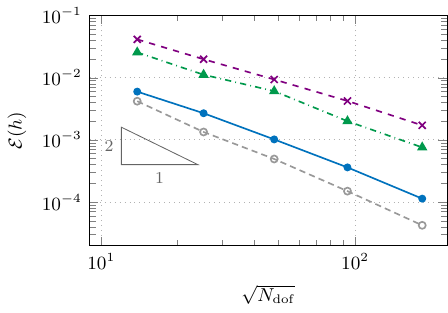}}\hfill
\subfloat[flux, Delaunay dual]{\includegraphics{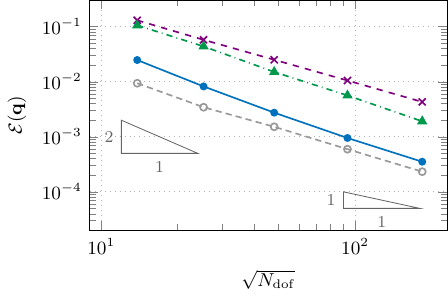}}
\caption{Steady flow in a Gardner soil, strong gravity ($\alpha L=5$): relative error of the nodal head and of the flux at the centroids against $\sqrt{N_{\mathrm{dof}}}$ on the polygonal meshes; HT, with exact Trefftz functions, is the reference, and the triangles mark the slopes of first and second order.}\label{tracycmp.f}
\end{figure}

The comparison with the other polygonal elements is first discussed for weak gravity, for which no figure is reported. In this case SF is close to HVEM in the head and somewhat less accurate in the flux, VEM is several times less accurate in both quantities, and HT gives practically the same results as HVEM, since the Kirchhoff potential is nearly harmonic. The picture changes when gravity dominates, as shown in Fig.~\ref{tracycmp.f} for the polygonal meshes. SF and VEM converge more slowly than HVEM and, on the finest meshes, their errors are up to an order of magnitude larger, both in the head and in the flux.

\begin{table}[t]
\footnotesize
\centering
\caption{Steady flow in a Gardner soil: relative error of the nodal head and of the flux at the centroids on the quadtree meshes of Fig.~\ref{sqgeo.f}b, refined by one, two or three levels, which leave $\nu=2$, $4$ or $8$ collinear sub-edges on the coarser elements; HT is the reference.}\label{hangcmp.tab}
\begin{tabular}{l *{6}{c} c *{6}{c}}
\toprule
 & \multicolumn{6}{c}{$10^4\,\mathcal E(h)$} & & \multicolumn{6}{c}{$10^3\,\mathcal E(\mathbf q)$}\\
\cmidrule(l{2pt}r{2pt}){2-7}\cmidrule(l{2pt}r{2pt}){9-14}
 & \multicolumn{3}{c}{$\alpha L=0.1$} & \multicolumn{3}{c}{$\alpha L=5$} & & \multicolumn{3}{c}{$\alpha L=0.1$} & \multicolumn{3}{c}{$\alpha L=5$}\\
\cmidrule(l{2pt}r{2pt}){2-4}\cmidrule(l{2pt}r{2pt}){5-7}\cmidrule(l{2pt}r{2pt}){9-11}\cmidrule(l{2pt}r{2pt}){12-14}
$\nu$ & 2 & 4 & 8 & 2 & 4 & 8 & & 2 & 4 & 8 & 2 & 4 & 8\\
\midrule
HVEM & 0.49 & 0.90 & 0.88 & 12 & 13 & 8.5 & & 0.67 & 1.1 & 1.0 & 0.82 & 1.3 & 0.98\\
SF & 0.48 & 0.99 & 1.0 & 22 & 22 & 22 & & 0.58 & 0.82 & 2.3 & 8.4 & 9.3 & 5.9\\
VEM & 0.51 & 1.2 & 1.0 & 25 & 50 & 57 & & 0.79 & 3.3 & 4.4 & 8.4 & 10 & 7.2\\
\midrule
HT, reference & 0.49 & 0.91 & 0.89 & 1.3 & 2.2 & 1.6 & & 0.67 & 1.1 & 1.0 & 0.82 & 1.3 & 1.2\\
\bottomrule
\end{tabular}
\end{table}

The behaviour with hanging nodes is tested on quadtree meshes refined towards the top side by one, two or three levels (Fig.~\ref{sqgeo.f}b). The coarser elements adjacent to the refined zone therefore have $\nu=2$, $4$ or $8$ collinear sub-edges, as in Fig.~\ref{stabgeo.f}d--e, and the hanging nodes are treated as ordinary vertices by all the virtual elements and by HT. The errors on the finest mesh of each sequence are reported in Table~\ref{hangcmp.tab}. With weak gravity the three virtual elements give similar heads, while in the flux HVEM is the least affected by the jump in refinement level. With strong gravity HVEM is the most accurate of the three, by almost an order of magnitude in the flux, and HT, which gives the same flux, is again more accurate in the head. It is worth noting that sizing the flux basis on the number of edges alone, as in Eq.~\eqref{eq:iso}, would change the flux of HVEM only slightly, because the heads at the hanging nodes are controlled by the finer elements that share them.

\subsection{Transient flow in a van Genuchten soil}\label{mms.s}

\begin{figure}[t]
\centering
\subfloat[transient, Eq.~\eqref{eq:mmst}]{\makebox[0.235\linewidth]{%
\includegraphics{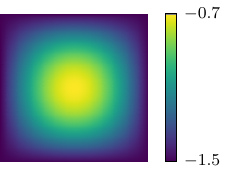}}}\hspace{0.08\linewidth}
\subfloat[steady, Eq.~\eqref{eq:mms2}]{\makebox[0.235\linewidth]{%
\includegraphics{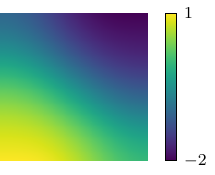}}}
\caption{Transient flow in a van Genuchten soil: the two manufactured solutions on the unit square, the transient one at $t=1$; the exact head is prescribed on the whole boundary.}\label{mms.f}
\end{figure}

The transient problem is tested through the manufactured solution
\begin{equation}\label{eq:mmst}
h(x,z,t)=-1.5+0.8\,t\,\sin(\pi x)\sin(\pi z),\qquad 0\le t\le1 ,
\end{equation}
defined on the unit square (Fig.~\ref{mms.f}a). The soil follows the van Genuchten--Mualem model with $\theta_r=0.05$, $\theta_s=0.40$, $h_\alpha=-1$, $n_{\mathrm{vG}}=2$, $K_s=1$ and $S_s=0$. The exact head is prescribed on the boundary, and the same time step $\Delta t=1/256$ is used on all the meshes.

All the elements converge at second order in the head on the mesh families of Fig.~\ref{sqgeo.f}c--f and on triangles. The errors of HVEM, SF and the primal elements are of the same size on every family, whereas VEM is less accurate on quadrilaterals and polygons. On triangles SF, VEM and T3 give the same result, since on this mesh the virtual elements reduce to the linear element. Conversely, HVEM, which on triangles coincides with the control-volume finite element scheme with the Kirchhoff mean, differs slightly from them, because it weights each pair of nodes by the mean conductivity between their heads rather than by its values at the quadrature points. All the schemes take three Newton iterations per step.

\begin{figure}[t]
\centering
{\mbox{%
  \raisebox{2pt}{\includegraphics{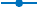}}\;{\footnotesize HVEM}}\quad\mbox{%
  \raisebox{2pt}{\includegraphics{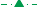}}\;{\footnotesize SF}}\quad\mbox{%
  \raisebox{2pt}{\includegraphics{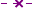}}\;{\footnotesize VEM}}\quad\mbox{%
  \raisebox{2pt}{\includegraphics{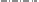}}\;{\footnotesize LU factorisation}}\quad}\\[1mm]
\subfloat[steady, Voronoi]{%
\includegraphics{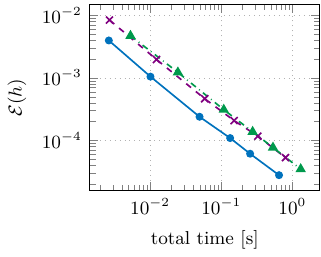}}\hfill
\subfloat[steady, Delaunay dual]{%
\includegraphics{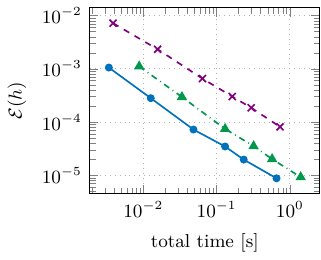}}\hfill
\subfloat[steady, Delaunay dual]{%
\includegraphics{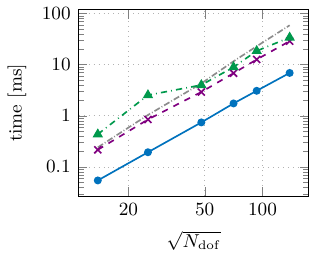}}\\
\subfloat[transient, Voronoi]{%
\includegraphics{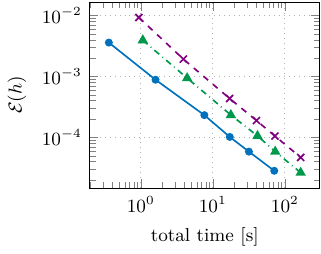}}\hfill
\subfloat[transient, Delaunay dual]{%
\includegraphics{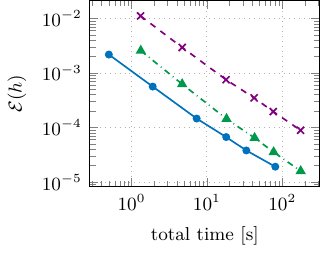}}\hfill
\subfloat[transient, Delaunay dual]{%
\includegraphics{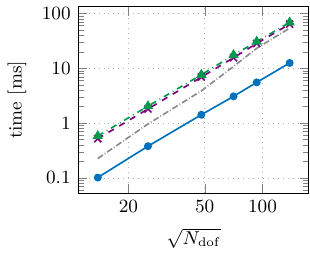}}
\caption{Transient flow in a van Genuchten soil: (a), (b), (d), (e) relative error of the nodal head against the computing time of the whole analysis, for the steady problem of Tracy with $\alpha L=0.1$ and for the transient solution of Eq.~\eqref{eq:mmst}; (c), (f) time of one assembly of residual and matrix for each element, and of one sparse factorisation (dash-dotted line), the same for all the elements.}\label{cpu.f}
\end{figure}

Since the accuracy of the elements is similar, their cost becomes the discriminating factor. Figure~\ref{cpu.f} compares the three virtual elements on the same meshes, for the steady problem of Tracy and for the transient one. HVEM assembles residual and tangent from nodal values only, since its operators are computed once as boundary integrals. SF and VEM instead evaluate the conductivity, the capacity and their derivatives at the quadrature points of every element at each iteration. As a consequence, one assembly of HVEM costs a fraction of theirs and less than the sparse factorisation. In the transient problem, where the many steps make the assembly dominant, HVEM reaches a given accuracy in less than half the time of the other two. In the steady problem the factorisation weighs more and the gain is smaller.

\begin{figure}[t]
\centering
{\mbox{%
  \raisebox{2pt}{\includegraphics{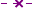}}\;{\footnotesize VEM}}\quad\mbox{%
  \raisebox{2pt}{\includegraphics{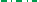}}\;{\footnotesize SF}}\quad\mbox{%
  \raisebox{2pt}{\includegraphics{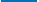}}\;{\footnotesize HVEM}}\quad\mbox{%
  \raisebox{2pt}{\includegraphics{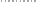}}\;{\footnotesize $\tau=1$}}\quad}\\[1mm]
\subfloat[steady, distorted quadrilaterals]{%
\includegraphics{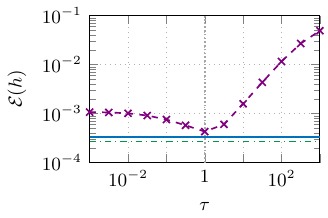}}\hfill
\subfloat[steady, Voronoi]{%
\includegraphics{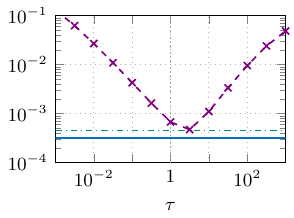}}\hfill
\subfloat[steady, Delaunay dual]{%
\includegraphics{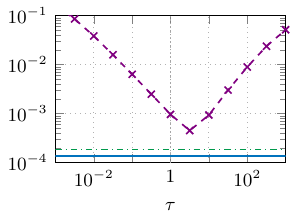}}\\
\subfloat[steady, distorted quadrilaterals]{%
\includegraphics{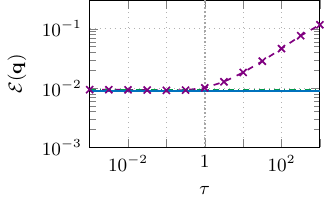}}\hfill
\subfloat[steady, Voronoi]{%
\includegraphics{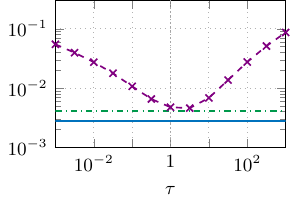}}\hfill
\subfloat[steady, Delaunay dual]{%
\includegraphics{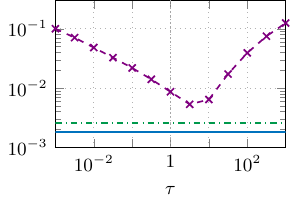}}\\
\subfloat[transient, Voronoi]{%
\includegraphics{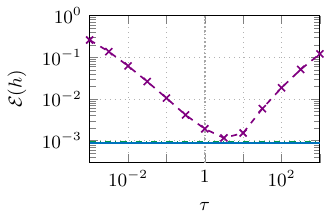}}\hspace{0.05\linewidth}
\subfloat[transient, Delaunay dual]{%
\includegraphics{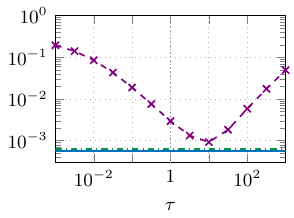}}
\caption{Transient flow in a van Genuchten soil: relative error of VEM against the factor $\tau$ on its stabilisation, $\tau=1$ being the value used in the comparisons, and of SF and HVEM, which have no such term; (a)--(f) steady solution of Eq.~\eqref{eq:mms2}, $n=24$, nodal head and flux at the centroids; (g), (h) transient solution of Eq.~\eqref{eq:mmst}, $n=16$, $\Delta t=1/64$, nodal head at $t=1$.}\label{stab.f}
\end{figure}

A last point concerns the stabilisation of VEM, whose scaling is left to the analyst. To assess its influence, the stabilisation term is multiplied by a factor $\tau$ ranging from $10^{-3}$ to $10^3$. Two problems are considered. The first is the steady manufactured solution of \cite{BergamaschiPutti1999} on the unit square (Fig.~\ref{mms.f}b),
\begin{equation}\label{eq:mms2}
h(x,z)=\cos\Bigl[k\pi\Bigl(\tfrac{x^3}3+\tfrac{z^3}3-\tfrac{x^2}2
-\tfrac{z^2}2\Bigr)\Bigr]-z,\qquad K=K_s\,e^{\alpha h},
\end{equation}
with $k=3$, $\alpha=1$, $K_s=1$, a constant water content and the source term computed from Eq.~\eqref{eq:mms2}. The second is the transient solution of Eq.~\eqref{eq:mmst}. Figure~\ref{stab.f} shows the error of VEM against $\tau$, together with the errors of SF and HVEM, which do not depend on $\tau$ and appear as horizontal lines. The error of VEM varies by more than two orders of magnitude. A small value of $\tau$ leaves the modes outside the polynomial space almost uncontrolled, whereas a large one makes the element too stiff. The optimal value depends on the mesh, on the quantity of interest and on the problem, and $\tau=1$, adopted in all the comparisons of this work, is close to it in every case, so that VEM is not penalised by this choice. Even at its best, VEM remains less accurate than the two stabilisation-free elements, whose errors are close to each other in all cases and require no tuning.

\subsection{Infiltration into a dry layered soil}\label{forsyth.s}

\begin{figure}[t]
\centering
\subfloat[geometry]{\makebox[0.32\linewidth]{%
\includegraphics{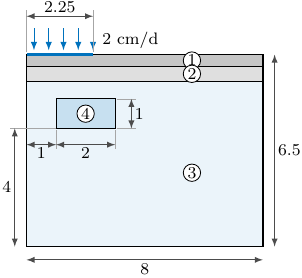}}}\hfill
\subfloat[quadrilaterals, $924$ nodes]{\makebox[0.32\linewidth]{%
  \includegraphics{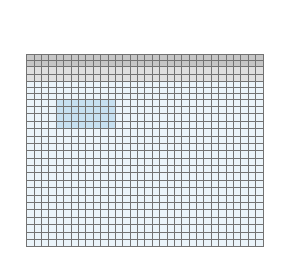}}}\hfill
\subfloat[Delaunay dual, $2284$ nodes]{\makebox[0.32\linewidth]{%
  \includegraphics{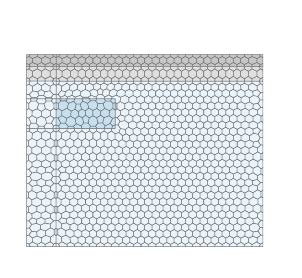}}}
\caption{Infiltration into a dry layered soil: (a) geometry in metres, with the four zones of Table~\ref{forsyth.tab} and the infiltration strip, the other sides being impervious; (b), (c) the rectangular mesh and the Delaunay dual of $0.25$~m, coloured by zone.}\label{forsyth.f}
\end{figure}

Problem~2 of Forsyth, Wu and Pruess \cite{Forsyth1995} concerns a section $8$~m wide and $6.5$~m deep made of four soils, whose geometry and parameters are given in Fig.~\ref{forsyth.f}a and Table~\ref{forsyth.tab}. Two layers lie near the surface, and a block more permeable than the surrounding soil is embedded below them. The soil is initially very dry, with $h_0=-7.34$~m, a head at which the conductivity of zone~$3$ is nearly six orders of magnitude lower than $K_s$. Water infiltrates for $30$~days at $2$~cm/d through the strip $0\le x\le2.25$~m of the surface, while the rest of the boundary is impervious. Since the initial head is uniform, the dry soil also drains slowly under gravity. Three kinds of mesh are used. The first consists of rectangles of about $0.5$, $0.25$ and $0.125$~m, whose lines pass through the interfaces and through the end of the strip (Fig.~\ref{forsyth.f}b). The others are a Delaunay triangulation of $0.25$~m and its dual, with $2284$ nodes (Fig.~\ref{forsyth.f}c). The reference solution is computed with T3 on a Delaunay triangulation of $0.05$~m, with $23819$ nodes.

\begin{table}[t]
\footnotesize
\centering
\caption{Infiltration into a dry layered soil: van Genuchten--Mualem parameters of the four zones~\cite{Forsyth1995}; $K_s$ follows from the permeability with $\rho g/\mu=9.81\cdot10^6$~m$^{-1}$s$^{-1}$, and $S_s=0$.}\label{forsyth.tab}
\begin{tabular}{c c c c c c}
\toprule
zone & $K_s$ [m/d] & $\theta_s$ & $\theta_r$ & $-1/h_\alpha$ [m$^{-1}$] & $n_{\mathrm{vG}}$\\
\midrule
1 & $7.91$ & $0.368$ & $0.102$ & $3.34$ & $1.982$\\
2 & $4.70$ & $0.351$ & $0.098$ & $3.63$ & $1.632$\\
3 & $4.15$ & $0.325$ & $0.086$ & $3.45$ & $1.573$\\
4 & $41.5$ & $0.325$ & $0.086$ & $3.45$ & $1.573$\\
\bottomrule
\end{tabular}
\end{table}

\begin{figure}[t]
\centering
\subfloat[quadrilaterals, $t=10$~d]{\includegraphics[width=0.36\linewidth]{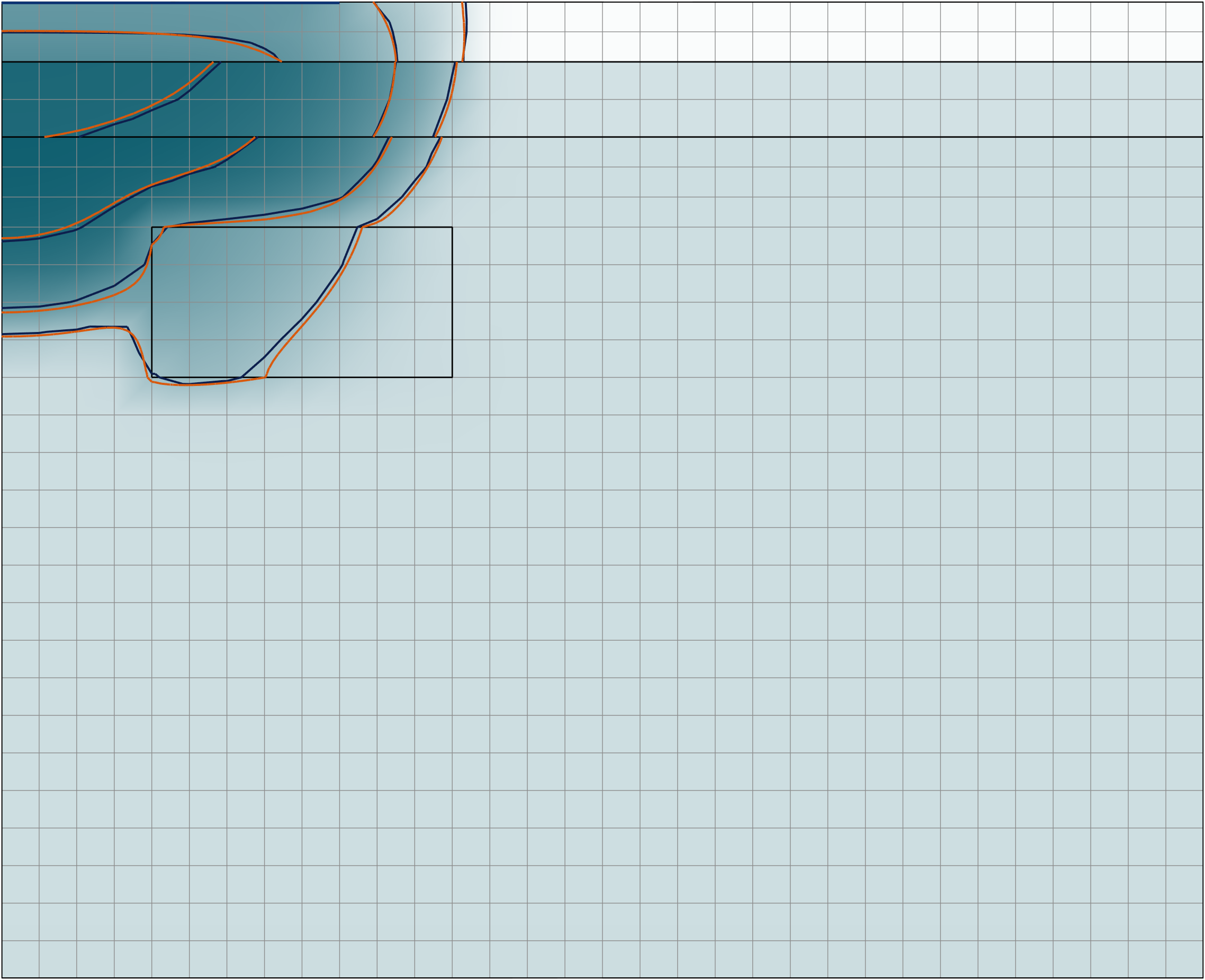}}\hspace{6mm}
\subfloat[quadrilaterals, $t=30$~d]{\includegraphics[width=0.36\linewidth]{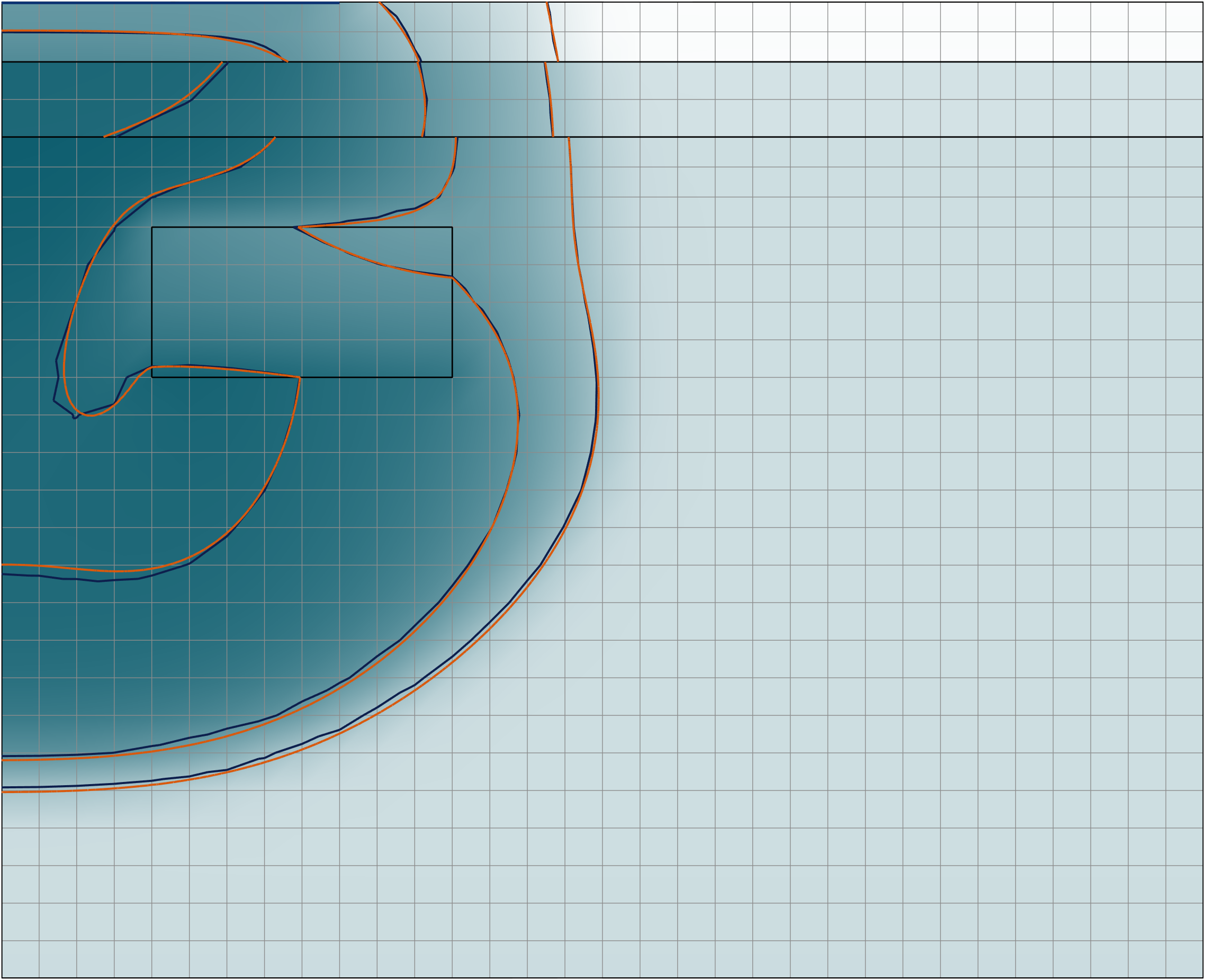}}\\[1mm]
\subfloat[Delaunay dual, $t=10$~d]{\includegraphics[width=0.36\linewidth]{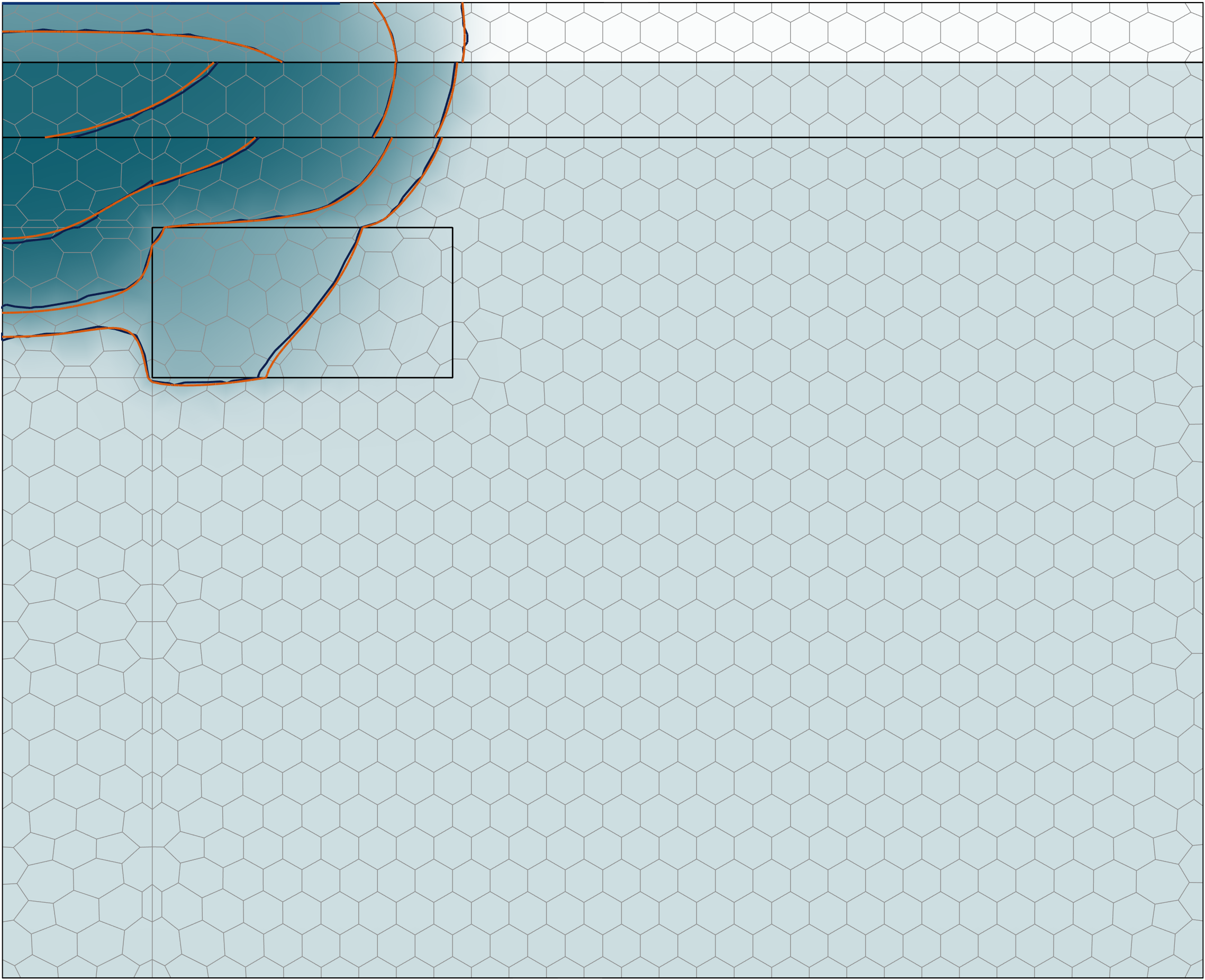}}\hspace{6mm}
\subfloat[Delaunay dual, $t=30$~d]{\includegraphics[width=0.36\linewidth]{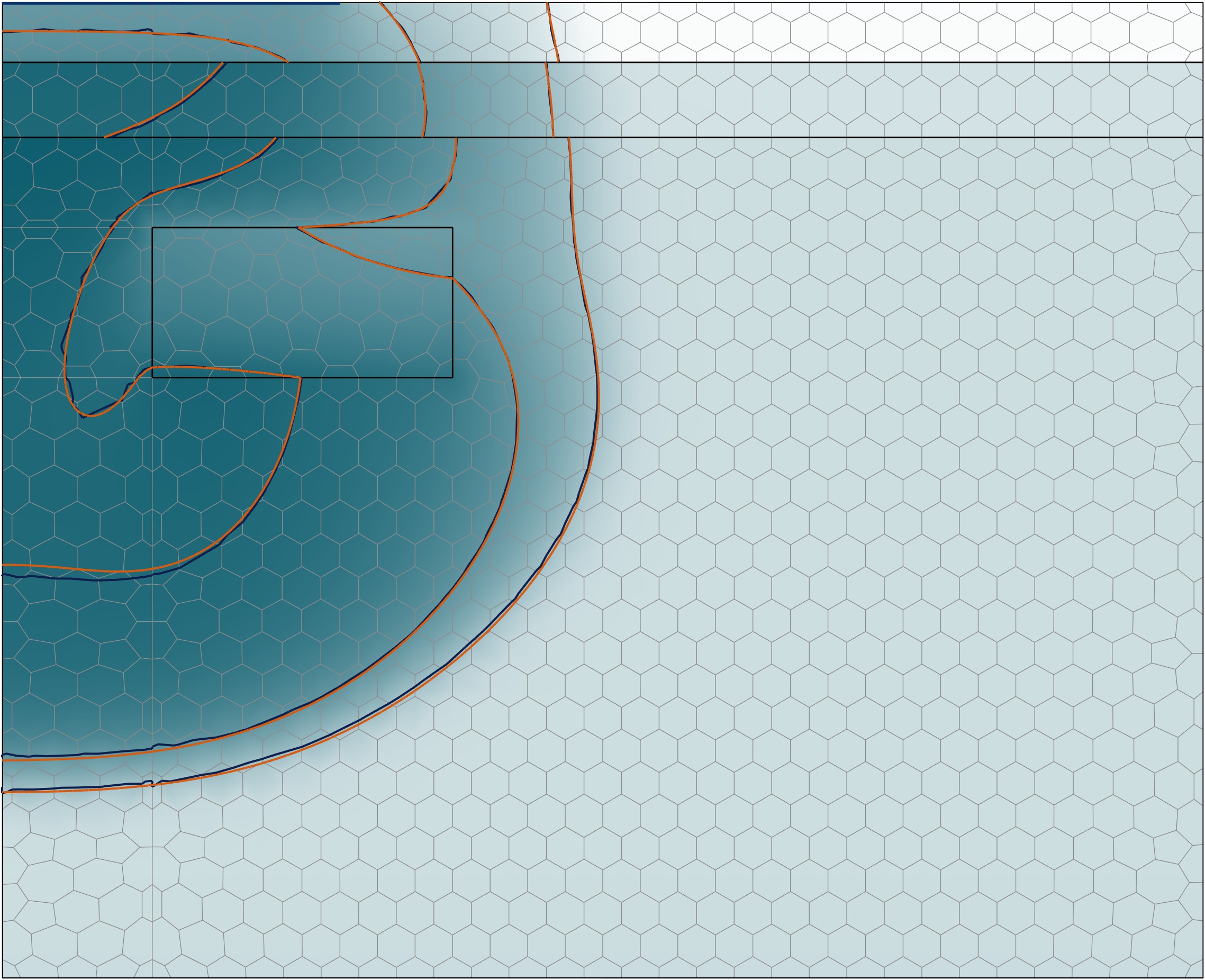}}\\[1mm]
  \includegraphics{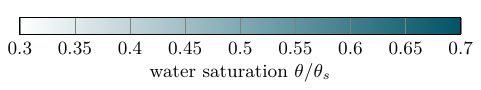}
\caption{Infiltration into a dry layered soil: water saturation computed by HVEM on the two meshes of Fig.~\ref{forsyth.f}, with its level lines $0.35$, $0.45$, $0.55$ and $0.65$ (dark) and those of the reference (orange).}\label{forsythmaps.f}
\end{figure}

\begin{figure}[t]
\centering
{\mbox{%
  \raisebox{2pt}{\includegraphics{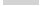}}\;{\footnotesize reference, T3 on $0.05$~m}}\quad\mbox{%
  \raisebox{2pt}{\includegraphics{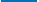}}\;{\footnotesize HVEM, quadrilaterals $0.125$~m}}\quad\mbox{%
  \raisebox{2pt}{\includegraphics{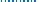}}\;{\footnotesize HVEM, Delaunay dual $0.25$~m}}\quad\mbox{%
  \raisebox{2pt}{\includegraphics{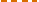}}\;{\footnotesize Abaqus, CPE4P on $0.0625$~m}}\quad}\\[1mm]
\includegraphics{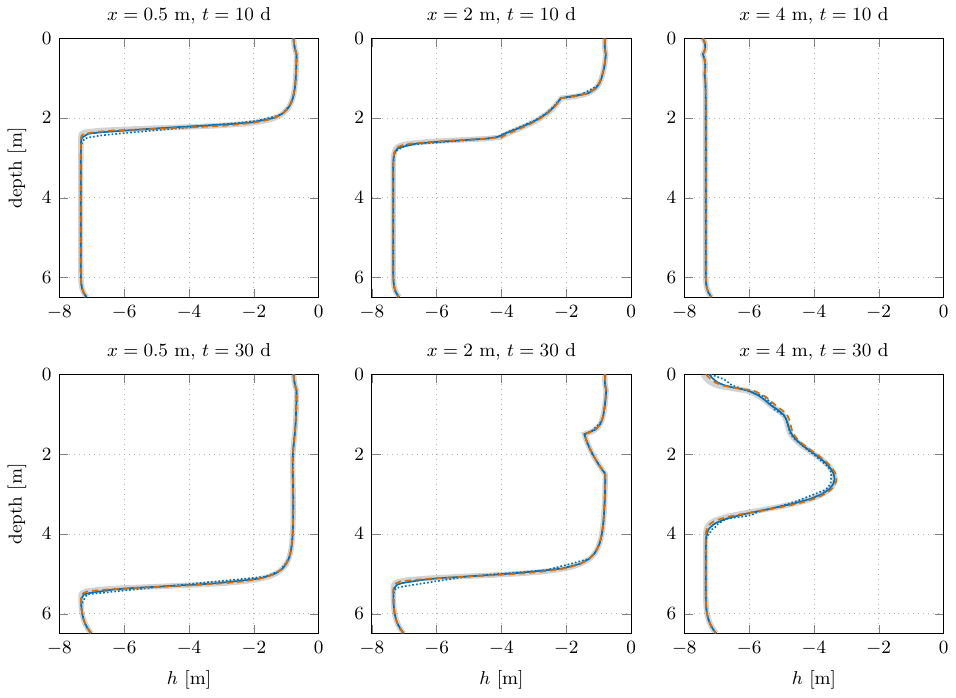}
\caption{Infiltration into a dry layered soil: pressure head along three verticals at $t=10$ and $30$~d, computed by HVEM on two meshes, by the reference and by Abaqus with CPE4P elements of $0.0625$~m; the vertical at $x=2$~m crosses the permeable block.}\label{forsythabq.f}
\end{figure}

Figure~\ref{forsythmaps.f} shows the water saturation computed by HVEM on the two meshes of $0.25$~m after $10$ and $30$~days, together with the level lines of the reference. The front advances faster through the permeable block, where the saturation remains lower because the block carries the same flux with less water, and spreads below it. On both meshes the level lines of HVEM follow those of the reference within a fraction of an element. The agreement is confirmed in Fig.~\ref{forsythabq.f}, which compares the head along three verticals with the independent solution of Abaqus. The three solutions are within a few centimetres of head of each other outside the front and across the block.

\begin{table}[t]
\footnotesize
\centering
\caption{Infiltration into a dry layered soil: cost, water budget, lowest head relative to the initial one and relative difference from the reference at $t=30$~d. The imbalance is in m$^2$, against an infiltrated volume of $1.35$~m$^2$; the reference gives $\min(h-h_0)=-0.28$~m. In the rows marked conservative, SF and VEM write the storage as the increment of the water content.}\label{forsythres.tab}
\begin{tabular}{l r r r r c c c}
\toprule
 & nodes & steps & iterations & CPU [s] & imbalance [m$^2$] & $\min(h-h_0)$ [m] & $\mathcal E(h)$ [\%]\\
\midrule
\multicolumn{8}{l}{\emph{quadrilaterals, 0.5 m}}\\
\quad HVEM & 270 & 142 & 684 & 0.7 & $5\cdot10^{-15}$ & $-0.20$ & $6.0$\\
\quad SF & 270 & 219 & 1070 & 3.2 & $-0.050$ & $-0.20$ & $6.2$\\
\quad VEM & 270 & 229 & 1122 & 3.2 & $-0.047$ & $-1.40$ & $6.3$\\
\quad Q4 & 270 & 150 & 736 & 0.7 & $-5\cdot10^{-15}$ & $-0.21$ & $5.6$\\
\multicolumn{8}{l}{\emph{quadrilaterals, 0.25 m}}\\
\quad HVEM & 924 & 203 & 994 & 3.2 & $-4\cdot10^{-15}$ & $-0.26$ & $1.8$\\
\quad SF & 924 & 317 & 1570 & 16 & $-0.042$ & $-0.26$ & $5.7$\\
\quad VEM & 924 & 329 & 1619 & 17 & $-0.041$ & $-0.91$ & $5.2$\\
\quad Q4 & 924 & 224 & 1098 & 4.3 & $-3\cdot10^{-15}$ & $-0.26$ & $1.7$\\
\multicolumn{8}{l}{\emph{quadrilaterals, 0.125 m}}\\
\quad HVEM & 3510 & 267 & 1315 & 18 & $2\cdot10^{-15}$ & $-0.27$ & $0.6$\\
\quad SF & 3510 & 424 & 2096 & 86 & $-0.038$ & $-0.27$ & $6.0$\\
\quad SF, conservative & 3510 & 282 & 1389 & 84 & $-2\cdot10^{-14}$ & $-0.27$ & $0.5$\\
\quad VEM & 3510 & 451 & 2242 & 93 & $-0.035$ & $-0.46$ & $5.8$\\
\quad VEM, conservative & 3510 & 256 & 1263 & 65 & $4\cdot10^{-15}$ & $-0.56$ & $0.9$\\
\quad Q4 & 3510 & 284 & 1399 & 23 & $-2\cdot10^{-14}$ & $-0.27$ & $0.5$\\
\multicolumn{8}{l}{\emph{Delaunay dual, 0.25 m}}\\
\quad HVEM & 2284 & 271 & 1306 & 13 & $4\cdot10^{-15}$ & $-0.27$ & $2.1$\\
\quad SF & 2284 & \multicolumn{6}{l}{stops at $t=0.19$~d}\\
\quad SF, conservative & 2284 & 290 & 1433 & 50 & $4\cdot10^{-15}$ & $-0.27$ & $2.0$\\
\quad VEM & 2284 & 351 & 1743 & 41 & $-0.041$ & $-1.33$ & $5.1$\\
\quad VEM, conservative & 2284 & 213 & 1055 & 35 & $2\cdot10^{-15}$ & $-1.18$ & $2.5$\\
\multicolumn{8}{l}{\emph{Delaunay triangulation, 0.25 m}}\\
\quad HVEM & 971 & 225 & 1095 & 5 & $5\cdot10^{-15}$ & $-0.26$ & $2.1$\\
\quad T3 & 971 & 226 & 1109 & 6 & $-1\cdot10^{-15}$ & $-0.26$ & $1.9$\\
\quad MHFE & 971 & 2529 & 12744 & 209 & $-7\cdot10^{-15}$ & $-0.27$ & $17.2$\\
\bottomrule
\end{tabular}
\end{table}

Table~\ref{forsythres.tab} collects, for all the meshes, the cost, the water budget, the lowest head reached and the difference from the reference at the end of the analysis. On the rectangular meshes HVEM converges to the reference as Q4 does, and both close the water budget to round-off. SF and VEM do not improve with refinement and remain a few percent away from the reference on every mesh. The reason is that they lose part of the infiltrated water, as their budget shows. This is the well known error of the capacitive form of the storage \cite{Celia1990A}, with which the two elements are published, and it grows with the time step but does not vanish under mesh refinement. VEM shows a further defect. Ahead of the front its head falls well below the initial value, far beyond the slow drainage of the dry soil, which the other schemes reproduce as the reference does.

On the triangulation HVEM is as accurate as T3, and on the dual mesh it is more accurate than VEM. SF does not complete the analysis on the dual mesh and stops at an early time with a positive head on the infiltration strip. This failure can be explained by looking at the capacity lumped on the nodes. In SF this is the capacity of the projected head at the quadrature points around the node, which at the front are drier than the node itself. At the edge of the strip it falls to a small fraction of the capacity of the node, which therefore saturates in a fraction of the time. HVEM, like the lumped primal element~\cite{Celia1990A}, evaluates the capacity at the node. MHFE, whose element heads are averaged to the nodes for the comparison, is farther from the reference and needs many more steps. As regards the cost, HVEM is several times faster than SF and VEM, which also require more steps, and slightly faster than Q4 and T3, whose assembly evaluates the conductivity at the Gauss points.

To confirm that these differences come from the storage and not from the elements, SF and VEM have been run once more with the storage written in the conservative form, as the increment of the water content, lumped on the nodes for SF and evaluated at the quadrature points for VEM, the rest of each element being unchanged (rows marked conservative in Table~\ref{forsythres.tab}). Both then close the water budget to round-off and converge to the reference as HVEM does, SF with the same accuracy and VEM with a somewhat larger difference, and SF completes the analysis on the dual mesh. The head of VEM still falls ahead of the front, an effect of its consistent capacity matrix rather than of the form of the storage. The cost is unchanged, since the conductivity and the capacity are still evaluated at the quadrature points of every element.

\subsection{Seepage face and local refinement}\label{aquifer.s}

\begin{figure}[t]
\centering
\includegraphics{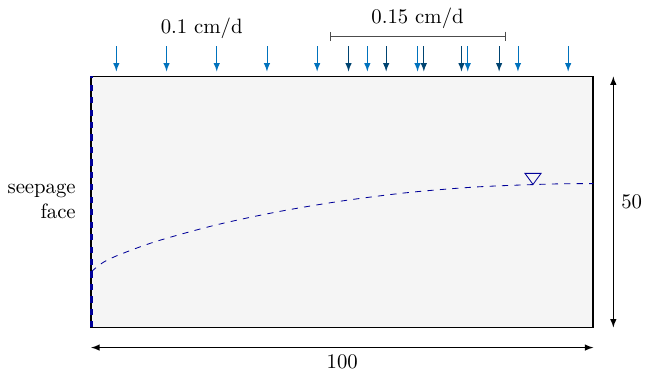}
\caption{Seepage face and local refinement: geometry in centimetres, rainfall and seepage face; the dashed line is the computed water table.}\label{aqgeo.f}
\end{figure}

A block of $100\times50$~cm under steady rainfall is analysed \cite{BergamaschiPutti1999}. The water is drained through its left side, which is a potential seepage face, as shown in Fig.~\ref{aqgeo.f}. The soil follows the model in Eq.~\eqref{eq:bp} with the logarithmic relative conductivity, with $\phi=0.30$, $S_{wr}=0.01$, $\alpha=0.015$~cm$^{-1}$, $\beta=2$, $\gamma=3$, $a=2$, $b=3.5$, $h_s=-10$~cm, $K_s=1$~cm/d and $S_s=0$. A rain of $0.1$~cm/d falls on the top side, increased to $0.15$~cm/d on the strip $47.5\le x\le82.5$~cm, while the base and the right side are impervious. The steady state is computed in a single implicit step.

\begin{figure}[t]
\centering
\subfloat[uniform, $55$ elements]{%
  \makebox[0.32\linewidth]{\includegraphics{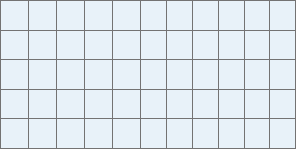}}}\hfill
\subfloat[graded, $127$ elements]{%
  \makebox[0.32\linewidth]{\includegraphics{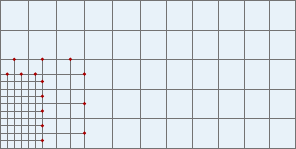}}}\hfill
\subfloat[Delaunay dual, $242$ elements]{%
  \makebox[0.32\linewidth]{\includegraphics{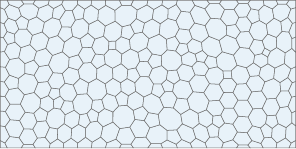}}}\\[1mm]
\subfloat[uniform]{\includegraphics[width=0.32\linewidth]{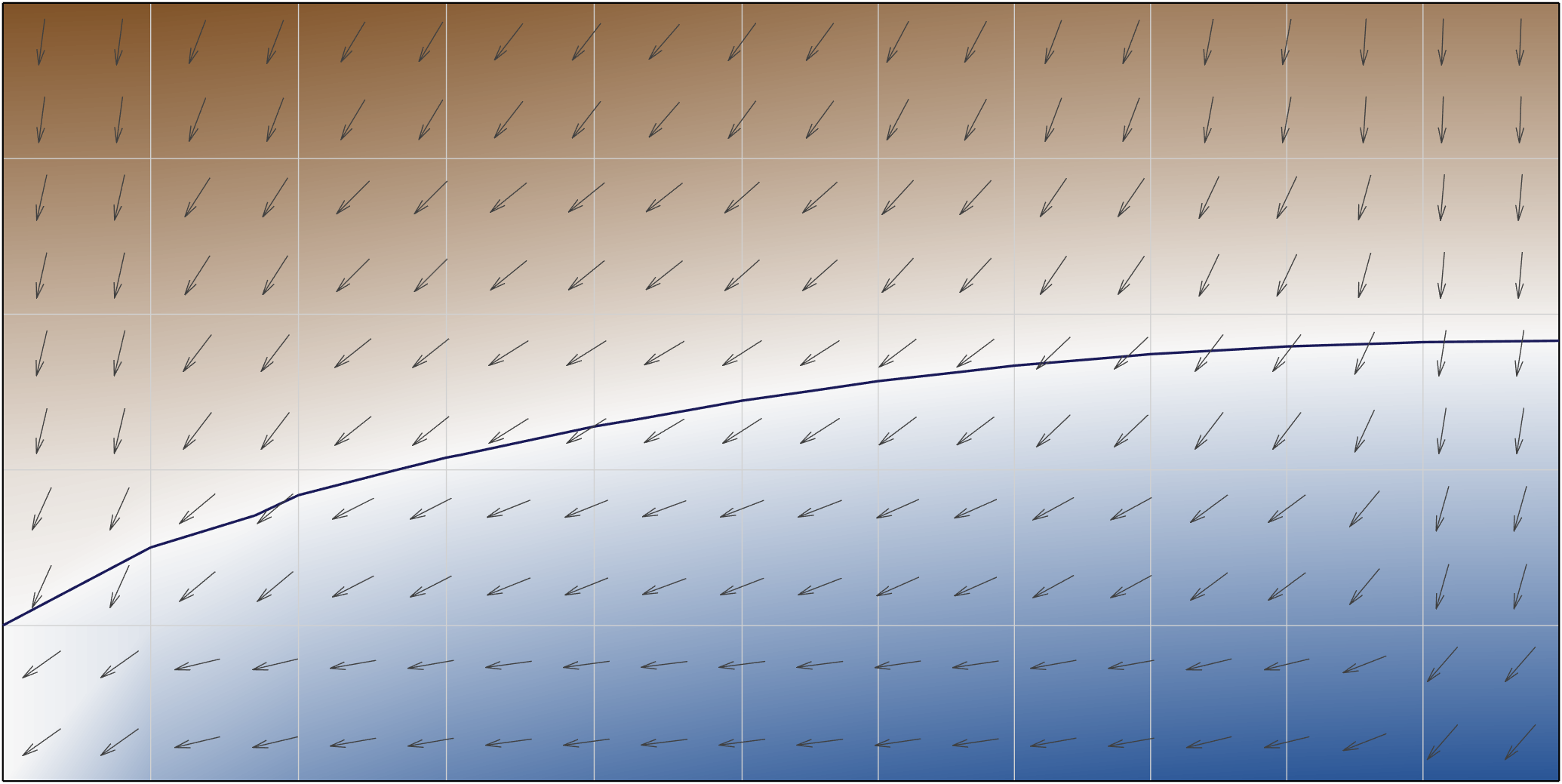}}\hfill
\subfloat[graded]{\includegraphics[width=0.32\linewidth]{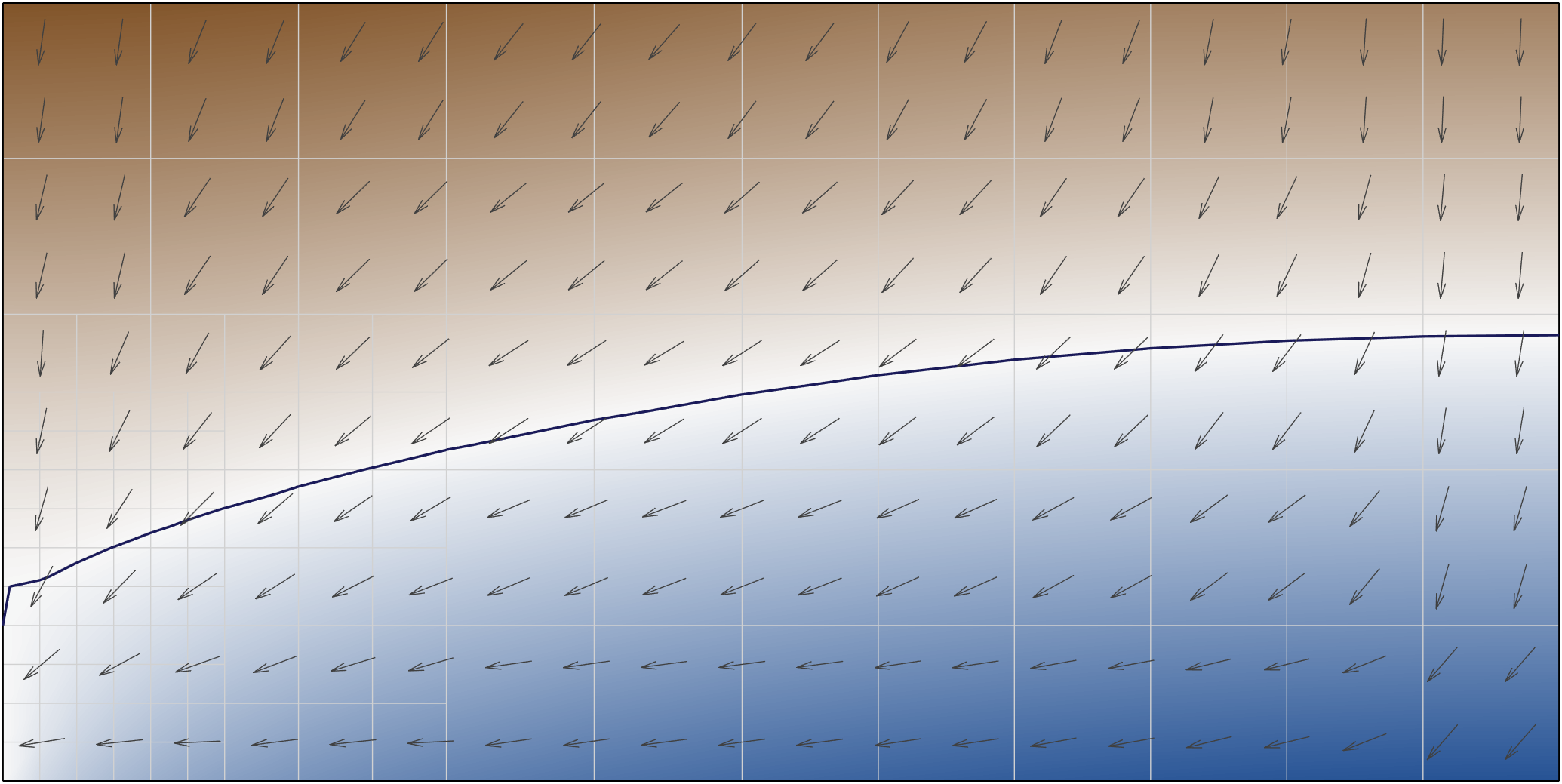}}\hfill
\subfloat[Delaunay dual]{\includegraphics[width=0.32\linewidth]{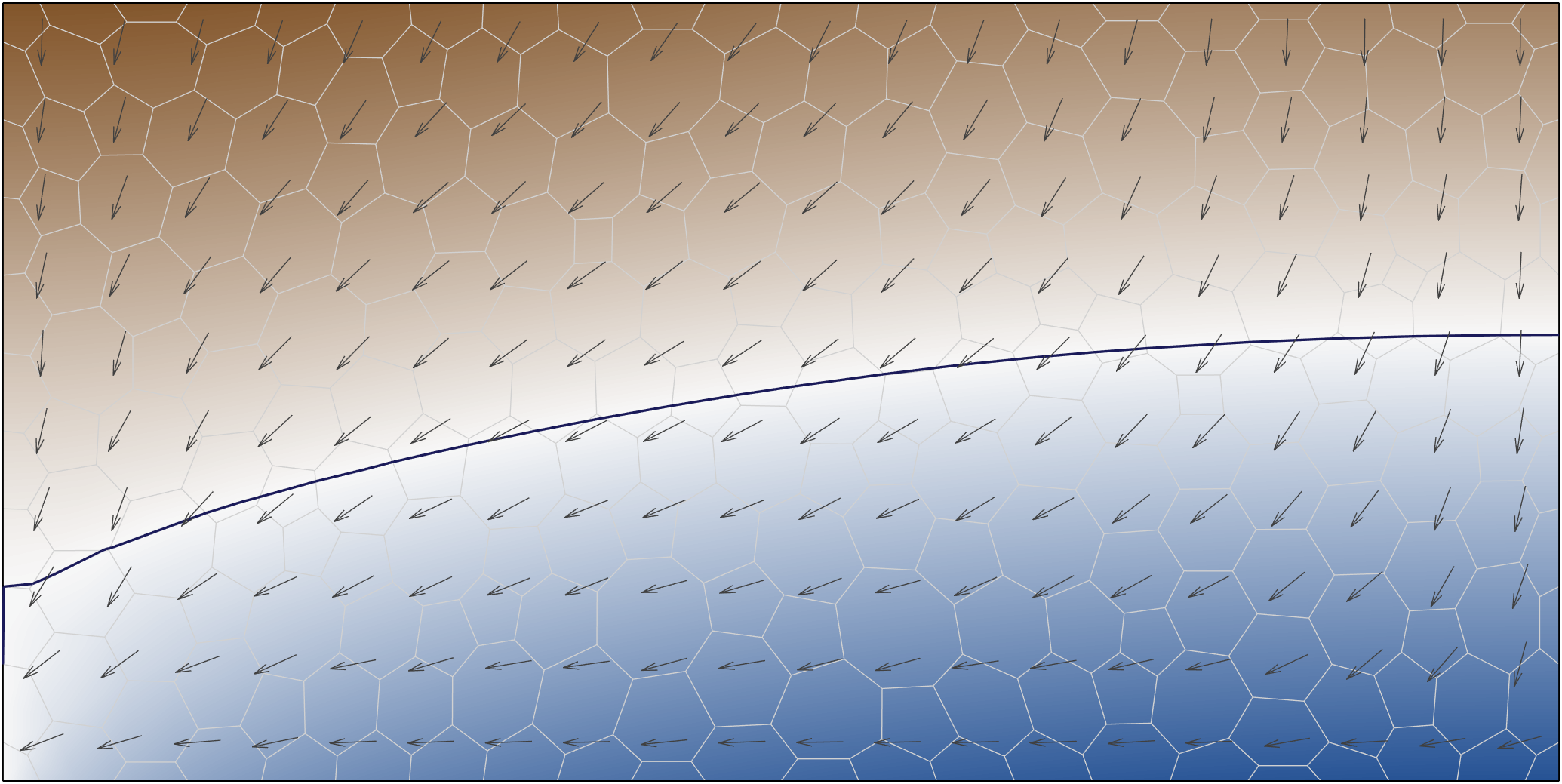}}\\[1mm]
  \includegraphics{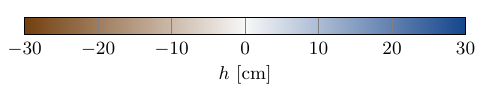}
\caption{Seepage face and local refinement: (a) the coarse uniform mesh, (b) the same refined by two levels towards the ditch, with the hanging nodes in red, (c) a coarse Delaunay dual; (d)--(f) head computed by HVEM on each of them, with the water table and the direction of the flux.}\label{aqhang.f}
\end{figure}

\begin{table}[t]
\footnotesize
\centering
\caption{Seepage face and local refinement: elevation of the exit point, where the head interpolated along the side vanishes, on uniform, graded and dual meshes and on the reference. In the first three groups the graded mesh is refined from the uniform mesh above it and has at the ditch the element size of the uniform mesh below it; Q4 accepts neither hanging nodes nor polygons.}\label{aqhang.tab}
\begin{tabular}{l rr cccc}
\toprule
 & & & \multicolumn{4}{c}{$z_{\mathrm{exit}}$ [cm]}\\
\cmidrule(l{2pt}r{2pt}){4-7}
mesh & elements & dofs & HVEM & SF & VEM & Q4\\
\midrule
uniform & 55 & 72 & $10.043$ & $10.043$ & $10.043$ & $10.043$\\
graded & 127 & 158 & $11.252$ & $11.268$ & $12.483$ & --\\
uniform & 800 & 861 & $11.130$ & $11.131$ & $12.472$ & $11.131$\\
\midrule
uniform & 210 & 242 & $10.090$ & $10.090$ & $10.095$ & $10.090$\\
graded & 486 & 546 & $11.345$ & $11.345$ & $11.350$ & --\\
uniform & 3200 & 3321 & $11.345$ & $11.345$ & $11.350$ & $11.345$\\
\midrule
uniform & 800 & 861 & $11.130$ & $11.131$ & $12.472$ & $11.131$\\
graded & 1904 & 2021 & $11.851$ & $11.851$ & $11.872$ & --\\
uniform & 12800 & 13041 & $11.852$ & $11.852$ & $11.873$ & $11.852$\\
\midrule
dual & 242 & 544 & $12.073$ & $12.503$ & $10.063$ & --\\
dual & 861 & 1840 & $11.379$ & $12.446$ & $12.505$ & --\\
dual & 3321 & 6880 & $11.888$ & $11.888$ & $11.879$ & --\\
\midrule
reference, Q4 & 51200 & 51681 & \multicolumn{4}{c}{$11.654$}\\
\bottomrule
\end{tabular}
\end{table}

\begin{figure}[t]
\centering
{\mbox{%
  \raisebox{2pt}{\includegraphics{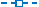}}\;{\footnotesize uniform meshes}}\quad\mbox{%
  \raisebox{2pt}{\includegraphics{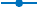}}\;{\footnotesize graded meshes}}\quad}\\[1mm]
\subfloat[nodal head]{\includegraphics{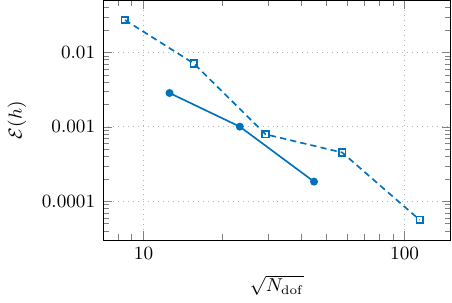}}\hfill
\subfloat[exit point]{\includegraphics{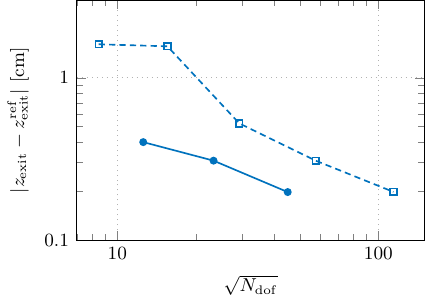}}
\caption{Seepage face and local refinement: (a) relative error of the head of HVEM and (b) error of the exit elevation, on the uniform and graded meshes of Table~\ref{aqhang.tab}; below about $0.2$~cm the exit elevation is not resolved by the reference.}\label{aqconv.f}
\end{figure}

The only strongly localised feature of the solution is the exit point of the seepage face, which makes it a natural candidate for local refinement. As discussed in \S\ref{sec:hang}, a virtual element needs no special treatment for this purpose. Figure~\ref{aqhang.f} shows the coarse uniform mesh, the mesh obtained from it by two levels of refinement towards the ditch and a coarse Delaunay dual, together with the heads computed by HVEM on each of them. Table~\ref{aqhang.tab} reports the elevation of the exit point given by all the elements on uniform, graded and dual meshes, with the reference in the last row. The exit point is located where the head, interpolated linearly on the triangulation of the nodes, changes sign along the side, so that it is not bound to the nodes. Each graded mesh places the exit point almost where the uniform mesh with the same element size at the ditch does, with several times fewer dofs. Since the free boundary moves across the mesh in discrete jumps, only the meshes of the same group can be compared with each other. The convergence of the head and of the exit elevation is shown in Fig.~\ref{aqconv.f} for the uniform and the graded meshes. At equal dofs the graded meshes are more accurate in the head, whereas the exit elevation is not resolved by the reference below about $0.2$~cm. Abaqus, on the uniform mesh of $3321$ dofs, places the exit point within a millimetre of the reference. SF and VEM agree with HVEM on the finer meshes. On the coarser ones, however, the exit point of VEM moves to a neighbouring node, and so does that of SF on the dual meshes.

\subsection{Rainfall infiltration in a homogeneous slope}\label{lu.s}

\begin{figure}[t]
\centering
\subfloat[geometry]{\makebox[0.47\linewidth]{%
\includegraphics{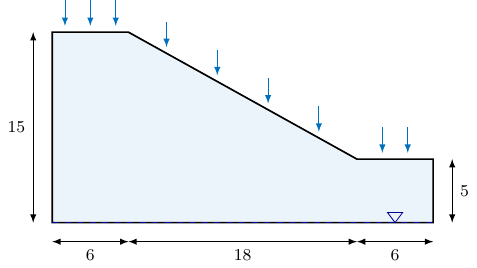}}}\hfill
\subfloat[daily rainfall]{%
\includegraphics{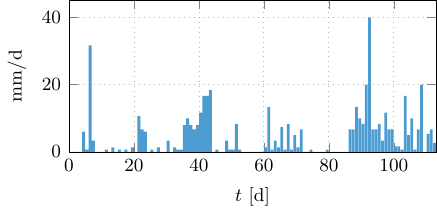}}
\caption{Rainfall infiltration in a homogeneous slope: (a) geometry in metres, with the soil of~\cite{Lu2012}, the water table at the base and the rain on the whole upper boundary; (b) daily rainfall, $480$~mm in all.}\label{lu.f}
\end{figure}

\begin{figure}[t]
\centering
\subfloat[triangles, $263$ nodes]{\makebox[0.48\linewidth]{%
  \includegraphics{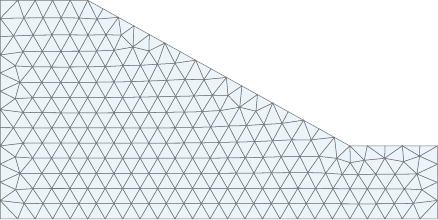}}}\hfill
\subfloat[boundary-fitted quadrilaterals, $338$ nodes]{\makebox[0.48\linewidth]{%
  \includegraphics{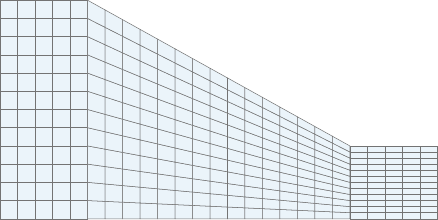}}}\\[1mm]
\subfloat[cut grid, $250$ nodes]{\makebox[0.48\linewidth]{%
  \includegraphics{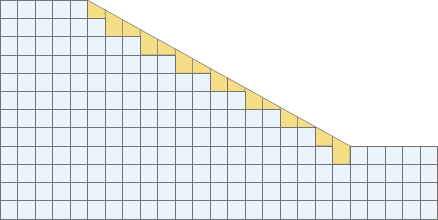}}}\hfill
\subfloat[Delaunay dual, $592$ nodes]{\makebox[0.48\linewidth]{%
  \includegraphics{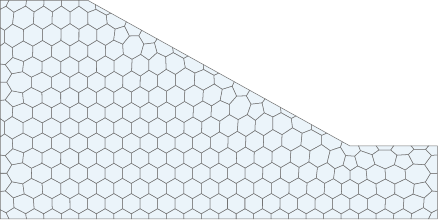}}}
\caption{Rainfall infiltration in a homogeneous slope: the four kinds of mesh at the size of $1.2$~m; in the cut grid only the elements crossed by the face, in yellow, become polygons, with up to six vertices.}\label{lumesh.f}
\end{figure}

A slope is herein analysed. The soil is homogeneous and the geometry, taken from Lu et al.~\cite{Lu2012}, is shown in Fig.~\ref{lu.f}. The face rises by $10$~m over $18$~m and meets the platform at the toe with a re-entrant corner. The soil is a silt with $\theta_r=0.034$, $\theta_s=0.46$, $n_{\mathrm{vG}}=1.37$, $h_\alpha=-0.625$~m, $K_s=1.39\cdot10^{-6}$~m/s and $S_s=10^{-3}$~m$^{-1}$. The water table lies at the base and the sides are impervious. Starting from hydrostatic equilibrium, rain falls on the whole upper boundary for $113$ days as a sequence of daily rainfalls of irregular intensity, $480$~mm in total, and the analysis proceeds with an adaptive step of up to six hours. The section is discretised by the four kinds of mesh shown in Fig.~\ref{lumesh.f}, namely a Delaunay triangulation of $1.2$~m and its dual, boundary-fitted quadrilaterals of the same size, as in a classical analysis, and cut grids of $1.2$, $0.6$ and $0.3$~m. The reference is T3 on a triangulation of $0.2$~m with $8791$ nodes. It is only $1.5$ times finer than the finest cut grid, so that on that grid the differences measure the resolution of the reference as much as the error of the schemes.

\begin{figure}[t]
\centering
\subfloat[cut grid, $0.6$~m, $t=45$~d]{\includegraphics[width=0.47\linewidth]{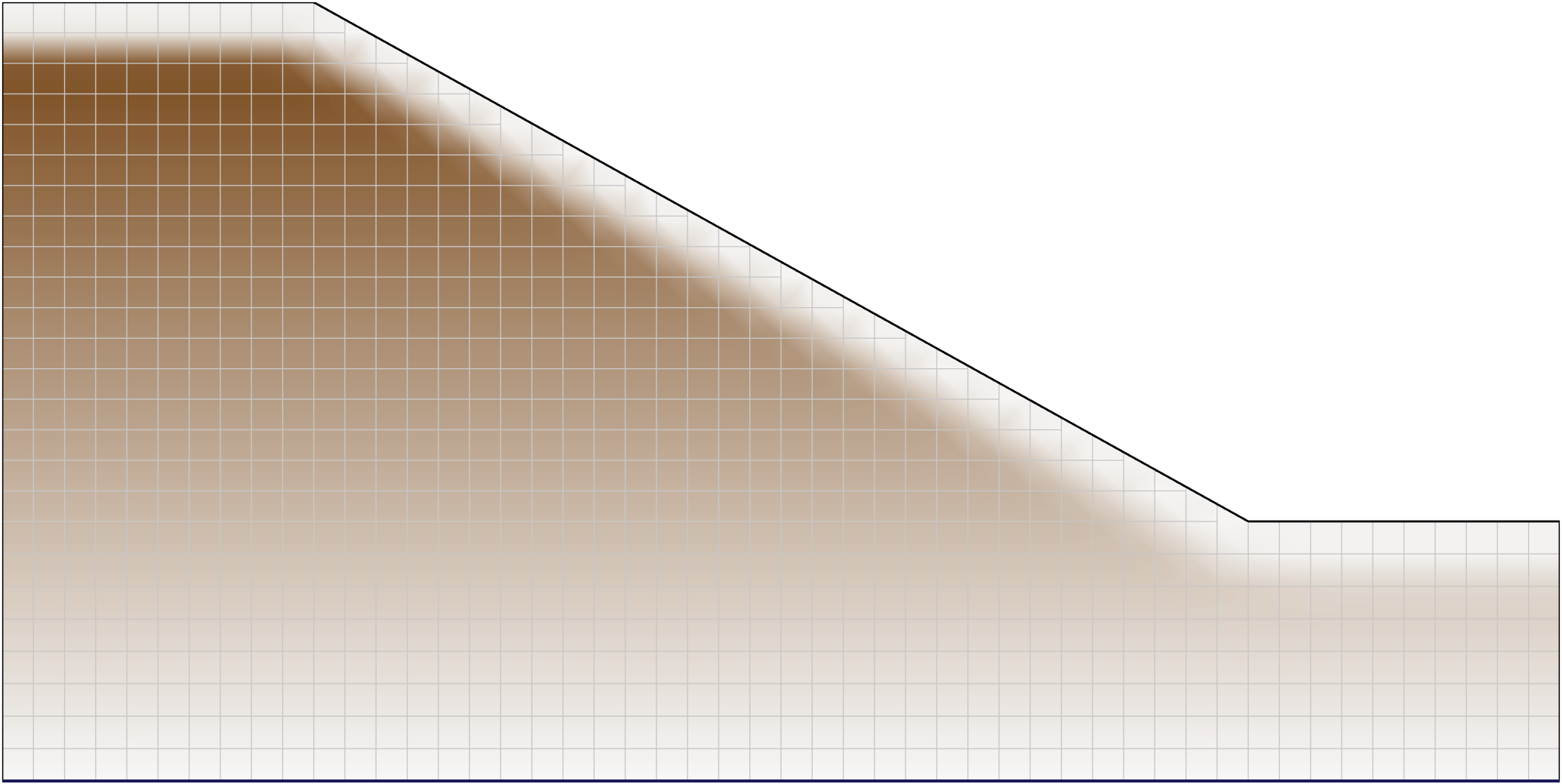}}\hfill
\subfloat[cut grid, $0.6$~m, $t=113$~d]{\includegraphics[width=0.47\linewidth]{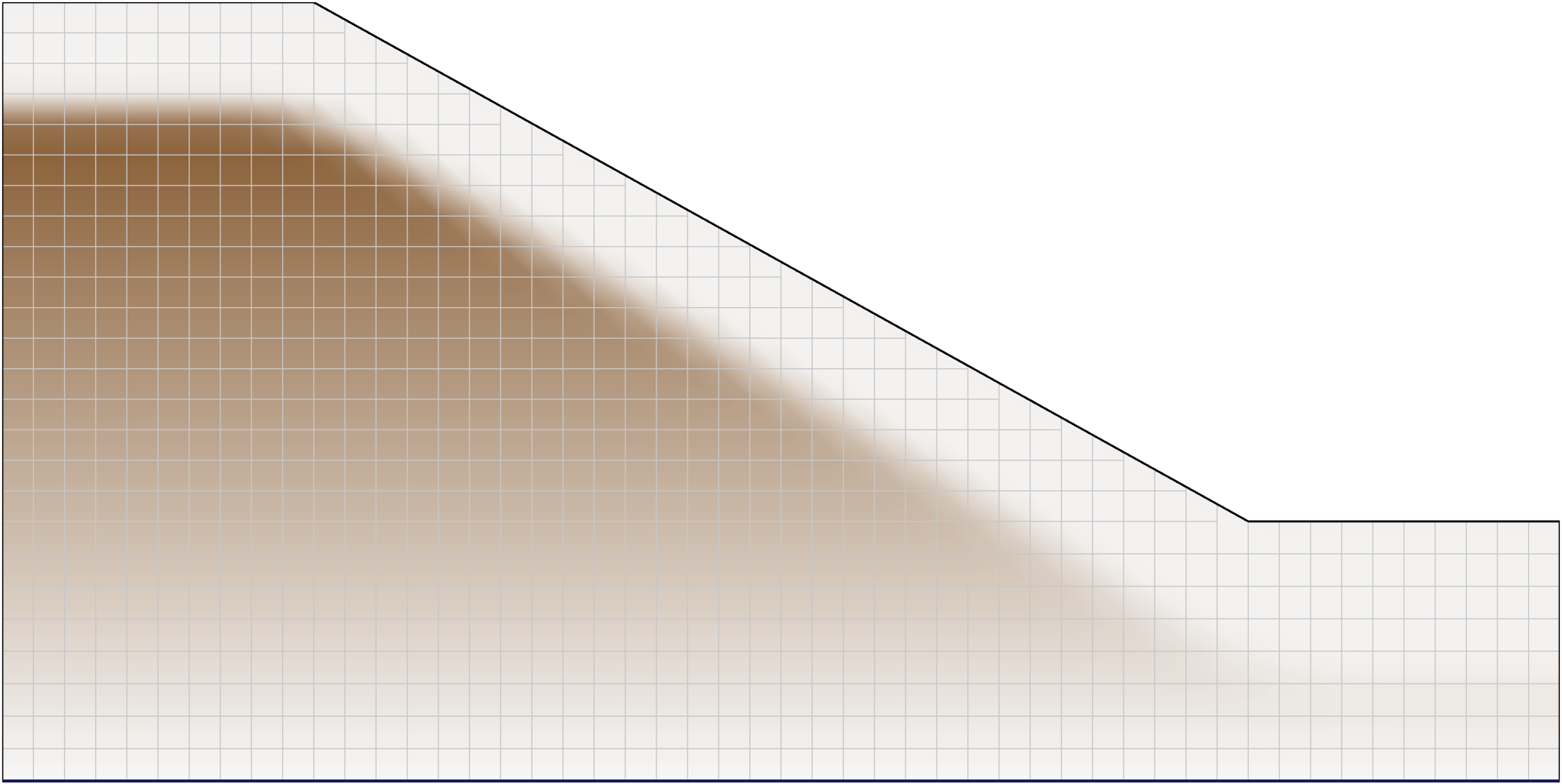}}\\[1mm]
\subfloat[Delaunay dual, $1.2$~m, $t=45$~d]{\includegraphics[width=0.47\linewidth]{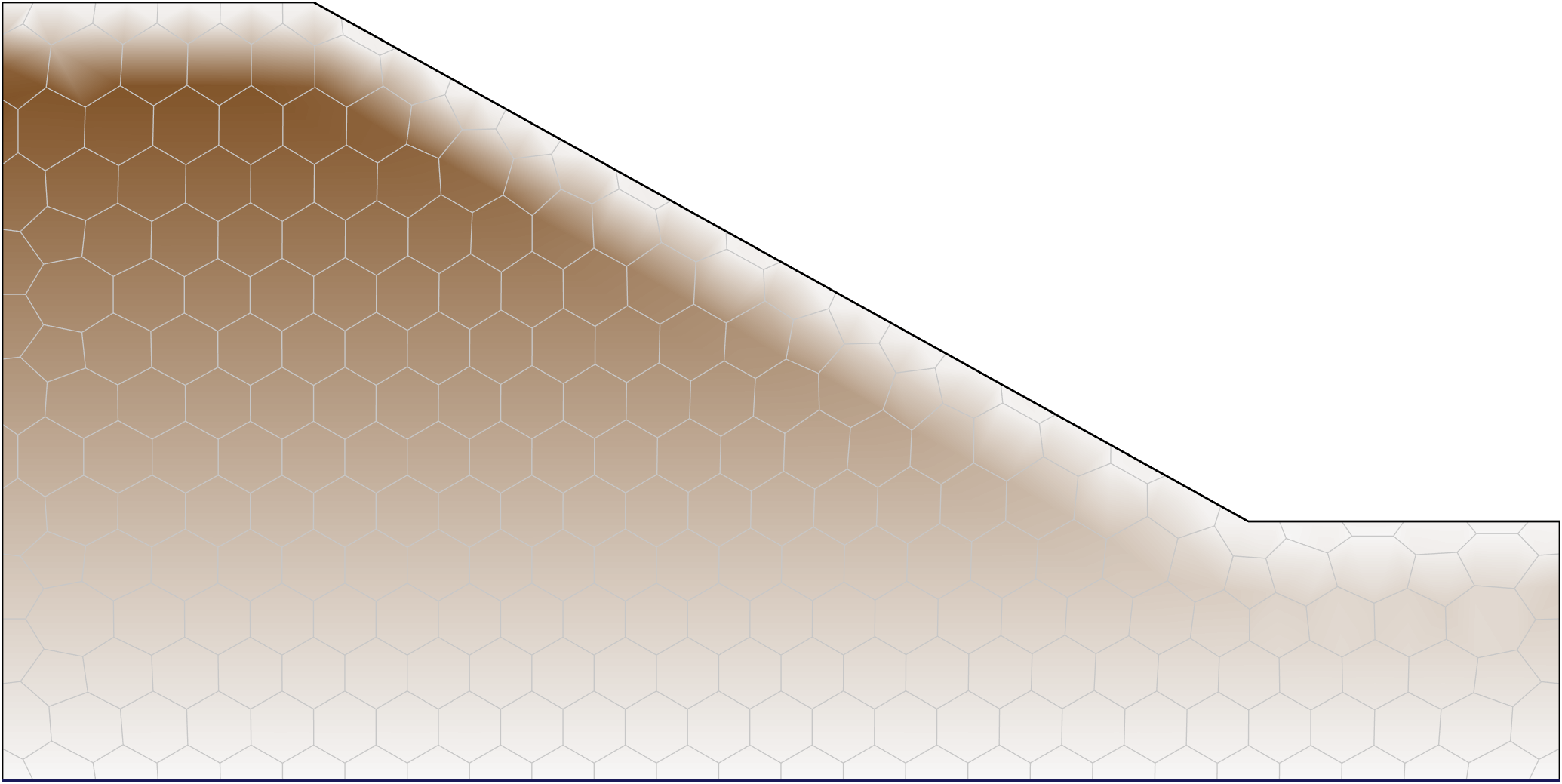}}\hfill
\subfloat[Delaunay dual, $1.2$~m, $t=113$~d]{\includegraphics[width=0.47\linewidth]{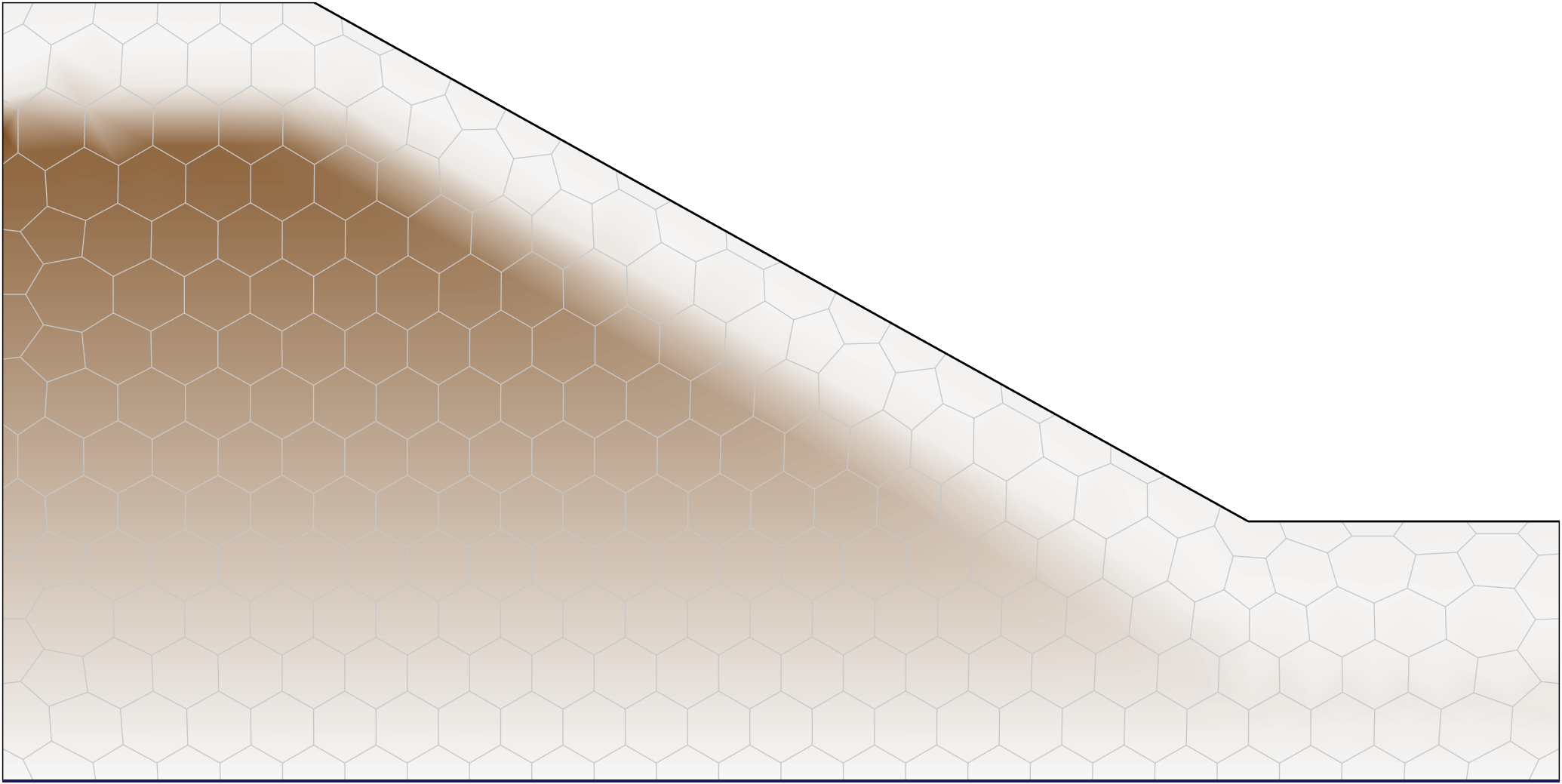}}\\[1mm]
  \includegraphics{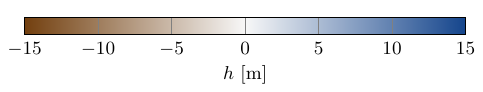}
\caption{Rainfall infiltration in a homogeneous slope: pressure head computed by HVEM on the cut grid of $0.6$~m and on the dual mesh, with the line $h=0$; the initial state is hydrostatic.}\label{lumaps.f}
\end{figure}

\begin{figure}[t]
\centering
{\mbox{%
  \raisebox{2pt}{\includegraphics{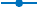}}\;{\footnotesize HVEM}}\quad\mbox{%
  \raisebox{2pt}{\includegraphics{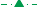}}\;{\footnotesize SF}}\quad\mbox{%
  \raisebox{2pt}{\includegraphics{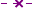}}\;{\footnotesize VEM}}\quad}\\[1mm]
\subfloat[]{\includegraphics{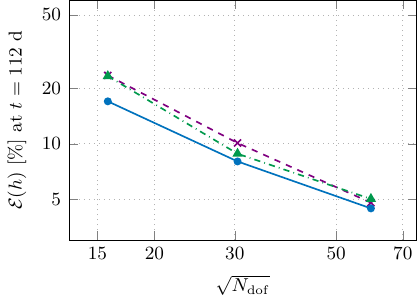}}\hfill
\subfloat[]{\includegraphics{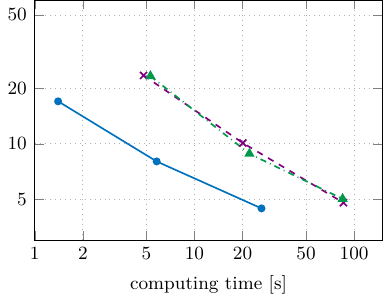}}
\caption{Rainfall infiltration in a homogeneous slope: relative difference of the head from the reference after $112$ days of rain on the cut grids of $1.2$, $0.6$ and $0.3$~m, against (a) $\sqrt{N_{\mathrm{dof}}}$ and (b) the computing time.}\label{luconv.f}
\end{figure}

\begin{table}[t]
\footnotesize
\centering
\caption{Rainfall infiltration in a homogeneous slope: cost, water budget, lowest head relative to the initial one and relative difference from the reference after $112$ days of rain. The imbalance is in m$^2$, against an inflow of $14.5$~m$^2$; the reference takes $487$ steps and $1755$ iterations.}\label{lu.tab}
\begin{tabular}{l r r r r c c c}
\toprule
 & nodes & steps & iterations & CPU [s] & imbalance [m$^2$] & $\min(h-h_0)$ [m] & $\mathcal E(h)$ [\%]\\
\midrule
\multicolumn{8}{l}{\emph{cut grid, 0.6 m}}\\
\quad HVEM & 925 & 487 & 1695 & 6 & $2\cdot10^{-15}$ & $0$ & $8.0$\\
\quad SF & 925 & 487 & 2000 & 22 & $-0.391$ & $0$ & $8.8$\\
\quad VEM & 925 & 487 & 1762 & 20 & $-0.501$ & $-2.62$ & $10.1$\\
\multicolumn{8}{l}{\emph{Delaunay dual, 1.2 m}}\\
\quad HVEM & 592 & 487 & 1702 & 3 & $-4\cdot10^{-14}$ & $-1.33$ & $13.8$\\
\quad SF & 592 & 487 & 1972 & 11 & $-0.332$ & $-4.30$ & $15.6$\\
\quad VEM & 592 & 487 & 1768 & 11 & $-0.448$ & $-4.20$ & $20.2$\\
\multicolumn{8}{l}{\emph{boundary-fitted quadrilaterals, 1.2 m}}\\
\quad HVEM & 338 & 487 & 1670 & 2 & $-1\cdot10^{-14}$ & $0$ & $15.2$\\
\quad SF & 338 & 488 & 1791 & 7 & $-0.297$ & $0$ & $20.2$\\
\quad VEM & 338 & 487 & 1751 & 7 & $-0.462$ & $-3.71$ & $19.5$\\
\quad Q4 & 338 & 487 & 1857 & 3 & $-4\cdot10^{-14}$ & $0$ & $13.9$\\
\multicolumn{8}{l}{\emph{Delaunay triangulation, 1.2 m}}\\
\quad T3 & 263 & 554 & 2262 & 6 & $9\cdot10^{-9}$ & $0$ & $14.4$\\
\quad MHFE & 263 & 665 & 3274 & 6 & $-3\cdot10^{-10}$ & $-0.51$ & $26.4$\\
\bottomrule
\end{tabular}
\end{table}

\begin{figure}[t]
\centering
{\mbox{%
  \raisebox{2pt}{\includegraphics{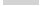}}\;{\footnotesize reference, T3 on $0.2$~m}}\quad\mbox{%
  \raisebox{2pt}{\includegraphics{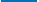}}\;{\footnotesize HVEM, cut grid $0.3$~m}}\quad\mbox{%
  \raisebox{2pt}{\includegraphics{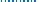}}\;{\footnotesize HVEM, Delaunay dual $1.2$~m}}\quad\mbox{%
  \raisebox{2pt}{\includegraphics{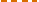}}\;{\footnotesize Abaqus, CPE4P, $400\times200$}}\quad}\\[1mm]
\includegraphics{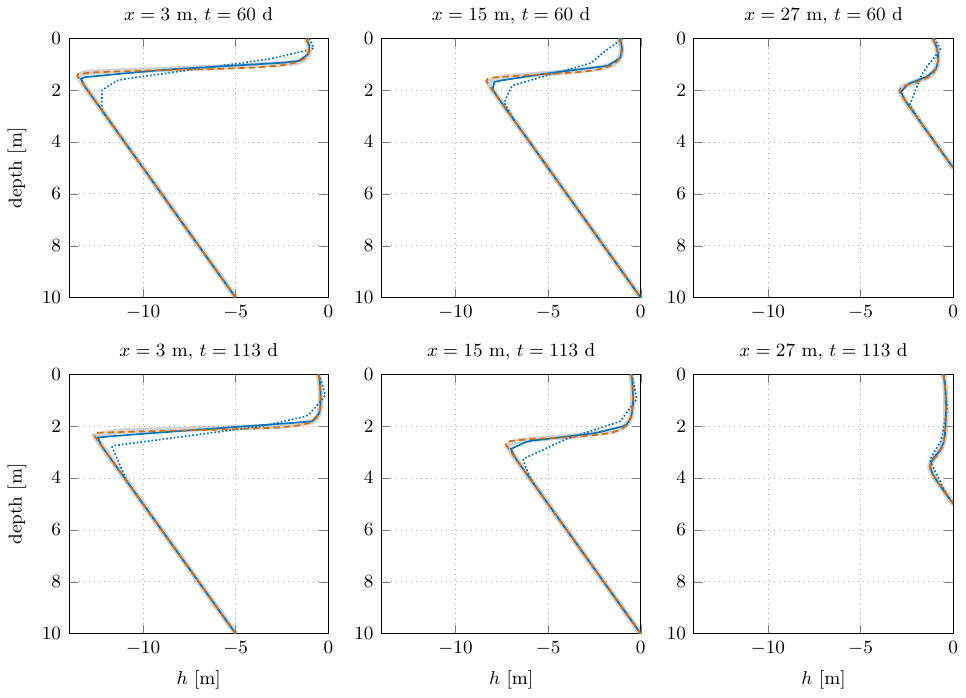}
\caption{Rainfall infiltration in a homogeneous slope: pressure head along three verticals, at the crest ($x=3$~m), on the face ($x=15$~m) and on the platform ($x=27$~m), computed by HVEM on two meshes, by the reference and by Abaqus on a boundary-fitted mesh of $400\times200$ CPE4P elements.}\label{luabq.f}
\end{figure}

Figure~\ref{lumaps.f} shows the pressure head computed by HVEM on the cut grid of $0.6$~m and on the dual mesh, after $45$ days and at the end of the record. The rain wets a layer parallel to the face, which is a few metres thick at the end of the record. Figure~\ref{luabq.f} compares the head along three verticals, at the crest, on the face and on the platform, with the independent solution of Abaqus. The reference and HVEM on the finest cut grid place the front close to Abaqus and follow its solution above and below the front.

The convergence on the cut grids is shown in Fig.~\ref{luconv.f}, and Table~\ref{lu.tab} collects the results on all the meshes. All the schemes complete the record and approach the reference at similar rates. They differ, however, in conservation, monotonicity and cost. HVEM and Q4 close the water budget to round-off, whereas SF and VEM lose a few percent of the inflow. Ahead of the front the head of VEM falls well below its initial value on the cut grids and on the quadrilaterals. On the dual mesh all the virtual elements show a fall of this kind at a few nodes, HVEM less than the others. This is the local loss of monotonicity anticipated in \S\ref{sec:cond} for polygons on which $\mathbf A^0_e$ has positive off-diagonal entries, and it does not affect the wetting of the slope. As regards the cost, HVEM takes about as many iterations as the reference on every mesh, fewer than SF and VEM, and is several times faster than them. MHFE, with its element heads averaged to the nodes, is farther from the reference than the other schemes.

\subsection{Rainfall infiltration in a layered slope}\label{sez.s}

\begin{figure}[t]
\centering
\subfloat[geometry]{%
\includegraphics{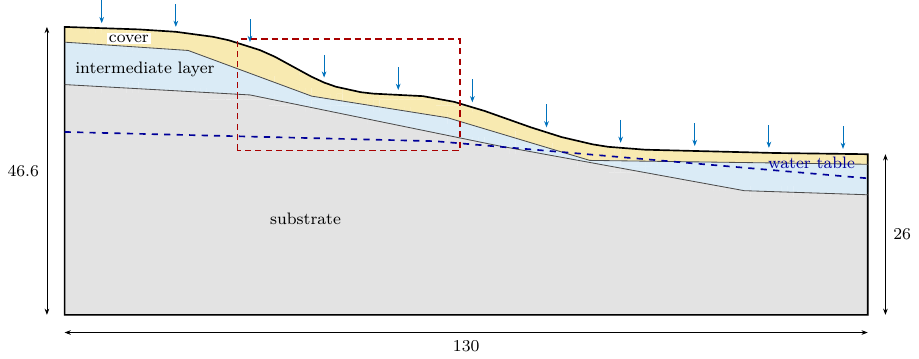}}\\[2mm]
\subfloat[cut grid, $1$~m]{\includegraphics{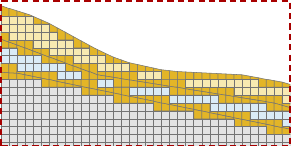}}\hfill
{%
\subfloat[refined twice in the band]{\includegraphics{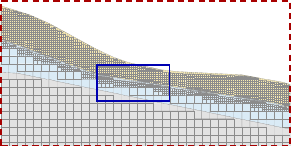}}}\hfill
\subfloat[detail of (c)]{\includegraphics{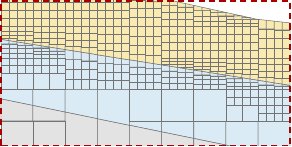}}
\caption{Rainfall infiltration in a layered slope: (a) geometry in metres, with the three layers, the initial water table and the rain; (b) the cut grid of $1$~m and (c) the same grid refined twice in the band along the surface, in the window marked in (a); (d) the detail marked in (c). The elements crossed by the surface and by the interfaces, darker, are polygons with up to eight vertices.}\label{sez.f}
\end{figure}

The last test concerns a slope of irregular profile, $130$~m long and $46.6$~m high, made of three layers separated by two inclined interfaces (Fig.~\ref{sez.f}a). A cover of thickness between $1.6$ and $4.4$~m lies on an intermediate layer, which in turn lies on a substrate. The three soils belong to the classes of Carsel and Parrish~\cite{CarselParrish1988}, and their conductivities decrease by about an order of magnitude from one layer to the next (Table~\ref{sez.soil.tab}). The water table is held on the two sides and is inclined, from $z=25$~m on the left to $17.5$~m on the right. Above it, the left side is impervious, while the right side, where the layers outcrop, is a seepage face. The initial head is hydrostatic, with the suction capped at $2.04$~m ($20$~kPa). A rain of $5$~mm/h, equal to half the saturated conductivity of the cover and slightly higher than that of the intermediate layer, falls on the whole surface for $72$~h, and the analysis is continued up to $96$~h. Where the surface saturates, the excess water runs off through a seepage condition, as on the right side.

In this test only HVEM is used, since the other schemes have been compared in the previous ones. The section is discretised by cut grids of $2$, $1$ and $0.5$~m (Fig.~\ref{sez.f}b) and by the cut grid of $1$~m refined once, twice and three times in a band formed by the cover and the first metre below it, where the wetting front and the perched saturated zones develop (Fig.~\ref{sez.f}c,d). Each element of the band is split by the horizontal and vertical lines through the centre of its bounding box. A square is thus divided into the four elements of a quadtree, and a cut element into the cut elements of the finer grid. The unrefined elements at the edge of the band keep the nodes of their refined neighbours as hanging nodes. The cut grids follow from the geometry alone, by intersecting a regular grid with the profile and the interfaces, and the polygons so produced are used as any other element. Pieces smaller than a fifth of a cell are merged with a neighbour of the same material, a simple form of the agglomeration used for virtual elements on cut meshes \cite{RangarajanSukumar2026}. On a section of this kind, with thin layers and interfaces that meet the surface at small angles, a boundary-fitted mesh would require far more care, and the ease of meshing indicates the advantage that the method can bring to realistic applications. The reference solution is computed on the grid refined three times, which has elements of $0.125$~m in the band and $44228$ nodes. Abaqus is run on the same mesh, using CPE4P quadrilaterals and splitting the other polygons into CPE6MP triangles, whose pore pressure is kept linear on every edge.

\begin{table}[t]
\footnotesize
\centering
\caption{Rainfall infiltration in a layered slope: van Genuchten--Mualem parameters of the three layers, classes of Carsel and Parrish~\cite{CarselParrish1988}, with $S_s=10^{-4}$~m$^{-1}$.}\label{sez.soil.tab}
\begin{tabular}{l l c c c c c}
\toprule
layer & class & $K_s$ [m/s] & $\theta_s$ & $\theta_r$ & $-1/h_\alpha$ [m$^{-1}$] & $n_{\mathrm{vG}}$\\
\midrule
cover        & loam            & $2.89\cdot10^{-6}$ & $0.43$ & $0.078$ & $3.6$ & $1.56$\\
intermediate & silt loam       & $1.25\cdot10^{-6}$ & $0.45$ & $0.067$ & $2.0$ & $1.41$\\
substrate    & silty clay loam & $1.94\cdot10^{-7}$ & $0.43$ & $0.089$ & $1.0$ & $1.23$\\
\bottomrule
\end{tabular}
\end{table}

\begin{figure}[t]
\centering
\subfloat[cut grid of $1$~m, $t=72$~h]{\includegraphics[width=0.48\linewidth]{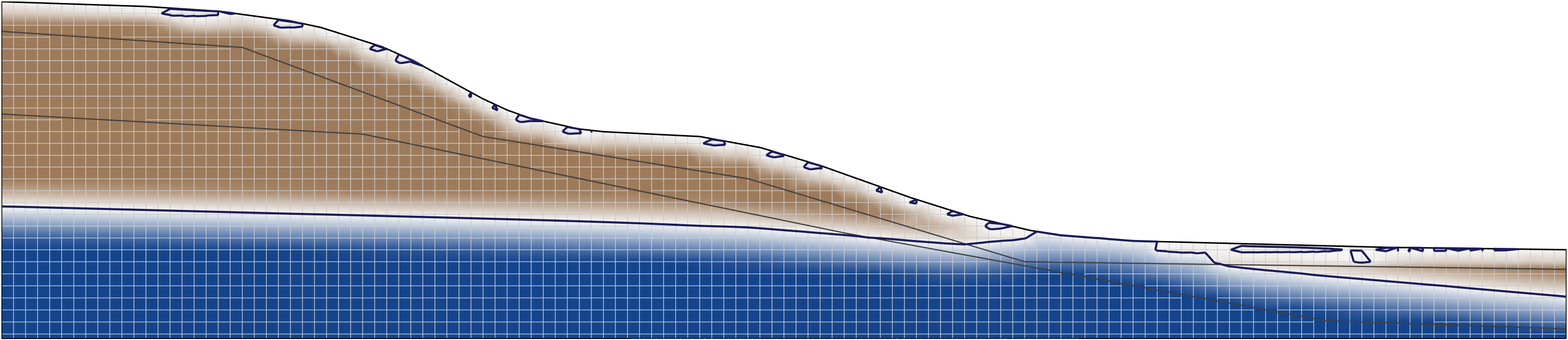}}\hfill
\subfloat[cut grid of $1$~m, $t=96$~h]{\includegraphics[width=0.48\linewidth]{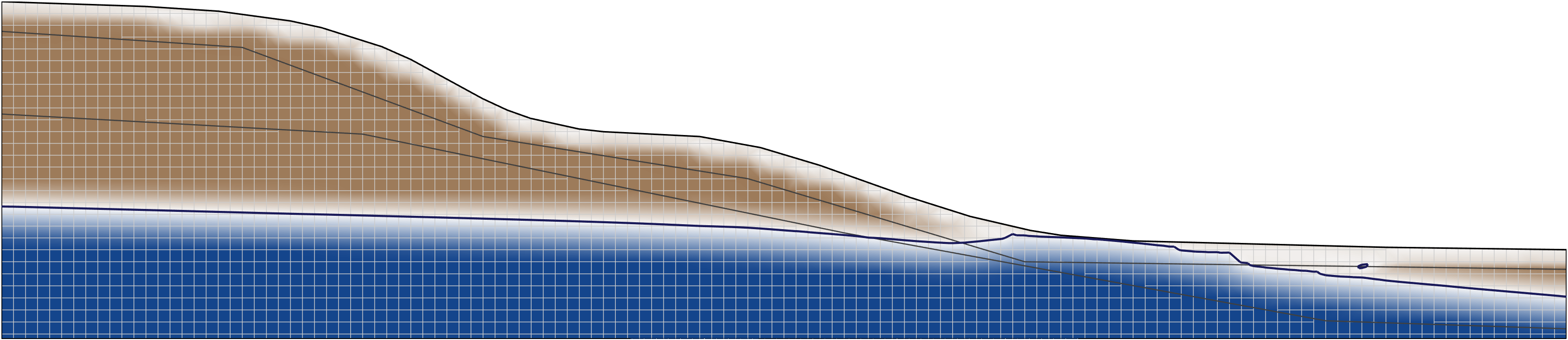}}\\[1mm]
\subfloat[band refined twice, $t=72$~h]{\includegraphics[width=0.48\linewidth]{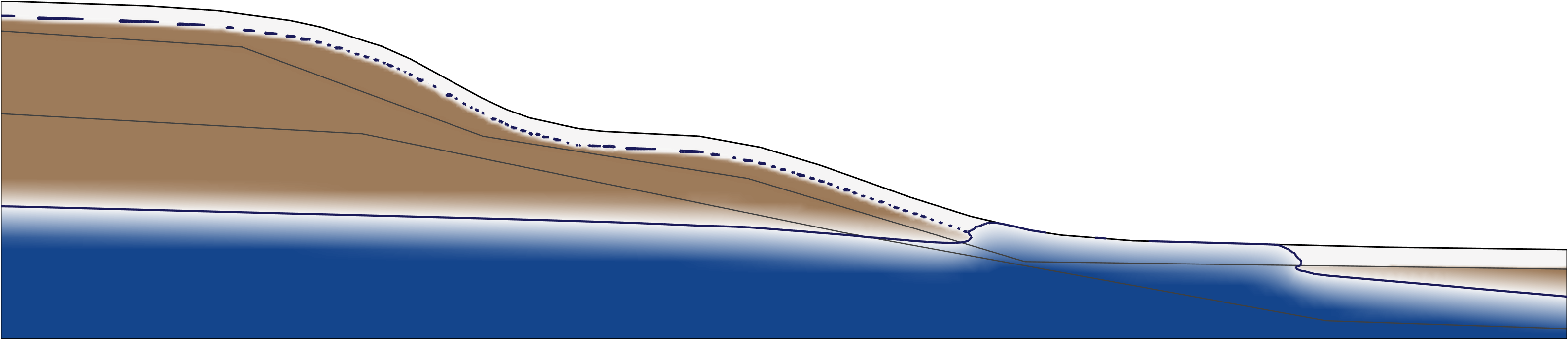}}\hfill
\subfloat[band refined twice, $t=96$~h]{\includegraphics[width=0.48\linewidth]{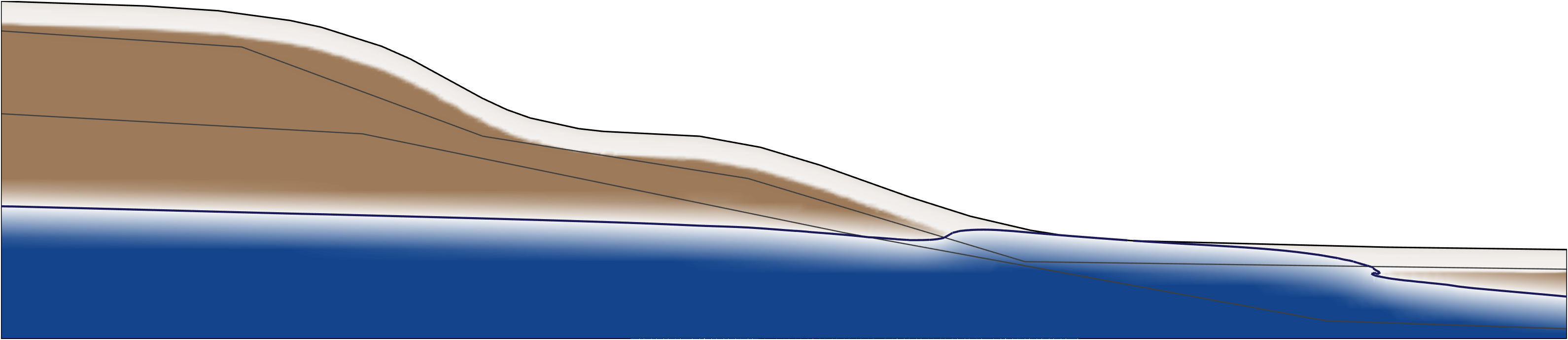}}\\[1mm]
  \includegraphics{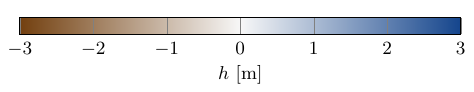}
\caption{Rainfall infiltration in a layered slope: pressure head computed by HVEM on the cut grid of $1$~m and on the same grid refined twice in the band, at the end of the rain and at the end of the analysis, above $z=14$~m, with the interfaces and the line $h=0$; the scale is limited to $\pm3$~m.}\label{sezmaps.f}
\end{figure}

\begin{table}[t]
\footnotesize
\centering
\caption{Rainfall infiltration in a layered slope: cost and water budget of HVEM on the uniform cut grids and on the grid of $1$~m refined in the band, and relative difference of the head from the reference above the initial water table at the four output times. The imbalance is in m$^2$.}\label{sez.tab}
\setlength{\tabcolsep}{4pt}
\begin{tabular}{l r r r r c cccc}
\toprule
 & & & & & & \multicolumn{4}{c}{$\mathcal E(h)$ [\%] at $t$ [h]}\\
\cmidrule(l{2pt}r{2pt}){7-10}
mesh & nodes & steps & iterations & CPU [s] & imbalance & $24$ & $48$ & $72$ & $96$\\
\midrule
cut grid, 2 m & 1397 & 193 & 858 & 11 & $7\cdot10^{-14}$ & $16.6$ & $18.5$ & $21.1$ & $21.8$\\
cut grid, 1 m & 5047 & 624 & 3096 & 185 & $-1\cdot10^{-13}$ & $11.2$ & $14.5$ & $19.0$ & $22.5$\\
cut grid, 0.5 m & 18991 & 2697 & 15975 & 2522 & $5\cdot10^{-8}$ & $7.5$ & $10.0$ & $10.8$ & $9.4$\\
1 m, band refined once & 7240 & 2781 & 16048 & 961 & $-9\cdot10^{-9}$ & $6.9$ & $9.6$ & $8.7$ & $9.4$\\
1 m, band refined twice & 15137 & 2754 & 16127 & 1543 & $3\cdot10^{-9}$ & $3.5$ & $3.7$ & $3.8$ & $3.1$\\
1 m, band refined three times & 44228 & 397 & 2255 & 395 & $-4\cdot10^{-13}$ & \multicolumn{4}{c}{reference}\\
\bottomrule
\end{tabular}
\end{table}

\begin{figure}[t]
\centering
{\mbox{%
  \raisebox{2pt}{\includegraphics{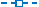}}\;{\footnotesize uniform cut grids, $2$, $1$, $0.5$~m}}\quad\mbox{%
  \raisebox{2pt}{\includegraphics{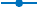}}\;{\footnotesize cut grid of $1$~m refined in the band, once and twice}}\quad}\\[1mm]
\subfloat[]{\includegraphics{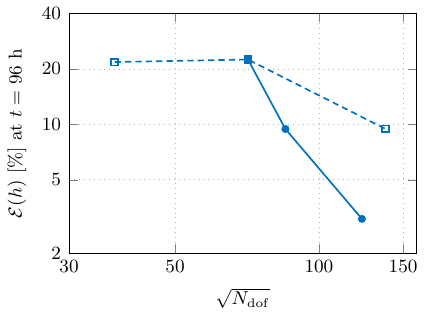}}\hfill
\subfloat[]{\includegraphics{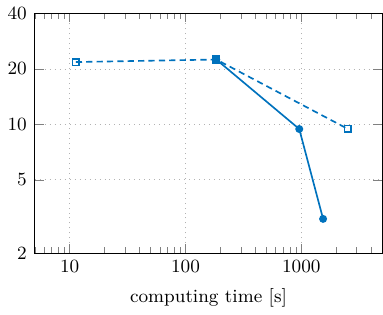}}
\caption{Rainfall infiltration in a layered slope: relative difference of the head of HVEM from the reference above the initial water table at $t=96$~h, on the uniform cut grids and on the grid of $1$~m refined in the band, against (a) $\sqrt{N_{\mathrm{dof}}}$ and (b) the computing time.}\label{sezconv.f}
\end{figure}

\begin{figure}[t]
\centering
{\mbox{%
  \raisebox{2pt}{\includegraphics{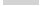}}\;{\footnotesize reference, HVEM, band refined three times}}\quad\mbox{%
  \raisebox{2pt}{\includegraphics{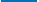}}\;{\footnotesize HVEM, band refined twice}}\quad\mbox{%
  \raisebox{2pt}{\includegraphics{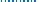}}\;{\footnotesize HVEM, cut grid $1$~m}}\quad\mbox{%
  \raisebox{2pt}{\includegraphics{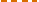}}\;{\footnotesize Abaqus, on the mesh of the reference}}\quad}\\[1mm]
\includegraphics{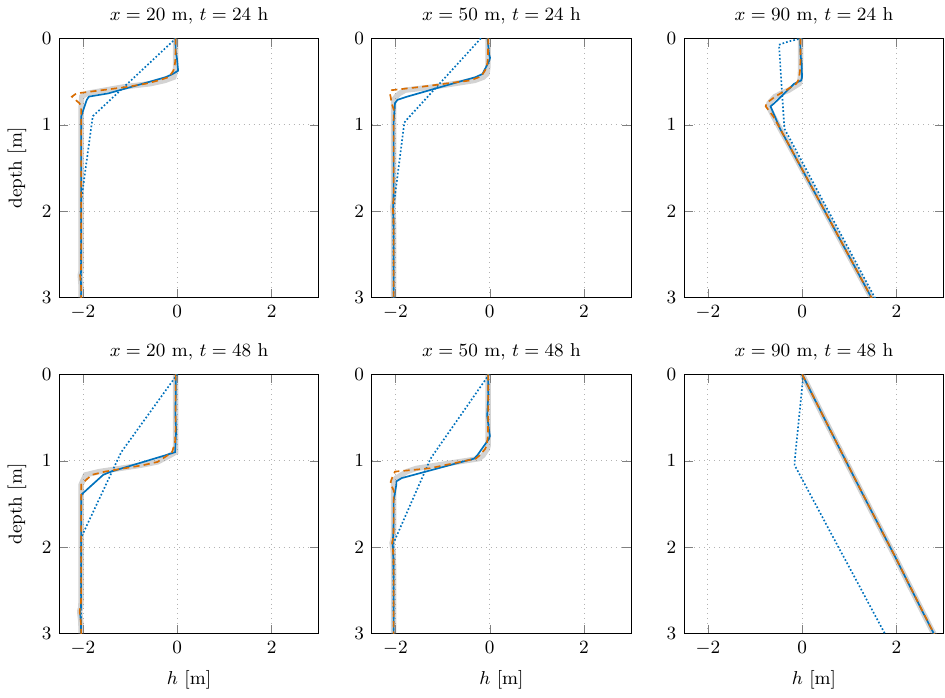}
\caption{Rainfall infiltration in a layered slope: pressure head along three verticals during the rain, computed by HVEM on the cut grid of $1$~m and on the grids refined in the band, and by Abaqus on the mesh of the reference, with CPE4P quadrilaterals and CPE6MP triangles.}\label{sezabq.f}
\end{figure}

Figure~\ref{sezmaps.f} shows the pressure head computed on the cut grid of $1$~m and on the grid refined twice, at the end of the rain and at the end of the analysis. The rain wets the cover from the surface, and in the wetted layer the head settles a few centimetres below zero, where the conductivity of the loam equals the rain rate. Where the cover is thin the front reaches its base. There the water collects on the intermediate layer, whose conductivity is lower than the rain rate, and the cover saturates up to the surface, so that part of the rain runs off. At the toe, where the substrate outcrops and the water table is shallow, the saturated zone rises to the surface. After the rain the wetted layer drains downwards and the ponded stretches dry out. The head along three verticals during the rain is compared in Fig.~\ref{sezabq.f} with the independent solution of Abaqus. The reference follows Abaqus within a few centimetres of head, and so does the grid refined twice in the band, whereas the cut grid of $1$~m smears the front and misses the saturation of the cover near the toe.

Table~\ref{sez.tab} reports the cost, the water budget and the difference from the reference above the initial water table at the four output times, and Fig.~\ref{sezconv.f} shows the difference at the end of the analysis against the number of dofs and against the computing time. Under uniform refinement the difference decreases slowly. The wetting front and the perched saturated zones are only a few decimetres thick, and the grids of $2$ and $1$~m have just one to four elements across the cover. Refining the band alone is far more effective. Refined once, the grid of $1$~m matches the accuracy of the uniform grid of $0.5$~m with about $40\%$ of its nodes and of its computing time. Refined twice, it is markedly more accurate, still with fewer nodes and less time than the uniform grid. The elements at the edge of the band, which have up to twelve vertices, are used as they are, and the water budget closes on every grid to the tolerance of the nonlinear iteration.

It is worth noting that the computing time depends on the number of steps more than on the size of the mesh. Newton's method converges slowly at the nodes approaching saturation, where the derivative of the van Genuchten--Mualem conductivity is unbounded for $n_{\mathrm{vG}}<2$, and the steps that fail are halved. The reference, whose band resolves the perched saturated zones over several elements, needs far fewer steps than the uniform grid of $0.5$~m and the grids refined once and twice, which resolve them only in part. Remedies such as the linearisations compared in \cite{ListRadu2016} or a change of primary variable \cite{BrennerCances2017} apply to HVEM as to any other element, and are not pursued here.

\clearpage
\section{Conclusions}\label{sec:concl}

A hybrid virtual element has been formulated for Richards' equation. The nodal heads are the only global dofs, and a Darcy flux is assumed within each polygon and condensed, with no stabilisation. The flux combines a divergence-free part with a particular field of linear divergence that absorbs storage and sources, so that the water balance of each element holds identically, at every iterate and for any time step. Darcy's law is written in each element for the Kirchhoff potential, and the gravity term uses the Kirchhoff mean of the conductivity between each pair of nodes, which makes hydrostatic equilibrium and gravity drainage exact on any mesh. The storage is lumped on the nodes with weights supplied by the assumed flux. The divergence-free fields are harmonic polynomials ordered by rotational index, which gives a closed-form rule for the number of flux parameters, derived on regular polygons and verified numerically on distorted polygons and on elements with hanging nodes, and every element operator reduces to a boundary integral computed once.

The numerical tests show that the element has no spurious modes on distorted polygons and on elements with hanging nodes, and that it converges at second order in the head and, in the flux, at second order on quadrilaterals and at a rate between first and second order on polygons. It is as accurate as the stabilisation-free virtual element where gravity is moderate, and it keeps its accuracy where gravity dominates, whereas the other virtual elements lose it. The water budget closes in every test, to round-off or to the tolerance of the iteration. The nonlinear iterations need only nodal values, so that the assembly costs about five times less than that of the other virtual elements, and the analyses are two to five times faster. Polygonal, cut and locally refined meshes are used with no change to the element, and local refinement reaches a given accuracy with a fraction of the dofs and of the computing time of uniform refinement.

The element shares its structure and its meshes with the hybrid virtual elements for solids. A coupled hydro-mechanical model built on the two formulations is being developed for the analysis of slopes under rainfall.

\appendix
\setcounter{table}{0}

\section{Flux fields of the polygonal elements}\label{app:fields}

This appendix lists the fields of the assumed flux of Eq.~\eqref{eq:flux} used on each polygon. The particular fields of Eq.~\eqref{eq:pfield} are the same on every element. In components along $x$ and $z$, with $\bX=(X,Z)$ the position relative to the centroid, they read
\begin{equation}\label{eq:pexplicit}
\mathbf p=\frac12\begin{bmatrix}X\\ Z\end{bmatrix},\qquad
\mathbf p_1=\frac18\begin{bmatrix}3X^2+Z^2\\ 2XZ\end{bmatrix},\qquad
\mathbf p_2=\frac18\begin{bmatrix}2XZ\\ X^2+3Z^2\end{bmatrix},
\end{equation}
and their divergences are $1$, $X$ and $Z$.

The divergence-free fields depend on the number of edges. They are written in the scaled coordinates $\boldsymbol\xi=\bX/\sqrt{\Ae}$ through the polynomials $a_m=r^m\cos m\vartheta$ and $b_m=r^m\sin m\vartheta$, which follow from $a_0=1$ and $b_0=0$ by the recurrence
\begin{equation}\label{eq:recur}
a_{m+1}=\xi_1\,a_m-\xi_2\,b_m,\qquad b_{m+1}=\xi_2\,a_m+\xi_1\,b_m .
\end{equation}
The block $\mathbf N_m$ of Eq.~\eqref{eq:Nm} holds the two fields
\begin{equation}\label{eq:blockfields}
\mathbf N_m=\begin{bmatrix}a_m&b_m\\ -b_m&a_m\end{bmatrix},
\end{equation}
which are the gradients of the harmonic polynomials $a_{m+1}/(m+1)$ and $b_{m+1}/(m+1)$ of Eq.~\eqref{eq:psi}. A polygon with $n$ edges, no two consecutive edges being aligned, uses the blocks $\mathbf N_0,\dots,\mathbf N_{\lfloor n/2\rfloor-1}$, so that polygons with $2j$ and $2j+1$ edges share the same fields (Table~\ref{fields.tab}). The flux parameters $\br$ multiply the fields in the order of the table, block after block.

\begin{table}[h]
\footnotesize
\centering
\caption{Divergence-free fields of the assumed flux on a polygon with $n$ edges, no two consecutive edges being aligned: number of flux parameters, blocks of Eq.~\eqref{eq:Nm}, and the two fields that each block adds to those of the rows above, in components along $\xi_1$ and $\xi_2$.}\label{fields.tab}
\setlength{\tabcolsep}{4pt}
\begin{tabular}{c c l c c}
\toprule
$n$ & $n_r$ & blocks & first field & second field\\
\midrule
$3$ & $2$ & $\mathbf N_0$ & $(1,\ 0)$ & $(0,\ 1)$\\[2pt]
$4$, $5$ & $4$ & $\mathbf N_0$, $\mathbf N_1$ & $(\xi_1,\ -\xi_2)$ & $(\xi_2,\ \xi_1)$\\[2pt]
$6$, $7$ & $6$ & $\mathbf N_0,\dots,\mathbf N_2$ & $(\xi_1^2-\xi_2^2,\ -2\xi_1\xi_2)$ & $(2\xi_1\xi_2,\ \xi_1^2-\xi_2^2)$\\[2pt]
$8$, $9$ & $8$ & $\mathbf N_0,\dots,\mathbf N_3$ & $(\xi_1^3-3\xi_1\xi_2^2,\ \xi_2^3-3\xi_1^2\xi_2)$ & $(3\xi_1^2\xi_2-\xi_2^3,\ \xi_1^3-3\xi_1\xi_2^2)$\\[2pt]
$10$, $11$ & $10$ & $\mathbf N_0,\dots,\mathbf N_4$ & $(\xi_1^4-6\xi_1^2\xi_2^2+\xi_2^4,\ 4\xi_1\xi_2^3-4\xi_1^3\xi_2)$ & $(4\xi_1^3\xi_2-4\xi_1\xi_2^3,\ \xi_1^4-6\xi_1^2\xi_2^2+\xi_2^4)$\\[2pt]
$2j$, $2j+1$ & $2j$ & $\mathbf N_0,\dots,\mathbf N_{j-1}$ & $(a_{j-1},\ -b_{j-1})$ & $(b_{j-1},\ a_{j-1})$\\
\bottomrule
\end{tabular}
\end{table}

Where a side of the polygon is made of collinear sub-edges, the list continues up to the block $\mathbf N_{m_{\max}}$ of Eq.~\eqref{eq:hangrule}. On the square of Fig.~\ref{hanggeom.f} with one side split into $\nu$ sub-edges, for instance, $n=\nu+3$: up to $\nu=5$ the fields are those of the table, whereas for $\nu=6$, $7$ and $8$ the last block is $\mathbf N_4$, $\mathbf N_5$ and $\mathbf N_6$, with $n_r=10$, $12$ and $14$ in place of $8$, $10$ and $10$ (Table~\ref{stability.tab}b).

\section*{Acknowledgements}
FSL acknowledges financial support from the TOWARD project (Teaming up On knoWledge-based multirisk and climate Adaptation and mitigation, eneRgy transition and smart Development; Action 1.1.2), under the Italian National Programme Research, Innovation and Competitiveness for the Green and Digital Transition 2021--2027.

\section*{Declaration of competing interest}
The author declares no known competing financial interests or personal relationships that could have appeared to influence the work reported in this paper.

\section*{Data availability}
The data of the numerical tests of Section~\ref{numerical.s}, namely the complete definition of each problem, the geometries, the meshes, the rainfall record, the reference solutions, the nodal heads computed by the element, the data behind the tables and figures and the input files of the independent Abaqus analyses, are openly available at \url{https://github.com/francescoliguori/hvem-variably-saturated-flow} under the CC~BY~4.0 licence.

\section*{Declaration of generative AI and AI-assisted technologies in the writing process}
During the preparation of this work the author used Claude (Anthropic) in order to improve the language and the organisation of the text and to assist in the preparation of figures and tables. After using this tool, the author reviewed and edited the content as needed and takes full responsibility for the content of the publication.

\clearpage
\bibliography{bibliography}

\end{document}